\documentclass[12pt,reqno]{amsart}
\usepackage[margin=1in]{geometry}
\usepackage{amsmath,amssymb,amsthm,graphicx,amsxtra, setspace}
\usepackage[utf8]{inputenc}
\usepackage{mathrsfs}
\usepackage{hyperref}
\usepackage{xcolor}
\usepackage{upgreek}
\usepackage{mathtools}
\usepackage{xcolor}
\usepackage{cite}
\allowdisplaybreaks

\usepackage[pagewise]{lineno}

\newtheorem{theorem}{Theorem}[section]

\newtheorem{lemma}[theorem]{Lemma}
\newtheorem{proposition}[theorem]{Proposition}
\newtheorem{conditions}[theorem]{Conditions}
\newtheorem{corollary}[theorem]{Corollary}
\newtheorem{definition}[theorem]{Definition}

\newtheorem{remark}[theorem]{Remark}

\numberwithin{equation}{section}

\let\originalleft\left
\let\originalright\right
\renewcommand{\left}{\mathopen{}\mathclose\bgroup\originalleft}
\renewcommand{\right}{\aftergroup\egroup\originalright}

\newcommand{\Addresses}{{% additional braces for segregating \footnotesize
		\footnote{
			\noindent \textsuperscript{1}Voronezh State University, Voronezh, Russian Federation\\
				\noindent \textsuperscript{2}Independent Researcher, India\\
				\noindent \textsuperscript{3}Lavrentyev Institute of Hydrodynamics, Novosibirsk, Russian Federation\\
				\noindent \textsuperscript{4}Institute Mathematics and Mathematical Modeling, Almaty, Kazakhstan\\
			\noindent \textsuperscript{5}Al-Farabi Kazakh National University (KazNU), Almaty, Kazakhstan
			\par\nopagebreak
			
			\noindent 	\textit{e-mail:} \texttt{Stanislav Antontsev:	antontsevsn@mail.ru}
			
			\noindent 	\textit{e-mail:} \texttt{Pardeep Kumar: s.pardeep82@gmail.com}

			\noindent 	\textit{e-mail:} \texttt{Ivan Kuznetsov: i.kuznetsov@g.nsu.ru}

			\textit{e-mail:} \texttt{Serik Aitzhanov: aitzhanovserik81@gmail.com}
			
			\noindent \textsuperscript{*}Corresponding author.
			
			\textit{Key words:} Boussinesq system, inverse problems,  integral over-determination conditions, contraction mapping principle. 	
			
			Mathematics Subject Classification (2020): Primary 35R30; Secondary 35Q35, 35Q30.

}}}

\begin{document}
	
	%	\linenumbers
	
	\title[Impulsive Inverse problems for 2D Boussinesq system]{Inverse problems for 2D Boussinesq system for incompressible viscous fluids with over-determination conditions discontinuous in time
		\Addresses}
	\author[S. Antontsev, P. Kumar, I. Kuznetsov, and S. Aitzhanov]{Stanislav Antontsev\textsuperscript{1}, Pardeep Kumar\textsuperscript{2}, Ivan Kuznetsov\textsuperscript{3,4*} and Serik Aitzhanov\textsuperscript{4,5}}

	\begin{abstract}
	In the present paper, we study the inverse problem for the Boussinesq system with the integral over-determination. The main idea is that the right-hand sides $f_n$ and $g_n$ of over-determination conditions for unknown velocity and temperature depend on $t$ and approximate discontinuity at $t=0$ as $n\to\infty$. This involves that velocity $\mathbf{v}_n$ and temperature $\vartheta_n$ also have a discontinuity at $t=0$ in the limit as $n\to\infty$. In order to pass to the limit and calculate new initial data for velocity and temperature, we apply scaling $\overline{t}=nt\in[0,1]$ and calculate final traces of the rescaled solution on an infinitesimal initial layer $[0,1]$.
	\end{abstract}
	
	\maketitle
	
	%\tableofcontents

	\section{\bf Introduction}\label{sec1}\setcounter{equation}{0}

			In the present paper, we proceed with the research initiated in \cite{Ant-2026}, where the inverse problem with the integral over-determination condition discontinuous in time for the Navier--Stokes equations (NSE) with nonlinear absorption was studied. Namely, the right-hand side of the over-determination condition for the velocity was approximated by the sequence of continuous functions. The limit led to the discontinuity of the velocity at the initial time. In the present paper, we apply the results obtained for the Boussinesq system.
			
			The Boussinesq approximation is applied to buoyancy-driven flow, or natural convection, in fluid dynamics. Density variations are ignored unless they are expressed in terms that are multiplied by the acceleration caused by gravity. The fundamental idea behind the Boussinesq approximation is that while there is very little variation in inertia between the two fluids, gravity is powerful enough to cause a noticeable difference in specific weight. Boussinesq flows are frequently found in the built environment (natural ventilation, central heating), in industry (dense gas dispersion, fume cupboard ventilation), and in nature (atmospheric fronts, oceanic circulation, katabatic winds). The approximation can be used to reduce the equations that describe such flows while maintaining a high level of accuracy in characterizing the flow behavior.
			
			The Boussinesq approximation is used in situations when the fluid's temperature fluctuates, causing a fluid flow and either mass or heat transfer. Mass, momentum, and energy conservation are all satisfied by the fluid. Along with notorious monograph \cite{boussinesq1897theorie}, this subject  was broadly described in a recent monograph \cite{tritton1988physical}. 
			From a mathematical point of view, the direct problem for the Boussinesq system  was broadly studied in a great variety of papers. Some of them were mentioned in \cite{feireisl2009singular}.    
			
			We study the following time-dependent inverse source problem for the heat convection system governing an incompressible isothermal viscous Navier--Stokes flow
			\begin{subequations}
				\label{eq:1o}
				\begin{align}
					&\label{eq:1ao} \partial_{t}\mathbf{v}+\left(\mathbf{v}\cdot \nabla \right) \mathbf{v}+\nabla\pi=\nu \Delta\mathbf{v}+\mathbf{g}(x,t)\vartheta
					+\varphi(t)\boldsymbol{\omega}(x),\ \   \left(x,t\right)\in Q_T,
					\\
					&\label{eq:divo}
					\operatorname{div}\mathbf{v}(x,t)=0,\ \  \left(x,t\right)\in Q_T,
					\\
					&\label{eq:heato}
					\partial_{t}\vartheta+\left(\mathbf{v}\cdot \nabla \right) \vartheta= \mu \Delta \vartheta+\psi(t)\eta(x),\ \ \left(x,t\right)\in Q_T,
				\end{align}
				which is supplemented with the initial conditions
				\begin{align}
					\mathbf{v}(x,0) =\mathbf{v}_{0}
					\left(x\right),\quad \vartheta(x,0)=
					\vartheta_0(x), \ \ x\in \Omega, \label{eq:2o}
				\end{align}
				the boundary conditions
				\begin{align}
					\mathbf{v}(x,t)=\mathbf{0}, \ \ \vartheta(x,t)=0,  \  \ \ \ 
					\left(x, t\right)\in\Gamma _{T}, \label{eq:3o}
				\end{align}
				and with the integral over-determination conditions
				\begin{equation}
					\label{eq:4o}  \int_{\Omega}
					\mathbf{v}(x,t)\cdot\boldsymbol{\omega}(x)\,dx=f(t),
					\ \ \ \ 
					\int_{\Omega} \vartheta(x,t)\eta(x)\,dx =g(t),
					\quad
					t\geq0.
				\end{equation}
			\end{subequations}
			Here, $\Omega$ is a bounded domain in $\mathbb{R}^2$ with a smooth boundary $\partial\Omega$, and $Q_T=\Omega\times \left(0,T\right) $ is a cylinder with lateral boundary $\Gamma_T=\partial\Omega\times(0,T)$, the vector-valued function $\mathbf{v}(x,t)$ is a velocity field, and the scalar-valued functions $\pi(x,t)$ and $\vartheta(x,t)$ are pressure and temperature, respectively. The positive constant $\nu$ denotes the kinematic viscosity, and the positive constant $\mu>0$ denotes the heat conductivity of fluids. In \eqref{eq:1o}, $\mathbf{v}_{0}(x)$, $\vartheta_{0}(x)$, $\mathbf{g}(x,t)$, $f(t)$, $g(t)$, $\boldsymbol{\omega}(x)$, and $\eta(x)$ are given functions. In particular, the functions  
			\[
			\mathbf{F}(x,t):=\varphi(t)\boldsymbol{\omega}(x) \ \ \text{ and } \ \ {G}(x,t):=
			\psi(t)\eta(x)
			\]
			denote the densities of the external force and the heat source with intensities $\varphi(t)$ and $\psi(t)$, respectively.

			It is known that the right-hand sides of the over-determination conditions should be continuous or belong to appropriate Sobolev spaces. The main novelty of the present paper is that we consider the right-hand sides of the over-determination conditions to be discontinuous at $t=0$. Namely, we approximate these discontinuous functions by sequences of continuous functions and pass to the limit.
			
			This work addresses a time-dependent impulsive inverse source problem to determine the functions $(\mathbf{v}, \vartheta,\varphi,\psi)$ that satisfy the system \eqref{eq:1o}, where the constants  $\nu$ and $\mu$, as well as the  functions $\mathbf{v}_{0}(x)$, $\vartheta_{0}(x)$, $\mathbf{g}(x,t)$, $f(t)$, $g(t)$, $\boldsymbol{\omega}(x)$, and $\eta(x)$ are   given.
			
			For a smooth right-hand side of the integral over-determination condition, the time-dependent inverse source problem for fluid flow was studied in \cite{Prilepko1989,Prilepko2000,Vasin1993,Kumar2022}, and the references therein. In \cite{Prilepko1989}, the authors proved the unique solvability result for an inverse problem for the linearized NSE with both integral and final over-determination data using the Banach fixed point theorem.
			The existence and uniqueness of generalized solutions to the time-dependent inverse source problem for the nonlinear Navier--Stokes system with integral over-determination were studied in \cite{Vasin1993}, see also \cite[Section 4.6]{Prilepko2000}. Recently, the well-posedness of an inverse problem for 2D convective Brinkman--Forchheimer (CBF) equations (or NSE with absorption terms) with integral over-determination data was established in \cite{Kumar2022} using an extension of the contraction mapping theorem. For a comprehensive introduction to inverse problems with over-determination conditions, we refer the reader to the monograph \cite{Prilepko2000} and references \cite{Cheb2011,FNN2009,AbyAit2015,ant-kho:2021}.
			
			Inverse problems for these systems have not been studied extensively. For instance, in \cite{AbyAit2015,Cheb2011,FNN2009}, the existence and uniqueness of solutions were established in the two-dimensional case. In \cite{FNN2009}, the authors proved the Lipschitz stability for the inverse source problem of the Boussinesq system in two dimensions. Here, they identified a time-independent external force in the system using observation data of velocity over an arbitrary subdomain and over a time interval, together with velocity and temperature data at a fixed positive time $t_0>0$ over the whole spatial domain. Also, recent results \cite{AK:2023,antontsev2024inverse} on viscoelastic fluid are of great importance.
			
			This problem \eqref{eq:1o} can be formulated in the equivalent way as the direct problem
			\[
			\varphi(t;[\mathbf{v},\vartheta])=\frac{1}{\|\boldsymbol\omega\|_{\mathbf{H}}^2}\Big(f'(t)+((\mathbf{v}\cdot\nabla)\mathbf{v},\boldsymbol{\omega})+\nu(\nabla \mathbf{v},
			\nabla\boldsymbol{\omega})-(\mathbf{g}\vartheta,\boldsymbol{\omega})\Big),
			\]
			\[
			\psi(t;[\mathbf{v},\vartheta])=\frac{1}{\|\eta\|_{L^2(\Omega)}^2}
			\Big(g'(t)+((\mathbf v\cdot\nabla)\vartheta,\eta)+\mu(\nabla\vartheta,\nabla\eta)\Big).
			\]
			
			The main idea is to approximate $f$ and $g$, which are discontinuous functions, with the help of the sequences $\{f_n\}_{n\geq n_0}$ and $\{g_n\}_{n\geq n_0}$. This means that $\{f_n'\}_{n\geq n_0}$ and  $\{g_n'\}_{n\geq n_0}$ approximate functionals. Automatically, the presence of such singular terms leads to the discontinuous solutions at $t=0$. These direct problems are linked with impulsive differential equations in the sense of J. Kurzweil \cite{Kurzweil1958}. The main aspect of J. Kurzweil's approach is to study the sequences of solutions of ODE, where the right-hand side is not an integrable Lebesgue function.
			
			The first example is the ODE:
			\[
			{x}_n'(t)=\theta_n'(t)x_n(t),\quad x_n(0)=x_0,
			\]
			where $\theta_n(t)$ approximates the Heaviside function $\theta(t)=\mathbf{1}_{t>0}$. 
			In the limit as $n\to\infty$, we get the condition on the discontinuity:
			\[
			x(0^+)=x_0e.
			\]
			The second example is the system of non-integrable equations. In order to formulate the condition on the discontinuity, we need to apply scaling $\overline{t}=nt\in[0,1]$ and formulate the problem on the infinitesimal initial layer. In the limit, the final trace of the scaled solution serves as the right-sided trace of a limit discontinuous solution.
			
			{Impulsive differential equations have become an active area of research due to their ability to model systems that undergo sudden changes at certain moments. They arise naturally in many applications where instantaneous effects or abrupt changes occur in the state of a system. The impulsive direct problems for fluid dynamic models were studied in \cite{AKS,AKS-2024}, impulsive pseudoparabolic equations in \cite{KS-NARWA-2022,KS-NA-TMA-2023}, semilinear impulsive pseudoparabolic equations in \cite{AKSS-2024}, and impulsive $p(x,t)$-parabolic equations in \cite{AKSS-2024-2}.
				
				In this paper, we employ a weak regular solutions approach to formulate and solve the time-dependent impulsive inverse source problem for the functions $(\mathbf{v}, \vartheta, \varphi,\psi)$. Following the approach developed in \cite{Ant-2026} for the impulsive inverse problem for the CBF equations in two dimensions, we establish the existence and uniqueness of weak regular solutions to the impulsive inverse source problem for the 2D Boussinesq system.}
			
			This paper is organized as follows. In Section \ref{Sect:Prelim}, we introduce the notation, several important facts, and the relevant assumptions. In Section \ref{sec3}, we formulate the main inverse problem \eqref{eq:1} and establish its relation to the corresponding functional direct problem. In the same section, we prove the existence and uniqueness result stated in Theorem \ref{thm1}. The main novelty is the estimates, which are uniform in $n\in\mathbb{N}$. In Section \ref{sec4}, we pass to the limit as $n\to\infty$. Theorem \ref{theorem3.1} presents the main result of the paper. Auxiliary procedures, such as scaling and shifting, are employed. 
			As $n\to\infty$, we deduce the linear direct problem in the infinitesimal initial layer, which gives the new initial data $\mathbf{v}_0({x})+(\overline{f}(1)-f(0))\boldsymbol{\omega}({x}),$ and $\vartheta_0({x})+(\overline{g}(1)-g(0))\eta({x}) $. 
			
\section{\bf Preliminaries}\label{Sect:Prelim}\setcounter{equation}{0}
We begin this section by introducing the preliminaries that will be used throughout the paper. We then provide the necessary conditions on the input data required to establish our main results.

\subsection{\bf Function spaces and preliminary results}

In this subsection, we introduce the main functional spaces and some useful inequalities.
We distinguish vectors from scalars by using boldface letters. For functions and function spaces, we will use this distinction as well. The symbol $C$ will denote a generic constant-generally a positive one whose value will not be specified; it can change from one inequality to another.

We denote by $\mathbf{L}^{2}(\Omega)$ the usual Lebesgue space of square integrable vector-valued functions on $\Omega$, and by $\mathbf{W}^{l,2}(\Omega)$ the Sobolev space of functions in $\mathbf{L}^{2}(\Omega)$ whose weak derivatives of order not greater than $l$ are in $\mathbf{L}^{2}(\Omega)$.

Let us introduce the following function spaces
widely used in Mathematical Fluid Mechanics, see \cite{Temam}:
\begin{align*}
	&\mathcal{V}:=\{\mathbf{v} \in\mathbf{C}_{0}^{\infty}(\Omega):\operatorname{div}\mathbf{v}=0\},
	\\
	&\mathbf{H}:=\mbox{closure of }\mathcal{V}\mbox{ in the norm of }\mathbf{L}^{2}(\Omega), 
	\\
	&\mathbf{V}^l:=\mbox{closure of }\mathcal{V}\mbox{ in the norm of }\mathbf{W}_{0}^{l,2}(\Omega),\, \ l\in\mathbb{N},
	\\
	&\widetilde{\mathbf{L}}^p:=\mbox{closure of }\mathcal{V}\mbox{ in the norm of }\mathbf{L}^{p}(\Omega), \ \ p\in (2,\infty].
\end{align*}
As usual, when $l=1$, we simply write $\mathbf{V}$ instead of $\mathbf{V}^{1}$. Let $(\cdot,\cdot)$ denote the inner product in the Hilbert space $\mathbf{H}$ (or equivalently in $L^2(\Omega)$) and $\langle\cdot,\cdot\rangle$ denotes the duality pairing between the spaces $\mathbf{V}$ (or $W_0^{1,2}(\Omega)$) and its dual $\mathbf{V}^{\prime}$ (or $W^{-1,2}(\Omega)$). 

Integration by parts gives 
\begin{align*}
	\big((\mathbf{u}\cdot\nabla)\mathbf{v},\mathbf{w}\big) =  -\big((\mathbf{u}\cdot\nabla)\mathbf{w},\mathbf{v}\big), \ \text{ for all } \mathbf{u},\mathbf{v},\mathbf{w}\in \mathbf{V},
\end{align*}
and 
\begin{align*}
	\big((\mathbf{u}\cdot\nabla)\mathbf{v},\mathbf{v}\big) = 0, \ \text{ for all } \mathbf{u},\mathbf{v} \in\mathbf{V}.
\end{align*}
Similarly, we have 
\begin{align*}
	\big((\mathbf{u}\cdot\nabla)\vartheta,\vartheta\big) = 0, \ \text{ for all } \mathbf{u} \in\mathbf{V}, \ \vartheta \in W_0^{1,2}(\Omega).
\end{align*}

We will use the following embedding inequalities from \cite{M:2011}.
\begin{lemma}\label{embed:1}
	Let $\Omega\subset\mathbb{R}^2$ be a bounded domain. If the boundary $\partial\Omega$ is assumed to be of class $C^{0,1}$, then there exist positive constants $C_{Ag}$ and $C_{CS}$, depending only on $\Omega$, such that the following estimates are valid:
	\begin{enumerate}
		\item Ladyzhenskaya's inequality
		\begin{align}
			\label{SE:case:p*>4}
			&\left\Vert  \mathbf{v}\right\Vert_{\widetilde{\mathbf{L}}^4}
			\leq {2}^{\frac{1}4}\left\Vert\mathbf{v}\right\Vert_{\mathbf{H}}^{\frac{1}2}
			\left\Vert \nabla \mathbf{v}\right\Vert_{\mathbf{H}}^{\frac{1}2},\,\, \text{ for all } \mathbf{v}\in \mathbf{V},
		\end{align}
		\item Agmon's inequality  \cite[Eq. (2.23)$_1$, p. 11]{Temam1995}, 
		\[
		\|\mathbf{v}\|_{\widetilde{\mathbf{L}}^\infty}\leq C_{Ag}
		\|\mathbf{v}\|_{\mathbf{H}}^{\frac12}
		\|\mathbf{v}\|_{\mathbf{V}^2}^{\frac12},
		\]
		\item Cattabriga--Solonnikov's inequality \cite[p. 113]{Seregin2014},
		\begin{align}
			&\|\mathbf{v}\|_{\mathbf{V}^2}\leq C_{CS}\|\Delta\mathbf{v}\|_{\mathbf{H}}.
		\end{align}
	\end{enumerate}
\end{lemma}

\begin{lemma}[\bf The Aubin--Lions--Simon Lemma] (\cite[Corol.\,4 \& 5]{Simon1987})
	\label{ALSlemma}
	Assume that $\mathbf{X}$, $\mathbf{B}$, and $\mathbf{Y}$ are three Banach spaces such that $\mathbf{X}\hookrightarrow \hookrightarrow \mathbf{B} \hookrightarrow \mathbf{Y}$, i.e., $\mathbf{X}$ is compactly embedded in $\mathbf{B}$, which is continuously embedded in $\mathbf{Y}$. Then, the following assertion holds true:
	\\[1ex]
	\indent For $1\leq p<\infty$ and  $1\leq q \leq\infty$, if some set $\mathcal{F}$ is bounded in $L^p(0,T;\mathbf{X})$, and 
	\begin{itemize}
		\item the set $ \partial \mathcal{F}/\partial t \stackrel{\scriptsize {\rm def}}{=} \bigl\{\partial_t \phi : \phi\in \mathcal{F} \bigr\}$ is bounded in $L^q(0,T;\mathbf{Y})$,
	\end{itemize}
	\quad or 
	\begin{itemize}
		\item\emph{$\int_0^{T-h}\|\phi(\cdot,t+h)-\phi(\cdot,t)\|_{\mathbf{Y}}^q\,dt\le C\omega(h)$, where $C$ does not depend on $\phi$, $\lim\limits_{h\to0^+}\omega(h)=0$,}
	\end{itemize}
	then $\mathcal{F}$ is relatively compact in $L^p(0,T;\mathbf{B})$. If $q>1$ and $p=\infty$, then $\mathcal{F}$ is relatively compact in $\mathrm{C}([0,T];\mathbf{B})$.
\end{lemma}
In the sequel, the symbols $\rightarrow$, $ {\rightharpoonup}  $, and $\overset{*}{\rightharpoonup}$  denote strong, weak, and weak$^*$ convergence, respectively.
For the definitions and properties of these notions of convergence, we refer the reader to \cite[Chapter 27]{Robinson}.

\subsection{\bf Conditions on the coefficients and the input data}
In this subsection, we introduce the conditions on the coefficients and the input data that will be used throughout the paper. These conditions play a crucial role in establishing the main result and in passing to the limit as $n\to\infty$. 
\begin{conditions}\label{C} 
	We impose the following conditions, ${\bf C}_{1}$--${\bf C}_{11}$:
	\begin{itemize}
		\item[${\bf C}_{1}.$] Functions $\boldsymbol{\omega} \in\mathbf{V} \cap \mathrm{C}^1(\overline\Omega)$ and $\eta\in W_0^{1,2}(\Omega)\cap \mathrm{C}^1(\overline\Omega)$ satisfy the following conditions
		\[
		\int_\Omega\eta^2(x)\,dx=\int_{\Omega}\boldsymbol{\omega}^2(x)\,dx=1.
		\]
		\item[${\bf C}_{2}.$] The only assumption on the positive parameter $T$ is $T>\frac1{n_0}$.
		\item[${\bf C}_{3}.$] The sequences  $\{f_n(t)\}_{n\geq n_0}$ and $\{g_n(t)\}_{n\geq n_0}$  of  $W^{1,2}(0,T)$  functions approximate the discontinuous functions $f(t)$ and $g(t)$, respectively, and satisfy the conditions $f_n(0)=f(0)$ and $g_n(0)=g(0)$, $f_n(t)\equiv f(t)$ and $g_n(t)\equiv g(t)$ for all $t\in[\frac1n,T_n]$. Therefore, we have
		\[
		\mathrm{supp}(f_n-f)=\mathrm{supp}(g_n-g)=[0,1/n].
		\]
		\item[${\bf C}_{4}.$] The inequalities 
		\[
		\sup_{t\in[0,T_n]}|f_n(t)|\le F,\,\,\sup_{t\in[0,T_n]}|g_n(t)|\le G,
		\,\,
		\int_{0}^{T_n}|f_n'(t)|\,dt\le F_1,\,\,\int_{0}^{T_n}|g_n'(t)|\,dt\le G_1,
		\]
		for every $n\in\mathbb{N}$, where $F, \, G, \, F_1$, and $G_1$ do not depend on $n$.
		Moreover, we have
		\begin{align*}
			\lim\limits_{n\to\infty}\int_{0}^{1/n}|f_n'(t)|^2\,dt
			=\lim\limits_{n\to\infty}\int_{0}^{1/n}|g_n'(t)|^2\,dt=\infty.
		\end{align*}
		\item[${\bf C}_{5}.$] The sequences $\{\overline{f}_n\}_{n\geq n_0}$ and $\{\overline{g}_n\}_{n\geq n_0}$ converge strongly in $\mathrm{C}([0,1])$ to $\overline{f}$ and $\overline{g}$, respectively, where $\overline{f}_n(\overline{t})=f_n({\overline{t}}/n)$ and $\overline{g}_n(\overline{t})=g_n({\overline{t}}/n)$.
		\item[${\bf C}_{6}.$] The sequences $\{\overline{f}_n'\}_{n\geq n_0}$ and $\{\overline{g}_n'\}_{n\geq n_0}$ converge strongly in $L^2(0,1)$ to $\overline{f}'$ and $\overline{g}'$, respectively.
		\item[${\bf C}_{7}.$] Since $\widetilde{f}_n(\widetilde{t})=f_n(\widetilde{t}+\frac1n)=
		f(\widetilde{t}+\frac1n)$ and $\widetilde{g}_n(\widetilde{t})=g_n(\widetilde{t}+\frac1n)=
		g(\widetilde{t}+\frac1n)$, the inequalities $\bigl(\int_0^T|\widetilde{f}_n'(\widetilde{t})|^2\,d\widetilde{t}\bigr)^{\frac{1}{2}}\le F_2$ and $\bigl(\int_0^T|\widetilde{g}_n'(\widetilde{t})|^2\,d\widetilde{t}\bigr)^{\frac{1}{2}}\le G_2$ hold for every $n\in\mathbb{N}$, where $F_2$ and $G_2$ do not depend on $n$.
		\item[${\bf C}_{8}.$] The sequences $\{\widetilde{f}_n'\}_{n\geq n_0}$ and 
		$\{\widetilde{g}_n'\}_{n\geq n_0}$ converge strongly in $L^2(0,T)$ to $f'$ and $g'$, respectively.
		\item[${\bf C}_{9}.$] The matching conditions hold
		\[
		\int_\Omega\mathbf{v}_0(x)\cdot\boldsymbol{\omega}(x)\,dx=f(0),\quad 
		\int_\Omega\vartheta_0(x)\eta(x)\,dx=g(0).
		\]
		\item[${\bf C}_{10}. $] $\mathbf{v}_0\in\mathbf{H}$, $\vartheta_0\in L^{2}(\Omega)$
		for a weak solution of problem \eqref{eq:1}.
		\item[${\bf C}_{10}^\prime.$] $\mathbf{v}_0\in\mathbf{V}$, $\vartheta_0\in W_0^{1,2}(\Omega)$ for a weak regular solution of problem \eqref{eq:1}.
		\item[${\bf C}_{11}. $] We assume that $\mathbf{g}(x,t) \in \mathrm{C}(Q_T)$ and $\exists g_0>0:$   $\max\limits_{Q_T}|\mathbf{g}(x,t)|\leq g_0$.
	\end{itemize}
\end{conditions}

\section{\bf Inverse and equivalent direct problems under study}\label{sec3}\setcounter{equation}{0}

This section begins by formulating the main inverse problem and establishing its relation to the corresponding nonlocal direct problem.  We also prove the existence of  weak and weak regular solutions  to the nonlocal direct problem.

\subsection{\bf Inverse problem}

In this subsection, we deal with the family of inverse problems for the Boussinesq system for incompressible viscous fluids depending on $n\in\mathbb{N}$. For a fixed $n$, we formulate the
inverse problem for the Boussinesq system 

\begin{subequations}
	\label{eq:1}
	\begin{align}
		&\label{eq:1a} \partial_{t}\mathbf{v}_{n}+\left(\mathbf{v}_n\cdot \nabla \right) \mathbf{v}_n+\nabla\pi_n=\nu \Delta\mathbf{v}_n+\mathbf{g}(x,t)\vartheta_n +\varphi_n(t)\boldsymbol{\omega}(x), \ \ \ \text{ in } Q_{T_n},
		\\
		&\label{eq:div}
		\operatorname{div}\mathbf{v}_n(x,t)=0, \ \ \ \text{ in } Q_{T_n},
		\\
		&\label{eq:heat}
		\partial_{t}\vartheta_{n}+\left(\mathbf{v}_n\cdot \nabla \right) \vartheta_n= \mu \Delta \vartheta_n+\psi_n(t)\eta(x), \ \ \ \text{ in } Q_{T_n},
	\end{align}
	which is supplemented with the initial conditions
	\begin{align}
		\mathbf{v}_n(x,0) =\mathbf{v}_{0}
		\left(x\right),\quad \vartheta_n(x,0)=
		\vartheta_0(x), \ \ \ \text{ in } \Omega, \label{eq:2}
	\end{align}
	the boundary conditions
	\begin{align}
		\mathbf{v}_n(x,t)=\mathbf{0},  \quad \vartheta_n(x,t)=0,
		\ \ \ \text{ in } \Gamma_{T_n}, \label{eq:3}
	\end{align}
	and with the integral over-determination conditions
	\begin{equation}
		\label{eq:4}  \int_{\Omega}
		\mathbf{v}_n(x,t)\cdot\boldsymbol{\omega}(x)\,dx=f_n(t),
		\quad
		\int_{\Omega} \vartheta_n(x,t)\eta(x)\,dx =g_n(t),
		\quad t\geq0.
	\end{equation}
\end{subequations}

Here, the vector-valued function $\mathbf{v}_n(x,t)$ is a velocity field, and the scalar-valued functions $\pi_n(x,t)$ and $\vartheta_n(x,t)$ are pressure and temperature, respectively. Let $Q_{T_n}=\Omega\times \left(0,T_n\right) $ be a cylinder with a lateral boundary $\Gamma_{T_n}=\partial\Omega\times (0,T_n)$, $T_n := T + \frac{1}{n}$ be a positive constant. The functions $\mathbf{v}_{0}(x)$, $\vartheta_{0}(x)$, $\mathbf{g}(x,t)$, $f_n(t)$, $g_n(t)$, $\boldsymbol{\omega}(x)$, and $\eta(x) $ are given. Namely, the functions  
\[
\mathbf{F}_n(x,t):=\varphi_n(t)\boldsymbol{\omega}(x) \ \ \text{ and }   \ \ {G}_n(x,t):=\psi_n(t)\eta(x)
\]
denote the densities of the external force and the heat source with intensities $\varphi_n(t)$ and $\psi_n(t)$, respectively.

\subsection{\bf Weak formulation}\label{Weak_formul}

A {\bf weak solution} to the inverse problem (\ref{eq:1}) is understood in the following sense:
\begin{subequations}
	\begin{definition}\label{def:ws}
		The collection of functions $(\mathbf{v}_n(x,t),\vartheta_n(x,t),\varphi_n(t),\psi_n(t))$ is called {\bf a weak solution} to the inverse problem \eqref{eq:1}, if:
		\begin{enumerate}
			\item $\mathbf{v}_n\in {L}^{\infty}(0,T_n;\mathbf{H})\cap {L}^{2}(0,T_n;\mathbf{V}),$  \ \ $\partial_{t}\mathbf{v}_{n}\in {L}^{2}(0,T_n;\mathbf{V}^\prime)$,   \ \ $\varphi_n(t)\in L^2(0,T_n)$;
			\item $\vartheta_n\in{L}^{\infty}(0,T_n;L^2(\Omega))\cap{L}^{2}(0,T_n;W_0^{1,2}(\Omega))$, \ \ $\partial_{t}\vartheta_{n}\in L^2(0,T_n;W^{-1,2}(\Omega))$, \\  
			$\psi_n(t)\in L^2(0,T_n)$;
			\item initial conditions are valid in the weak sense:  
			\[
			\lim\limits_{t\to0^+}
			\int_\Omega(\mathbf{v}_n(x,t)-\mathbf{v}_{0}(x))\cdot\boldsymbol{\phi}(x)\,dx=0,\quad
			\lim\limits_{t\to0^+}
			\int_\Omega(\vartheta_n(x,t)
			-\vartheta_{0}(x))\zeta(x)\,dx=0,
			\]
			for every $\mathbf{\phi}\in\mathbf{V}$ and $\zeta\in W_0^{1,2}(\Omega)$;
			\item  integral over-determination conditions \eqref{eq:4} holds for all $t\in [0,T_n]$;
			\item For every $\boldsymbol{\phi}\in\mathbf{V}$,
			\begin{align}
				\begin{split}
					&\dfrac{d}{dt}\int_{\Omega}\mathbf{v}_n\cdot \boldsymbol{\phi} \, d x+\int_{\Omega}\left( \mathbf{v}_n\cdot {\nabla}\right)\mathbf{v}_n \cdot \boldsymbol{\phi} \, d x
					+\nu \int_{\Omega} \nabla \mathbf{v}_n : \nabla\boldsymbol{\phi} \, d x \\
					&\quad=\varphi_n(t) \int_{\Omega}\boldsymbol{\omega}(x) \cdot \boldsymbol{\phi} \, d x
					+\int_{\Omega}\vartheta_n\mathbf{g}(x,t)\cdot \boldsymbol{\phi} \, d x
				\end{split}
				\label{we:equa}
			\end{align}
			holds in the distribution sense on $(0,T_n)$;
			\item For every $\zeta\in W_0^{1,2}(\Omega)$,
			\begin{equation}
				\label{we:eqheat}% 
				\dfrac{d}{dt} \int_{\Omega}\vartheta_n\zeta \, d x + \int_{\Omega}\left(\mathbf{v}_n \cdot \nabla\right)\vartheta_n \zeta \, d x
				+\mu \int_{\Omega} {\nabla}\vartheta_n \cdot  {\nabla}\zeta \, d x
				=\psi_n(t) \int_{\Omega} \eta (x) \zeta \, d x
			\end{equation}
			holds in the distribution sense on $(0,T_n)$.
		\end{enumerate}
	\end{definition}
\end{subequations}
A {\bf weak regular solution} to the inverse problem \eqref{eq:1} is understood in the following sense:
\begin{subequations}
	\begin{definition}\label{def:wrs}
		The collection of functions $(\mathbf{v}_n(x,t),\vartheta_n(x,t),\varphi_n(t),\psi_n(t))$ is called a {\bf weak regular solution} to the inverse problem \eqref{eq:1}, if it is a {\bf weak solution} to the inverse problem \eqref{eq:1} and satisfies additional conditions:
		\begin{enumerate}
			\item $\mathbf{v}_n\in \mathrm{C}([0,T_n];\mathbf{H})\cap{L}^{\infty}(0,T_n;\mathbf{V})\cap{L}^{2}(0,T_n;\mathbf{V}^2)$;  
			\item $\vartheta_n\in \mathrm{C}([0,T_n];L^2(\Omega))\cap{L}^{\infty}(0,T_n;W_0^{1,2}(\Omega))\cap {L}^{2}(0,T_n;W^{2,2}(\Omega))$;
			\item initial conditions are valid in a strong sense 
			\[
			\lim\limits_{t\to0^+}
			\|\mathbf{v}_n(\cdot,t)-\mathbf{v}_{0}(\cdot)\|_{\mathbf{H}}=0,\quad
			\lim\limits_{t\to0^+}
			\|\vartheta_n(\cdot,t)-\vartheta_{0}(\cdot)\|_{L^2(\Omega)}=0.
			\]
		\end{enumerate}
	\end{definition}
\end{subequations}
\begin{remark} The pressure $\pi_n(x,t)$, as usual, was not included in the definition of a weak solution. It can be uniquely recovered from equation \ref{eq:1} by using de~Rhaam's lemma, after the existence of $(\mathbf{v}_n(x,t),\vartheta_n(x,t),\varphi_n(t),\psi_n(t))$ as in \cite{AOKh:2020:Nonl}.
\end{remark}
\begin{remark}
	Along with equations \eqref{we:equa} and \eqref{we:eqheat}, we also need 
	\begin{subequations}
		\label{eq2.1}
		\begin{align}
			&\int_{Q_{T_n}}\Big(-\mathbf{v}_n \cdot \partial_t \boldsymbol{\phi} +(\mathbf{v}_n \cdot \nabla)\mathbf{v}_n\cdot\boldsymbol{\phi} + \nu \,
			\nabla\mathbf{v}_n:
			\nabla \boldsymbol{\phi}
			-\vartheta_n\mathbf{g}\cdot\boldsymbol{\phi}
			-\varphi_n(t)\boldsymbol{\omega}({x})\cdot\boldsymbol{\phi}\Big)\, d{x} dt\nonumber
			\\
			&\quad
			-\int_\Omega \mathbf{v}_0({x}) \cdot \boldsymbol{\phi}({x},0)\, d{x} =0, \label{eq2.1a}
			\\
			&\int_{Q_{T_n}}\Big(-\vartheta_n\partial_t\zeta 
			+(\mathbf{v}_n \cdot \nabla)\vartheta_n\zeta
			+ \mu \nabla\vartheta_n\cdot\nabla \zeta -\psi_n(t)\eta({x})\zeta \Big)\, d{x} dt \nonumber
			\\& \quad
			- \int_\Omega \vartheta_0({x})\zeta({x},0)\, d{x}=0. \label{eq2.1b}
		\end{align}
	\end{subequations}
\end{remark}
\begin{remark} We introduce two types of solutions, namely, \textbf{weak} and \textbf{weak regular solutions}. This is very important, since in Theorem \ref{theorem3.1}, formally, a shifted solution lacks regularity and is a \textbf{weak solution}. The same can be said about shifted initial data, which is a final trace $\overline{\mathbf{v}}_n(\cdot,1)\in\mathbf{H}$ and $\overline{\vartheta}_n(\cdot,1)\in L^2(\Omega)$ of the rescaled solution. Here, a shifted inherits regularity properties of an original solution and is a  weak solution.
\end{remark}

\subsection{\bf An equivalent nonlocal direct problem}
Taking into account Conditions \ref{C}, 
we have
\begin{subequations}
	\begin{align}
		\label{00:eq:ht}
		& \varphi_n(t;[\mathbf{v}_n,\vartheta_n])=f_n'(t)+\left(\left(\mathbf{v}_n\cdot
		{\nabla}\right)\mathbf{v}_n,\boldsymbol{\omega}\right)+ 
		\nu
		\left( {\nabla}\mathbf{v}_n,{\nabla}\boldsymbol{\omega}
		\right)-\left(\mathbf{g}\vartheta_n, \boldsymbol{\omega}
		\right),
		\\
		&\label{00:eq:jt}
		\psi_n(t;[\mathbf{v}_n,\vartheta_n])=g_n'(t)
		+
		\left(\left(\mathbf{v}_n\cdot \nabla\right)\vartheta_n,\eta\right)
		+
		\mu\left({\nabla}\vartheta_n,{\nabla}\eta\right).
	\end{align}
\end{subequations}
Plugging \eqref{00:eq:ht} and \eqref{00:eq:jt} into (\ref{eq:1a}) and (\ref{eq:heat}), respectively, we obtain the following nonlocal direct problem with unknown solution $(\mathbf{v}_n,  \pi_n,\vartheta_n)$:
\begin{subequations}
	\label{eq:1+}
	\begin{align} 
		& \label{eq:1+a}
		\partial_t \mathbf{v}_{n}
		+\left(\mathbf{v}_n\cdot \nabla \right) \mathbf{v}_n-\nu \Delta\mathbf{v}_n
		+\nabla\pi_n=\mathbf{g}(x,t) \vartheta_n+ \varphi_n(t;[\mathbf{v}_n,\vartheta_n])\boldsymbol{\omega}(x), \ \ \ \text{ in } Q_{T_n}, 
		\\
		&\operatorname{div}\mathbf{v}_n(x,t)=0, \ \ \ \text{ in } Q_{T_n},
		\\
		&\partial_t \vartheta_{n}+\left(\mathbf{v}_n\cdot \nabla \right)\vartheta_n
		-\mu\Delta\vartheta_n=\psi_n(t;[\mathbf{v}_n,\vartheta_n])\eta(x), \ \ \ \text{ in } Q_{T_n},
	\end{align}
\end{subequations}
which is supplemented with the initial and boundary conditions \eqref{eq:2} and \eqref{eq:3}. 

The following lemma holds.
\begin{lemma}\label{le:1}
	Assume that Conditions \ref{C} are fulfilled. Then the unique solvability of the inverse problem \eqref{eq:1} is equivalent to the unique solvability of the nonlocal problem \eqref{eq:1+}.
\end{lemma}
The proof of Lemma~\ref{le:1} is similar to the proof of the analogous lemma in \cite{AK:2023}; therefore, we omit the proof.
\subsection{\bf Existence of weak and weak regular solutions} \label{Sec:LE1}
In this subsection, we prove the existence of weak and weak regular solutions of the nonlocal direct problem
\eqref{eq:1+}, which, due to Lemma~\ref{le:1}, holds true for the inverse problem \eqref{eq:1}.
\begin{theorem}\label{thm1}
	Under Conditions \ref{C}, for every $T_n>0$, problem \eqref{eq:1+} has at least one weak solution in the cylinder $Q_{T_n}$. Accordingly, the inverse problem \eqref{eq:1} has at least one weak solution, and the following estimates hold
	\begin{subequations}
		\label{estim-1+}
		\begin{align}
			&\label{estim-1a+}\left\Vert \mathbf{v}_n\right\Vert_{{L}^{\infty}\left(0,T_n; \mathbf{H}\right)}
			+\left\Vert \mathbf{v}_n\right\Vert_{{L}^{2}\left(0,T_n;\mathbf{V}\right)}
			+\left\Vert \varphi_n\right\Vert_{L^1(0,T_n)} +\left\Vert  \partial_{t} \mathbf{v}_{n}\right\Vert_{{L}^1\left(0,T_n; \mathbf{V}^{\prime}\right)} \leq C_1,
			\\
			&\label{estim-1b+}\left\Vert \vartheta_n\right\Vert_{{L}^{\infty}(0,T_n; L^{2}(\Omega))}
			+\left\Vert \vartheta_n\right\Vert_{L^2(0,T_n;W_0^{1,2}(\Omega))} +\left\Vert \psi_n\right\Vert_{L^1(0,T_n)}
			\nonumber
			\\
			&\quad +\left\Vert \partial_{t}\vartheta_{n}\right\Vert_{L^1(0,T_n;W^{-1,2}(\Omega))}\leq C_2,
		\end{align}
		where $C_1$ and $C_2$ are constants independent of $n$.
	\end{subequations}
	Moreover, if the initial data have additional regularity, i.e., $\mathbf{v}_0\in\mathbf{V}$ and $\vartheta_0\in W_0^{1,2}(\Omega)$, then problem \eqref{eq:1+} admits a weak regular solution with the estimates 
	\begin{subequations}
		\label{estim-1}
		\begin{align}
			&\label{estim-1a}\left\Vert \mathbf{v}_n\right\Vert_{{L}^{\infty}\left(0,T_n; \mathbf{V}\right)}
			+\left\Vert \mathbf{v}_n\right\Vert_{{L}^{2}\left(0,T_n;\mathbf{V}^2\right)}
			\leq C_3,
			\\
			&\label{estim-1b}\left\Vert \vartheta_n\right\Vert_{{L}^{\infty}(0,T_n; W_0^{1,2}(\Omega))}+
			\left\Vert  \vartheta_n \right\Vert_{L^2(0,T_n;W_0^{2,2}(\Omega))}\leq C_4,
		\end{align}
	\end{subequations}
	where $C_3$ and $C_4$ are positive constants depending on the data of the problem.
\end{theorem}
\begin{corollary}
	\label{cor3.1}
	The families $\{ \mathbf{v}_n\}$ and $\{\vartheta_n\}$ are integral equi-continuous. 
\end{corollary}
The proof of this corollary follows the estimates of time derivatives. We split the proof of Theorem \ref{thm1} into several subsections. We use Galerkin's method to prove the result: construct a sequence of Galerkin's approximations, obtain a priori estimates, and pass to a limit.

\subsection{\bf Galerkin's approximations}\label{Subs:GA}
A solution to the problem \eqref{eq:1+} is constructed as the limit of the sequence of finite-dimensional approximations:
\begin{align*}
	&\mathbf{v}_n(x,t)=\lim\limits_{m\rightarrow\infty}\mathbf{v}_n^{(m)}(x,t), \ \  \, \mathbf{v}_n^{(m)}(x,t)= \sum_{l=1}^{m}c_{lmn}(t)\boldsymbol{\phi}_l(x)\in\mathcal{X}^{m}=\text{span}\left\{\boldsymbol{\phi}_{1},\dots,\boldsymbol{\phi}_{m}\right\}, 
	\\
	&\vartheta_n(x,t)=\lim\limits_{m\rightarrow\infty}\vartheta_{n}^{(m)}(x,t), \ \ \, 
	\vartheta_{n}^{(m)}(x,t)= \sum_{l=1}^{m}d_{lmn}(t){\zeta}_{l}(x)\in \mathcal{Y}^{m}=\text{span}\left\{{\zeta}_{1},\dots,\zeta_{m}\right\},
\end{align*}
 where $\left\{ \boldsymbol{\phi}_{k}\right\}_{k\in\mathbb{N}}$ is the family in $\mathbf{L}^2(\Omega)$ formed by functions of $\mathbf{H}$ 
\[
-\Delta \boldsymbol{\phi}_{k}(x)=\varkappa_k\boldsymbol{\phi}_{k}(x), \ \ x \in \Omega,\ \  \left. \boldsymbol{\phi}_{k}(x)\right|_{\partial\Omega}=0,
\]
whose linear combinations are dense in $\mathbf{V}$, see \cite{Lions:1969}, and $\left\{{\zeta}_{k}\right\}_{k\in\mathbb{N}} $ is the orthonormal basis of $L^{2}(\Omega)$ composed of eigenfunctions of the following spectral problem for the Laplace operator such that $\zeta_{k}(x)\in W_0^{1,2}(\Omega)$ and  
\[
-\Delta{\zeta}_{k}(x)=\lambda_k\zeta_{k}(x), \ \ x \in \Omega,\ \  \left. {\zeta}_{k}(x)\right|_{\partial\Omega}=0.
\] 
The coefficients $c_{lmn}(t)$ and  $d_{lmn}(t)$,  $l=1,\ldots,m$ are defined as the solutions of the following Cauchy problem for the system of $2m$ ordinary nonlinear differential equations
\begin{subequations}
	\label{eq:pr2+}
	\begin{align}\label{eq:pr2}
		\begin{split}
			&\dfrac{d}{dt} (\mathbf{v}_n^{(m)},\boldsymbol{\phi}_l)+\big(( \mathbf{v}_n^{(m)} \cdot {\nabla})\mathbf{v}_n^{(m)},\boldsymbol{\phi}_l \big)
			+\nu({\nabla} \mathbf{v}_n^{(m)}, {\nabla}\boldsymbol{\phi}_l )
			\\
			&\quad=\varphi_n^{(m)}(t)(\boldsymbol{\omega}(x), \boldsymbol{\phi}_l )+(\mathbf{g}(\cdot,t)\vartheta_{n}^{(m)}, \boldsymbol{\phi}_l ),
		\end{split}
		\\
		\label{eq:pr22}
		\begin{split}
			&\dfrac{d}{dt}(\vartheta_{n}^{(m)},\zeta_l )+\big((\mathbf{v}_n^{(m)} \cdot \nabla)\vartheta_{n}^{(m)},  \zeta_l \big)+\mu ({\nabla}\vartheta_{n}^{(m)}, {\nabla}\zeta_l )
			=\psi_n^{(m)}(t)(\eta,\zeta_l ),
		\end{split}
		\\
		\label{eq:pr0}
		&\mathbf{v}_n^{(m)}(x,0)=\mathbf{v}_0^{(m)}(x), \ \ \ \vartheta_{n}^{(m)}(x,0)=\vartheta_{0}^{(m)}(x), \ \ \ x\in \Omega,\quad l=1,2,\ldots,m,
	\end{align}
	where
	\begin{align}
		\label{eq:n:ht}
		\begin{split}
			\varphi_n^{(m)}(t)&=\|\boldsymbol{\omega}^{(m)}\|_{\mathbf{H}}^{-2}\Big(f_n'(t)+\big( ( \mathbf{v}_n^{(m)} \cdot {\nabla})\mathbf{v}_n^{(m)},\boldsymbol{\omega}\big)
			+\nu\big({\nabla} \mathbf{v}_n^{(m)}, {\nabla}\boldsymbol{\omega}
			\big)
			\\
			&\quad-\big( \mathbf{g}(\cdot,t)\vartheta_{n}^{(m)}, \boldsymbol{\omega}\big)
			\Big),
		\end{split}
		\\
		\label{eq:n:jt}
		\begin{split}
			\psi_n^{(m)}(t)
			&=\|\eta^{(m)}\|_{L^2(\Omega)}^{-2}\Big(g_n'(t)+
			\big((\mathbf{v}_n^{(m)} \cdot
			\nabla)\vartheta_{n}^{(m)},\eta\big)
			+\mu
			\big({\nabla}\vartheta_{n}^{(m)},
			{\nabla}\eta\big)\Big),\end{split}
	\end{align}
\end{subequations}
with $\omega^{(m)}=\sum_{l=1}^{m}\big(\boldsymbol{\omega},\boldsymbol\phi_{l}\big)\boldsymbol{\phi}_{l}$, \ \ $\eta^{(m)}=\sum_{l=1}^{m}\big(\eta,\zeta_{l}\big)\zeta_{l},$ under the conditions
\[
\|\omega^{(m)}\|_{\mathbf{H}}^2,\,\|\eta^{(m)}\|_{L^2(\Omega)}^2>\frac12.
\]
Here 
\[ 
\mathbf{v}_0^{(m)}= \sum_{l=1}^{m}\left(\mathbf{v}_0,
\boldsymbol\phi_{l}\right)\boldsymbol\phi_{l} \quad \text{and} \quad \vartheta_{0}^{(m)}= \sum_{l=1}^{m}\left(\vartheta_0,\zeta_{l}\right) \zeta_{l}
\]
are sequences in $\mathbf{H} \cap\mathbf{V}$ and ${L}^2(\Omega)\cap W_0^{1,2}(\Omega)$, respectively, such that
\begin{align}
	\label{eq:init} 
	&\mathbf{v}_0^{(m)} \xrightarrow[m\to\infty]{} \mathbf{v}_0 \quad \text{strongly}    \text{ in }    \mathbf{H}\cap\mathbf{V},
	\\
	\label{eq:Tinit} 
	&\vartheta_{0}^{(m)} \xrightarrow[m\to\infty]{} \vartheta_0 \quad \text{strongly}   \text{ in }   {L}^2(\Omega)\cap W_0^{1,2}(\Omega).
\end{align}

According to the general theory of ordinary differential equations, the Cauchy problem (\ref{eq:pr2+}) has a continuous solutions $c_{lmn}(t)$ and $d_{lmn}(t)$ on $(0,T_m) $. By a priori estimates, which we shall establish below, the solutions can be extended to $[0,T_m]\subset[0,T]$, where $T$ is an arbitrary positive number.

\subsection{\bf A priori estimates} \label{fir:est}
\begin{lemma}\label{1Globest:1} Suppose that Conditions \ref{C} are satisfied. Let $\mathbf{v}_0(x) \in \mathbf{H}$ and $\vartheta_0(x) \in L^2(\Omega)$. Then, for all $T>0$, the following a priori estimates hold for all $t\in(0,T_n)$:
	\begin{subequations}
		\label{fir-est1}
		\begin{align}
			\label{fir-est1a}
			&\Vert \mathbf{v}_n^{(m)}\Vert_{{L}^{\infty}\left(0,T_n; \mathbf{H}\right)}^{2}+ \Vert \vartheta_{n}^{(m)} \Vert_{{L}^{\infty}
				(0,T_n; L^2(\Omega))}^{2}+
			\Vert \mathbf{v}_n^{(m)}\Vert_{{L}^2(0,T_n;
				\mathbf{V})}^{2}+ \Vert \vartheta_{n}^{(m)}
			\Vert_{L^2(0,T_n,W_0^{1,2}(\Omega))}^{2}\leq C,
			\\
			\label{fir-est1b}
			&\Vert \varphi_n^{(m)}\Vert_{L^1(0,T_n)}+
			\Vert \psi_{n}^{(m)}\Vert_{L^1(0,T_n)}\leq C\left(1+\int_0^{T_n}|f_n'(t)|\,dt+\int_0^{T_n}|g_n'(t)|\,dt\right),
			\\
			\label{fir-est1c}
			&\Vert \varphi_n^{(m)}\Vert_{L^2(0,T_n)}^{2}+
			\Vert \psi_{n}^{(m)} \Vert_{L^2(0,T_n)}^{2}\leq C\left(1+\int_0^{T_n}|f_n'(t)|^2\,dt+\int_0^{T_n}|g_n'(t)|^2\,dt\right),
		\end{align}
	\end{subequations}
	where the constant $C$ does not depend on $n$ and $m$.
\end{lemma}
\begin{proof}
	Multiply equation (\ref{eq:pr2}) by $c_{lmn}(t)$ and equation (\ref{eq:pr22}) by $d_{lmn}(t)$, and sum up from 1 to $m$. The sum of these two equations is the following:
	\begin{align}\label{en:eq}
		\begin{split}& \frac{1}{2}\frac{d}{dt}
			\left(\Vert \mathbf{v}_n^{(m)}\Vert _{\mathbf{H}}^{2}
			+\Vert \vartheta_{n}^{(m)} \Vert_{L^2(\Omega)}^{2}\right)+\nu \Vert\mathbf{\nabla}
			\mathbf{v}_n^{(m)} \Vert_{\mathbf{H}}^{2}
			+\mu \Vert\nabla
			\vartheta_{n}^{(m)} \Vert_{\mathbf{L}^2(\Omega)}^{2}
			\\
			&
			=\varphi_n^{(m)}(t)f_n(t)+
			\big( \mathbf{g}(\cdot,t)\vartheta_{n}^{(m)},\mathbf{v}_n^{(m)} \big)
			+\psi_n^{(m)}(t)g_n(t)=: I_{11}+I_{12}+I_{13}.
		\end{split}
	\end{align}
	First, we apply H\"older's and Young's inequalities to \eqref{eq:n:ht} and \eqref{eq:n:jt}, in order to obtain
	\begin{align*}
		%\label{est:fn}
		\begin{split}
			|I_{11}|&=\big|{{\varphi_n^{(m)}(t)}}\big||f_n(t)|\leq 2{F}
			\Big(|f_n'(t)|
			+\nu \Vert
			\nabla \mathbf{v}_{n}^{(m)}\Vert_{\mathbf{H}}\left\Vert
			\nabla\boldsymbol\omega\right\Vert _{\mathbf{H}}
			+g_0 \Vert
			\vartheta_{n}^{(m)} \Vert_{L^2(\Omega)}\left\Vert \boldsymbol\omega\right\Vert _{\mathbf{H}}
			\\
			&\quad+\Vert \mathbf{v}_{n}^{(m)}\Vert_{\mathbf{H}}
			\Vert\nabla\mathbf{v}_{n}^{(m)}\Vert_{\mathbf{H}}
			\left\Vert\boldsymbol\omega\right\Vert_{\widetilde{\mathbf{L}}^\infty}\Big)
			\\
			&
			\leq 2{F}\Big(|f_n'(t)|+\frac{\nu\delta_1}{2}\Vert \nabla \mathbf{v}_{n}^{(m)}\Vert_{\mathbf{H}}^2
			+\frac{\nu}{2\delta_1}\Vert \nabla\boldsymbol\omega \Vert_{\mathbf{H}}^2 
			+\frac{g_0}2 \Vert
			\vartheta_{n}^{(m)}\Vert_{L^2(\Omega)}^2+\frac{g_0}2\Vert \boldsymbol\omega \Vert_{\mathbf{H}}^2
			\\
			&\quad
			+\frac{\nu\delta_1}{2}\Vert \nabla \mathbf{v}_{n}^{(m)}\Vert_{\mathbf{H}}^2
			+\frac{1}{2\nu\delta_1}\Vert \mathbf{v}_{n}^{(m)}\Vert_{\mathbf{H}}^2
			\left\Vert\boldsymbol\omega\right\Vert_{\widetilde{\mathbf{L}}^\infty}^2\Big),
		\end{split}
		\\
		%\label{est:gTv}
		\begin{split}
			|I_{12}|&=\big|\big( \mathbf{g}(\cdot,t)\vartheta_{n}^{(m)}, \mathbf{v}_n^{(m)} \big)\big|
			\leq
			g_0 \Vert
			\vartheta_{n}^{(m)} \Vert_{L^2(\Omega)}\Vert  \mathbf{v}_n^{(m)}\Vert_{\mathbf{H}}
			\leq\frac{g_0}{2}\left(\Vert  \mathbf{v}_n^{(m)}\Vert_{\mathbf{H}}^2+\Vert
			\vartheta_{n}^{(m)}\Vert_{L^2(\Omega)}^2\right),
		\end{split}
		\\
		%\label{est:jn}
		\begin{split}
			|I_{13}|&=\big|\psi_n^{(m)}\big|\left\vert g_n \right\vert\leq 2{G}\Big(|g_n'(t)|+\mu \Vert \nabla \vartheta_n^{(m)} \Vert_{\mathbf{L}^2(\Omega)}\left\Vert \nabla\eta\right\Vert_{\mathbf{L}^2(\Omega)} 
			\\
			&\quad+\Vert\mathbf{v}_n^{(m)}\Vert_{\mathbf{H}}
			\Vert\nabla\vartheta_n^{(m)}\Vert_{\mathbf{L}^2(\Omega)}\Vert\eta\Vert _{L^\infty(\Omega)}\Big)
			\\
			&\leq 2{G}\Big(|g_n'(t)|+\frac{\mu\delta_2}2 \Vert \nabla \vartheta_n^{(m)} \Vert_{\mathbf{L}^2(\Omega)}^2
			+\frac{\mu}{2\delta_2} \Vert \nabla\eta \Vert _{\mathbf{L}^2(\Omega)}^2 
			+\frac{\mu \delta_2}2
			\Vert \nabla \vartheta_n^{(m)} \Vert_{\mathbf{L}^2(\Omega)}^2
			\\
			&\quad+\frac1{2\mu \delta_2}\Vert\mathbf{v}_n^{(m)} \Vert_{\mathbf{H}}^2\left\Vert\eta\right\Vert _{L^\infty(\Omega)}^2\Big),
		\end{split}
	\end{align*}
	respectively. Let $Y_1(t)= \Vert \mathbf{v}_n^{(m)}\Vert _{\mathbf{H}}^{2}
	+\Vert \vartheta_{n}^{(m)} \Vert_{L^2(\Omega)}^{2}$, $\delta_1=\frac{1}{4F}$, and $\delta_2=\frac{1}{4G}$. We arrive at
	\begin{align*}
		%\label{est:fn}
		\begin{split}
			&\frac{d}{d t} Y_1(t)+\nu \Vert \nabla\mathbf{v}_n^{(m)} \Vert _{\mathbf{H}}^{2}+\mu \Vert \nabla\vartheta_n^{(m)} \Vert _{\mathbf{L}^2(\Omega)}^{2}
			\leq 4G|g_n'(t)|+4F|f_n'(t)|\\
			&\quad+\frac{2F\nu}{\delta_1}\Vert \nabla\boldsymbol\omega \Vert_{\mathbf{H}} 
			+\max\bigg\{\frac{2 F\|\boldsymbol\omega\|_{\widetilde{\mathbf{L}}^\infty}^2}{ \nu \delta_1}+g_0+\frac{2G\left\Vert\eta\right\Vert_{L^\infty(\Omega)}^2}{\mu\delta_2},2g_0 F+1\bigg\}Y_1(t)
			\\
			& \qquad+2{g_0F}\left\Vert \boldsymbol\omega\right\Vert_{\mathbf{H}}^2+\frac{2\mu G}{\delta_2}\left\Vert \nabla\eta\right\Vert_{\mathbf{L}^2(\Omega)}^2,
		\end{split}
	\end{align*}
	Using Gronwall's inequality, we obtain 
	\[
	\sup\limits_{t\in(0,T_n)}\left(\Vert \mathbf{v}_n^{(m)}\Vert_{\mathbf{H}}^{2} +\Vert\vartheta_n^{(m)}\Vert_{L^2(\Omega)}^{2}\right)+ \Vert\nabla \mathbf{v}_n^{(m)} \Vert_{{L}^2(Q_{T_n})}^{2}+ \Vert\nabla \vartheta_n^{(m)} \Vert_{{L}^2(Q_{T_n})}^{2}\leq  C,
	\]
	where $C$ does not depend on $n$. 
	
	Let us estimate the norms of $\varphi_n^{(m)}$ and $\psi_n^{(m)}$. From equation \eqref{eq:n:ht}, using H\"older's inequality, we obtain
	\begin{align*}
		&
		\int_0^{T_n}|\varphi_n^{(m)}(t)|\,dt\leq \int_0^{T_n}|f_n'(t)|\,dt
		+\nu\int_0^{T_n}\Vert\nabla \mathbf{v}_{n}^{(m)}\Vert_{\mathbf{H}}\Vert
		\nabla\boldsymbol\omega\Vert _{\mathbf{H}} \, d t
		\\
		&\qquad+g_0\int_0^{T_n} \Vert
		\vartheta_{n}^{(m)}\Vert_{L^2(\Omega)}\left\Vert \boldsymbol\omega\right\Vert _{\mathbf{H}}\,dt
		+\int_0^{T_n} \Vert \mathbf{v}_{n}^{(m)}\Vert_{\mathbf{H}}
		\Vert\nabla\mathbf{v}_{n}^{(m)}\Vert_{\mathbf{H}}
		\Vert\boldsymbol\omega \Vert _{\widetilde{\mathbf{L}}^\infty}\,dt
		\\
		&\leq \int_0^{T_n}|f_n'(t)|\,dt
		+\nu \sqrt{T_n} \Vert\nabla \mathbf{v}_{n}^{(m)}\Vert_{L^2(0,{T_n};\mathbf{H})}\Vert
		\nabla\boldsymbol\omega\Vert_{\mathbf{H}}
		+g_0 T_n \Vert
		\vartheta_{n}^{(m)}\Vert_{L^\infty(0,T_n;L^2(\Omega))} \left\Vert \boldsymbol\omega\right\Vert_{\mathbf{H}}
		\\
		&\qquad+\sqrt{T_n} \Vert \mathbf{v}_{n}^{(m)}\Vert_{L^\infty(0,{T_n};\mathbf{H})}
		\Vert\nabla\mathbf{v}_{n}^{(m)}\Vert_{L^2(0,{T_n};\mathbf{H})}
		\left\Vert\boldsymbol\omega\right\Vert_{\widetilde{\mathbf{L}}^\infty}.
	\end{align*}
	Similarly, from \eqref{eq:n:jt}, we get \begin{align*}
		&\int_0^{T_n}|\psi_n^{(m)}(t)|\,dt\leq \int_0^{T_n}|g_n'(t)|\,dt+\mu \sqrt{T_n}\Vert \nabla \vartheta_n^{(m)}\Vert_{L^2(0,{T_n};\mathbf{L}^2(\Omega))}
		\Vert \nabla\eta \Vert_{\mathbf{L}^2(\Omega)} 
		\\
		&\qquad+\sqrt{T_n} \, \Vert \mathbf{v}_{n}^{(m)}\Vert_{L^\infty(0,{T_n};\mathbf{H})} \Vert\nabla\vartheta_n^{(m)}\Vert_{L^2(0,{T_n};\mathbf{L}^2(\Omega))}\Vert\eta\Vert _{L^\infty(\Omega)},
	\end{align*}
	and the estimate \eqref{fir-est1b} follows. Moreover, from \eqref{eq:n:ht} and \eqref{eq:n:jt}, we also obtain
	\begin{align*}
		&
		\left(\int_0^{T_n}|\varphi_n^{(m)}(t)|^2\,dt\right)^{\frac12}\leq \left(\int_0^{T_n}|f_n'(t)|^2\,dt\right)^{\frac12}
		+\nu \Vert\nabla \mathbf{v}_{n}^{(m)}\Vert_{L^2(0,{T_n};\mathbf{H})}\Vert
		\nabla\boldsymbol\omega\Vert_{\mathbf{H}}
		\\
		&\quad+g_0 \sqrt{T_n} \, \Vert
		\vartheta_{n}^{(m)}\Vert_{L^\infty(0,T_n;L^2(\Omega))} \left\Vert \boldsymbol\omega\right\Vert _{\mathbf{H}}
		+  \Vert \mathbf{v}_{n}^{(m)}\Vert_{L^\infty(0,T_n;\mathbf{H})}
		\Vert\nabla\mathbf{v}_{n}^{(m)}\Vert_{L^2(0,{T_n};\mathbf{H})}
		\left\Vert\boldsymbol\omega\right\Vert_{\widetilde{\mathbf{L}}^\infty},
	\end{align*}
	and
	\begin{align*}
		&\left(\int_0^{T_n}|\psi_n^{(m)}(t)|^2\,dt\right)^{\frac12}\leq \left(\int_0^{T_n}|g_n'(t)|^2\,dt\right)^{\frac12}+\mu\Vert \nabla \vartheta_n^{(m)}\Vert_{L^2(0,{T_n};\mathbf{L}^2(\Omega))}\left\Vert \nabla\eta\right\Vert _{\mathbf{L}^2(\Omega)} 
		\\
		&\quad+   \Vert \mathbf{v}_{n}^{(m)}\Vert_{L^\infty(0,T_n;\mathbf{H})}
		\Vert\nabla\vartheta_n^{(m)}\Vert_{L^2(0,{T_n};\mathbf{L}^2(\Omega))}\left\Vert\eta\right\Vert _{L^\infty(\Omega)},
	\end{align*}
	respectively. Therefore, we get \eqref{fir-est1c}.
\end{proof}
\begin{remark}
	The estimates $\int_0^{T_n}|\varphi_n^{(m)}(t)|\,dt$ and $\int_0^{T_n}|\psi_n^{(m)}(t)|\,dt$ are uniform in $m$ and $n$. On the other hand, the estimates $\int_0^{T_n}|\varphi_n^{(m)}(t)|^{2}\,dt$ and $\int_0^{T_n}|\psi_n^{(m)}(t)|^{2}\,dt$ are uniform in $m$, but not uniform in $n$.  {However, the estimates $\int_{1/n}^{T_n}|\varphi_n^{(m)}(t)|^{2}\,dt$ and $\int_{1/n}^{T_n}|\psi_n^{(m)}(t)|^{2}\,dt$ are uniformly bounded in both $m$ and $n$.}
\end{remark}
\begin{remark}\label{Remark 3.5}
	The Ladyzhenskaya inequality and the estimates \eqref{fir-est1a} imply that $\mathbf{v}_n^{(m)} \in L^4(0,{T_n};\widetilde{\mathbf{L}}^4)$ and $\vartheta_n^{(m)} \in L^4(0,{T_n};{L}^4(\Omega))$. Indeed,
	\begin{align*}
		\int_0^{T_n}\Vert\mathbf{v}_{n}^{(m)}(\cdot,t)\Vert_{\widetilde{\mathbf{L}}^4}^4 \, d t & \leq  \int_0^{T_n} 2 \Vert\mathbf{v}_{n}^{(m)}(\cdot,t)\Vert_{\mathbf{H}}^{2} \Vert \nabla \mathbf{v}_{n}^{(m)}(\cdot,t)\Vert_{\mathbf{H}}^{2}  \, d t
		\\
		& 
		\leq 2  \, \underset{t\in(0,T_n)}{\rm esssup} \, \Vert\mathbf{v}_{n}^{(m)}(\cdot,t)\Vert_{\mathbf{H}}^{2} \,  \int_0^{T_n}\Vert \nabla \mathbf{v}_{n}^{(m)}(\cdot,t)\Vert_{\mathbf{H}}^{2} \, d t\leq C,
	\end{align*}
	and similarly,
	\begin{align*}
		\int_0^{T_n}\Vert\vartheta_{n}^{(m)}(\cdot,t)\Vert_{L^4(\Omega)}^4  \, dt\leq 2  \, \underset{t\in(0,T_n)}{\rm esssup} \, \Vert \vartheta_{n}^{(m)}(\cdot,t)\Vert_{L^2(\Omega)}^{2} \,  \int_0^{T_n}\Vert \nabla \vartheta_{n}^{(m)}(\cdot,t)\Vert_{\mathbf{L}^2(\Omega)}^{2} \, d t\leq C,
	\end{align*}
	where $C$ does not depend on $m$ and $n$.
\end{remark}
\begin{lemma}\label{Lemma 3.5}
	Suppose that Lemma \ref{1Globest:1} holds. Then the following a priori estimates hold for all $t\in
	(0,T_n)$:
	\begin{equation}
		\label{eq4.19}
		\begin{split}
			&\Vert \partial_{t}\mathbf{v}_{n}^{(m)}\Vert_{{L}^1
				\left(0,T_n; \mathbf{V}^\prime\right)}^{2}\leq C<\infty,\quad \Vert \partial_{t}\mathbf{v}_{n}^{(m)}\Vert_{{L}^2
				\left(0,T_n; \mathbf{V}^\prime\right)}^{2}\leq C\left(1+\int_0^{T_n}|f_n'(t)|^2\,dt\right),
		\end{split}
	\end{equation}
	where $C$ does not depend on $m$ and $n$.
\end{lemma}
\begin{proof} 
	Equation \eqref{eq:pr2} yields the following integral:
	\begin{align}
		\label{eq:lim1:nz1}
		\begin{split}
			&
			\big| \langle \partial_{t}\mathbf{v}_{n}^{(m)}(\cdot,t),\boldsymbol{\phi}(\cdot)\rangle \big|=
			-\int_\Omega(\mathbf{v}_n^{(m)}\cdot{\nabla})\mathbf{v}_n^{(m)}\cdot\boldsymbol{\phi}\,dx
			-\nu\int_\Omega\nabla\mathbf{v}_{n}^{(m)}:{\nabla}\boldsymbol{\phi}\,dx
			\\
			&\quad +\int_\Omega \vartheta_{n}^{(m)} \mathbf{g}(x,t)\cdot\boldsymbol{\phi}\,dx+
			\varphi_n^{(m)}(t)\int_{\Omega}\boldsymbol{\omega}(x)\cdot\boldsymbol{\phi}\,dx,
		\end{split}
	\end{align}
	{where $\boldsymbol{\phi}$ is a linear combination of the eigen functions $\{\boldsymbol{\phi}_l\}$, $l=1,...,m$.}
	We estimate the terms on the right-hand side of \eqref{eq:lim1:nz1} using H\"older's and Ladyzhenskaya's inequalities as
	\begin{align}
		\label{est:I21}
		\begin{split}
			\left|\int_\Omega(\mathbf{v}_n^{(m)}
			\cdot{\nabla})\boldsymbol{\phi}\cdot\mathbf{v}_n^{(m)}\,dx\right|&\leq 
			\|\mathbf{v}_n^{(m)}\|_{\widetilde{\mathbf{L}}^4}^2\|\nabla\boldsymbol{\phi}\|_{\mathbf{H}}
			\\&  \leq 
			\sqrt{2}\|\mathbf{v}_n^{(m)}\|_{\mathbf{H}}\|\nabla\mathbf{v}_n^{(m)}\|_{\mathbf{H}}\|\nabla\boldsymbol{\phi}\|_{\mathbf{H}},
		\end{split}
		\\
		\label{est:I22}
		\begin{split}
			\Big| \nu\int_\Omega\nabla\mathbf{v}_{n}^{(m)}:{\nabla}\boldsymbol{\phi}\,dx \Big|&\leq 
			\nu\|\nabla\mathbf{v}_n^{(m)}\|_{\mathbf{H}}\|\nabla\boldsymbol{\phi}\|_{\mathbf{H}},
		\end{split}
		\\
		\label{est:I23}
		\begin{split}
			\Big| \int_\Omega \vartheta_{n}^{(m)} \mathbf{g}(x,t)\cdot\boldsymbol{\phi}\,dx \Big| &\leq 
			g_0
			\|\vartheta_{n}^{(m)}\|_{{L}^2(\Omega)}\|\boldsymbol{\phi}\|_{\mathbf{H}},
		\end{split}
		\\
		\label{est:I24}
		\begin{split}
			\Big|\varphi_n^{(m)}(t)\int_{\Omega}\boldsymbol{\omega}(x) \cdot \boldsymbol{\phi}\,dx \Big| &\leq 
			\big|\varphi_n^{(m)}(t) \big| \|\boldsymbol{\omega}\|_{\mathbf{H}}\|\boldsymbol{\phi}\|_{\mathbf{H}}.
		\end{split}
	\end{align}
	These estimates lead to the following estimate:
	\begin{align}\label{Time}
		\big| \langle \partial_{t}\mathbf{v}_{n}^{(m)}(\cdot,t),\boldsymbol{\phi}(\cdot)\rangle \big|\leq \Upsilon_n^{(m)}(t)\|\nabla\boldsymbol{\phi}\|_{\mathbf{H}},
	\end{align}
	where 
	\begin{align*}
		\Upsilon_n^{(m)}(t)=\sqrt{2}\|\mathbf{v}_n^{(m)}\|_{\mathbf{H}}\|\nabla\mathbf{v}_n^{(m)}\|_{\mathbf{H}}+\nu\|\nabla\mathbf{v}_n^{(m)}\|_{\mathbf{H}}+\frac{g_0}{\sqrt{\varkappa_1}}
		\|\vartheta_{n}^{(m)}\|_{{L}^2(\Omega)}+\frac{1}{\sqrt{\varkappa_1}}|\varphi_n^{(m)}(t)\|\boldsymbol{\omega}\|_{\mathbf{H}}.
	\end{align*}
	The relation \eqref{Time} yields
	\begin{align*}
		\|\partial_{t}\mathbf{v}_{n}^{(m)}(\cdot,t) \|_{\mathbf{V}^\prime} \leq \Upsilon_n^{(m)}(t).
	\end{align*}
	We can easily see that $\Upsilon_n^{(m)}(t)$ is uniformly bounded in $L^1(0,T_n)$ in the parameters $m$ and $n$, but it is bounded in $L^2(0,T_n)$ uniformly only in the parameter $m$. This implies that $\partial_{t}\mathbf{v}_{n}^{(m)}$ is uniformly bounded in the space $L^1(0,T_n;\mathbf{V}^\prime)$ in the parameters $m$ and $n$, but it is  uniformly bounded in $L^2(0,T_n;\mathbf{V}^\prime)$ only in the parameter $m$ not in $n$. 
\end{proof}
\begin{lemma}\label{Lemma 3.6} Suppose that Lemma \ref{1Globest:1} holds. Then the following a priori estimates hold for all $t\in(0,T_n)$:
	\begin{equation}\label{sec-est13}
		\begin{split}
			&\Vert\partial_t\vartheta_{n}^{(m)}\Vert_{L^1(0,T_n;W^{-1,2}(\Omega))}\leq C, \ \ 
			\Vert \partial_t\vartheta_{n}^{(m)}\Vert_{L^2(0,T_n;W^{-1,2}(\Omega))}
			\leq C\Big(1+\Big(\int_0^{T_n}|g_n'(t)|^2\,dt\Big)^{\frac12}\Big),
		\end{split}
	\end{equation}
	where $C$ does not depend on $m$ and $n$.
\end{lemma}
\begin{proof}
	Applying H\"older's, Sobolev's, and Poincar\'e's inequalities, we obtain from the equation \eqref{eq:pr22}:
	\begin{align}
		\label{en:Teq2}
		\big|\langle \partial_t\vartheta_{n}^{(m)}(\cdot,t),\zeta(\cdot)\rangle\big|&=\Big|-\mu(\nabla
		\vartheta_{n}^{(m)},\nabla\zeta)+
		\big((\mathbf{v}_n^{(m)}\cdot\nabla)\zeta,\vartheta_{n}^{(m)}\big)
		+\psi_n^{(m)}(t)(\eta(x),\zeta) \Big|
		\nonumber \\& \leq
		\mu\|\nabla\vartheta_{n}^{(m)}\|_{\mathbf{L}^2(\Omega)}\|\nabla\zeta\|_{\mathbf{L}^2(\Omega)} +\|\mathbf{v}_n^{(m)}\|_{\widetilde{\mathbf{L}}^4}
		\|\vartheta_{n}^{(m)}\|_{L^4(\Omega)}\|\nabla\zeta\|_{\mathbf{L}^2(\Omega)}\nonumber\\& \quad+|\psi_n^{(m)}(t)|\|\eta\|_{L^2(\Omega)}\|\zeta\|_{L^2(\Omega)}
		\nonumber \\& \leq 
		\Xi_n^{(m)}(t)\|\nabla\zeta\|_{\mathbf{L}^2(\Omega)},
	\end{align}
	{ where $\zeta$ is a linear combination of the eigen functions $\{\zeta_l\}$, $l=1,...,m$,} and 
	\begin{align*}
		\Xi_n^{(m)}(t)= \big(\mu + C\|\nabla \mathbf{v}_{n}^{(m)}\|_{\mathbf{H}}\big)\|\nabla\vartheta_{n}^{(m)}\|_{\mathbf{L}^2(\Omega)}+\frac{1}{\sqrt{\lambda_1}}|\psi_n^{(m)}(t)|\|\eta\|_{L^2(\Omega)}.
	\end{align*}
	Inequality \eqref{en:Teq2} yields
	\begin{align}\label{Timervartheta} 
		\|\partial_t\vartheta_{n}^{(m)} (\cdot,t)\|_{W^{-1,2}(\Omega)} \leq \Xi_n^{(m)}(t).
	\end{align}
	It is deduced that \(\Xi_n^{(m)}\in L^1(0,T_n)\) uniformly in $m$ and \(n\), whereas \( \Xi_n^{(m)}\in L^{2}(0,T_n)\) uniformly in $m$, but not in $n$. Therefore, \eqref{Timervartheta} yields that \(\partial_t\vartheta_{n}^{(m)} \in L^1(0,T_n;W^{-1,2}(\Omega))\) uniformly in $m$ and \(n\), while \( \partial_t\vartheta_{n}^{(m)} \in L^{2}(0,T_n;W^{-1,2}(\Omega))\) uniformly in $m$, but not in $n$.
\end{proof}

\begin{lemma}\label{lem:estT2} Let Conditions \ref{C} be satisfied and let $\mathbf{v}_0(x) \in \mathbf{V}$. Then the following estimates hold:
	\begin{subequations}
		\begin{align}
			&\label{sec-est13+}
			\underset{t\in(0,T_n)}{\rm sup} \, \left\Vert\nabla\mathbf{v}_{n}^{(m)}(\cdot,t)\right\Vert_{\mathbf{H}}^{2}
			+
			\Vert \Delta\mathbf{v}_{n}^{(m)}\Vert_{L^2( Q_{T_n})}^{2}
			\leq C<\infty,
			\\
			\label{sec-est13++}
			& \underset{t\in(0,T_n)}{\rm sup} \,\left\Vert\mathbf{v}_{n}^{(m)}(\cdot,t)\right\Vert_{\widetilde{\mathbf{L}}^4}
			\leq C<\infty,
		\end{align}
		where $C$ does not depend on $m$ and $n$.
	\end{subequations}
\end{lemma}
\begin{proof}
	The proof is similar to the analogical results in \cite{AK:2023}; therefore, we omit the proof. Multiply (\ref{eq:pr2}) by $\varkappa_l c_{lmn}(t)$, and sum up the results by $l$ from $1$ to $m$:
	\begin{align}
		\label{en:eqT2+}
		\begin{split}
			\frac{1}{2}\frac{d}{dt}\Vert \nabla\mathbf{v}_{n}^{(m)}\Vert_{\mathbf{H}}^{2}+\nu\Vert \Delta\mathbf{v}_{n}^{(m)}\Vert_{\mathbf{H}}^{2}
			&=\big((\mathbf{v}_n^{(m)}\cdot\nabla)\mathbf{v}_{n}^{(m)},\Delta\mathbf{v}_{n}^{(m)}\big)+(\mathbf{g}\vartheta_n^{(m)},-\Delta\mathbf{v}_{n}^{(m)})
			\\
			&\quad+\varphi_n^{(m)}(t)(\nabla \boldsymbol{\omega},\nabla\mathbf{v}_{n}^{(m)})=: \sum\limits_{i=1}^3 I_i.
		\end{split}
	\end{align}
	We estimate the terms on the right-hand side of \eqref{en:eqT2+}. Using H\"older's, Agmon's, Cattabriga-Solonnikov's, and Young's inequalities, we obtain 
	\begin{align}
		&\text{the estimate for the first integral}\nonumber
		\\
		\label{inq:J3:d3}
		\begin{split}
			&|I_1|=\Big|\big(
			(\mathbf{v}_n^{(m)} \cdot \nabla)\mathbf{v}_n^{(m)},
			\Delta\mathbf{v}_n^{(m)}\big)\Big|
			\leq 
			\|\mathbf{v}_n^{(m)}\|_{\widetilde{\mathbf{L}}^\infty}\| \nabla\mathbf{v}_n^{(m)}\|_{\mathbf{H}}
			\|\Delta\mathbf{v}_n^{(m)}\|_{\mathbf{H}}
			\\
			&\leq C_{Ag}\|\mathbf{v}_n^{(m)} \|_{\mathbf{H}}^\frac{1}{2}\| \nabla\mathbf{v}_n^{(m)} \|_{\mathbf{H}}
			\|\Delta\mathbf{v}_n^{(m)} \|_{\mathbf{H}}^\frac{3}{2}
			\\
			&\leq \frac{\nu}{4} \|\Delta\mathbf{v}_n^{(m)} \|_{\mathbf{H}}^2+ \frac{C_{Ag}}{\nu}\|\mathbf{v}_n^{(m)} \|_{\mathbf{H}}^2 \| \nabla\mathbf{v}_n^{(m)} \|_{\mathbf{H}}^4,
		\end{split}
		\\
		&\text{the estimate for the second integral}\nonumber
		\\
		\begin{split}
			&|I_2|=\Big|(\mathbf{g}\vartheta_n^{(m)},-\Delta\mathbf{v}_{n}^{(m)})\Big|
			\leq g_0\|\vartheta_n^{(m)}\|_{L^2(\Omega)}\|\Delta\mathbf{v}_{n}^{(m)}\|_{\mathbf{H}}
			\\
			&\leq 
			\frac{\nu}{4}\|\Delta\mathbf{v}_{n}^{(m)}\|_{\mathbf{H}}^2
			+\frac{g_0^2}{\nu}\|\vartheta_n^{(m)}\|_{L^2(\Omega)}^2,
		\end{split}
		\\
		&\text{the estimate for the third integral}\nonumber
		\\
		\begin{split}
			&|I_3|=\Big|\varphi_n^{(m)}(t)(\nabla \boldsymbol{\omega},\nabla\mathbf{v}_{n}^{(m)})\Big| \leq
			|\varphi_n^{(m)}(t)|\|\nabla \boldsymbol{\omega}\|_{\mathbf{H}}\|\nabla\mathbf{v}_{n}^{(m)}\|_{\mathbf{H}}
			\\
			&\leq \frac{|\varphi_n^{(m)}(t)|}2\left(\|\nabla\boldsymbol{\omega}\|_{\mathbf{H}}^2+\|\nabla\mathbf{v}_{n}^{(m)}\|_{\mathbf{H}}^2\right).
		\end{split}
	\end{align}
	Substituting the above estimates  in \eqref{en:eqT2+}, we obtain
	\begin{align*}
		&   \frac{d}{dt}\Vert \nabla\mathbf{v}_{n}^{(m)}(\cdot,t)\Vert_{\mathbf{H}}^{2}+\nu\Vert \Delta\mathbf{v}_{n}^{(m)}(\cdot,t)\Vert_{\mathbf{H}}^{2}
		\leq\Big(\frac{2 C_{Ag}}{\nu}\|\mathbf{v}_n^{(m)}(\cdot,t)\|_{\mathbf{H}}^{2}
		\|\nabla\mathbf{v}_n^{(m)}(\cdot,t)\|_{\mathbf{H}}^{2}
		\\
		&
		+
		|\varphi_n^{(m)}(t)|\Big)\|\nabla\mathbf{v}_n^{(m)}(\cdot,t)\|_{\mathbf{H}}^{2}
		+|\varphi_n^{(m)}(t)|\|\nabla\boldsymbol{\omega}\|_{\mathbf{H}}^{2}+\frac{2g_0^2}{\nu}
		\|\vartheta_n^{(m)}(\cdot,t)\|_{L^2(\Omega)}^{2},
	\end{align*}
	for all $t\in [0,T_n]$. By Gr\"onwall's inequality, we get the uniform estimate \eqref{sec-est13+}. The Ladyzhenskaya inequality yields
	\begin{align*}
		\underset{t\in(0,T_n)}{\rm sup} \,\left\Vert\mathbf{v}_{n}^{(m)}(\cdot,t)\right\Vert_{\widetilde{\mathbf{L}}^4} & \leq   \, \, \underset{t\in(0,T_n)}{\rm sup} \, \Big(2^\frac{1}{4} \left\Vert\mathbf{v}_{n}^{(m)}(\cdot,t)\right\Vert_{\mathbf{H}}^\frac{1}{2} \left\Vert \nabla \mathbf{v}_{n}^{(m)}(\cdot,t)\right\Vert_{\mathbf{H}}^\frac{1}{2} \Big)
		\\
		& \leq 2^\frac{1}{4} \,   \underset{t\in(0,T_n)}{\rm sup} \, \left\Vert\mathbf{v}_{n}^{(m)}(\cdot,t)\right\Vert_{\mathbf{H}}^\frac{1}{2} \, \underset{t\in(0,T_n)}{\rm sup} \, \left\Vert \nabla \mathbf{v}_{n}^{(m)}(\cdot,t)\right\Vert_{\mathbf{H}}^\frac{1}{2} \leq C,
	\end{align*}
	which completes the proof of the lemma.
\end{proof}
\begin{lemma}\label{lem:estT21} Let Conditions \ref{C} be satisfied and let $\vartheta_0(x)\in W_0^{1,2}(\Omega)$. Then the following estimates hold: 
	\begin{subequations}
		\begin{align}\label{sec-est20}
			&\, \underset{t\in(0,T_n)}{\rm sup} \, \left\Vert\nabla\vartheta_{n}^{(m)}(\cdot,t)\right\Vert_{\mathbf{L}^2(\Omega)}^{2}
			+\Vert \Delta\vartheta_{n}^{(m)}\Vert_{L^2( Q_{T_n})}^{2}\leq C,
			\\
			\label{sec-est19}
			& \, \underset{t\in(0,T_n)}{\rm sup} \,
			\left\Vert\vartheta_{n}^{(m)}(\cdot,t)\right\Vert_{L^4(\Omega)}
			\leq C,
		\end{align}
	\end{subequations}
	where $C$ does not depend on $m$ and $n$.
\end{lemma}
\begin{proof}
	Multiplying (\ref{eq:pr22}) by $\lambda_l d_{lmn}(t)$ and summing up the results by $l$ from $1$ to $m$, we have
	\begin{align}
		\label{en:eqT2}
		\begin{split}
			&\frac{1}{2}\frac{d}{dt}\Vert \nabla\vartheta_{n}^{(m)}\Vert_{\mathbf{L}^2(\Omega)}^{2}+\mu \Vert \Delta\vartheta_{n}^{(m)}\Vert_{L^2(\Omega)}^{2}=\big((\mathbf{v}_n^{(m)}\cdot\nabla) \vartheta_{n}^{(m)}, -\Delta\vartheta_{n}^{(m)}\big)
			\\
			&\qquad+\psi_n^{(m)}(t)(\eta, -\Delta\vartheta_{n}^{(m)}).
		\end{split}
	\end{align}
	We deal with the following estimates: 
	\begin{align}
		\begin{split}
			&\big|\big((\mathbf{v}_n^{(m)}\cdot\nabla)\vartheta_{n}^{(m)},  -\Delta\vartheta_{n}^{(m)}\big)\big|\leq \Vert \Delta\vartheta_{n}^{(m)}\Vert_{L^2(\Omega)} \Vert \mathbf{v}_n^{(m)}\Vert_{\widetilde{\mathbf{L}}^4} \Vert \nabla\vartheta_{n}^{(m)}\Vert_{\mathbf{L}^4(\Omega)} 
			\\
			&\leq 
			C \Vert \Delta\vartheta_{n}^{(m)}\Vert_{L^2(\Omega)}^{\frac{3}{2}} \, \,  \underset{t\in(0,T_n)}{\rm sup} \, \Vert  \mathbf{v}_n^{(m)}(\cdot,t)\Vert_{\widetilde{\mathbf{L}}^4}
			\Vert\nabla \vartheta_{n}^{(m)}\Vert_{\mathbf{L}^2(\Omega)}^{\frac{1}{2}}
			\\
			&\leq\frac{\mu}{2}\Vert \Delta\vartheta_{n}^{(m)}\Vert_{L^2(\Omega)}^2
			+\frac{C}{\mu} \,  \underset{t\in(0,T_n)}{\rm sup} \, \Vert  \mathbf{v}_n^{(m)}(\cdot,t)\Vert_{\widetilde{\mathbf{L}}^4}^4 \Vert\nabla\vartheta_{n}^{(m)}\Vert_{\mathbf{L}^2(\Omega)}^2,
		\end{split}
		\\
		\label{est:J4}
		\begin{split}
			&\big|\psi_{n}^{(m)}(t) (\nabla\eta,\nabla\vartheta_{n}^{(m)} )\big| \leq 
			|\psi_{n}^{(m)}(t)|\left\Vert\nabla\eta\right\Vert_{\mathbf{L}^2(\Omega)} \Vert\nabla\vartheta_{n}^{(m)}\Vert_{\mathbf{L}^2(\Omega)}
			\\
			&\leq\frac{|\psi_n^{(m)}(t)|}2\Big(\Vert\nabla\vartheta_{n}^{(m)}\Vert_{\mathbf{L}^2(\Omega)}^2+\Vert\nabla\eta\Vert_{\mathbf{L}^2(\Omega)}^2\Big).
		\end{split}
	\end{align}
	Substituting the above estimates in \eqref{en:eqT2}, we obtain
	\begin{align}\label{Gr2}
		& \frac{d}{dt}\Vert \nabla\vartheta_{n}^{(m)}(\cdot,t)\Vert_{\mathbf{L}^2(\Omega)}^{2}+\mu \Vert \Delta\vartheta_{n}^{(m)}(\cdot,t)\Vert_{L^2(\Omega)}^{2}
		\nonumber 
		\\
		&\leq \Big(\frac{C}{\mu}
		\, \underset{t\in(0,T_n)}{\rm sup} \, \Vert  \mathbf{v}_n^{(m)}(\cdot,t)\Vert_{\widetilde{\mathbf{L}}^4}^4 +|\psi_n^{(m)}(t)|\Big)\Vert\nabla\vartheta_{n}^{(m)}(\cdot,t)\Vert_{\mathbf{L}^2(\Omega)}^2 \nonumber\\&\qquad+ |\psi_n^{(m)}(\cdot,t)|\Vert\nabla\eta\Vert_{\mathbf{L}^2(\Omega)}^2,
	\end{align}
	for all $t\in [0,T_n]$. An application of Gr\"onwall's inequality in \eqref{Gr2}, and then using the uniform estimates \eqref{fir-est1b} and \eqref{sec-est13++} in the resulting estimate, we obtain the uniform estimate \eqref{sec-est20}.
	
	Ladyzhenskaya's inequality gives the following estimate:
	\begin{align*}
		\underset{t\in(0,T_n)}{\rm sup} \, \left\Vert\vartheta_{n}^{(m)}(\cdot,t)\right\Vert_{L^4(\Omega)} & \leq  \, \,  \underset{t\in(0,T_n)}{\rm sup} \, \Big(2^\frac{1}{4} \left\Vert\vartheta_{n}^{(m)}(\cdot,t)\right\Vert_{L^2(\Omega)}^\frac12 \left\Vert \nabla \vartheta_{n}^{(m)}(\cdot,t)\right\Vert_{\mathbf{L}^2(\Omega)}^\frac12 \Big)
		\\
		& 
		\leq 2^\frac{1}{4}  \, \underset{t\in(0,T_n)}{\rm sup} \, \left\Vert\vartheta_{n}^{(m)}(\cdot,t)\right\Vert_{L^2(\Omega)}^\frac12 \, \underset{t\in(0,T_n)}{\rm sup} \, \left\Vert \nabla\vartheta_{n}^{(m)}(\cdot,t)\right\Vert_{\mathbf{L}^2(\Omega)}^\frac12 \leq C,
	\end{align*}
	which completes the proof of the lemma.
\end{proof}

\subsection{Passage to the limit as
	\texorpdfstring{$m\to\infty$}{m tends to infinity}}\label{Subs:PTL}
By the Banach--Alaoglu theorem and reflexivity, and  the energy estimates obtained in the Lemmas \ref{1Globest:1}--\ref{lem:estT21},  we can extract subsequences of $\{{\mathbf{v}}_n^{(m)}\}$, $\{{\vartheta}_n^{(m)}\}$, $\{{\varphi}_n^{(m)}\}$, and $\{\psi_n^{(m)}\}$, still labeled by $m$, such that the following convergences hold:
\begin{subequations}
	\begin{align}
		\label{wc:conv:L2p}
		\mathbf{v}_{n}^{(m)} \xrightharpoonup[m\to\infty]{}  \mathbf{v}_n\quad&\text{weakly in}\ L^{2}(0,T_n;\mathbf{V} \cap \mathbf{V}^2), 
		\\
		\label{wc:convT:L2p}
		\vartheta_{n}^{(m)} \xrightharpoonup[m\to\infty]{} \vartheta_n\quad&\text{weakly in}\ L^{2}(0,T_n;{W_0^{1,2}}(\Omega)\cap{W^{2,2}}(\Omega)),
		\\
		\label{wc:*conv}
		\mathbf{v}_{n}^{(m)} \xrightharpoonup[m\to\infty]{*}\mathbf{v}_n\quad&\text{weakly$^\ast$ in}\ 
		L^{\infty}(0,T_n;{\mathbf{V}}),
		\\
		\label{wc:convT:*W2}
		\vartheta_{n}^{(m)} \xrightharpoonup[m\to\infty]{*} \vartheta_{n}\quad&\text{weakly$^\ast$ in }L^{\infty}(0,T_n;
		{W_0^{1,2}}(\Omega)),
		\\
		\label{wc:convf:L2}
		\varphi_n^{(m)} \xrightharpoonup[m\to\infty]{} \varphi_n\quad&\text{weakly in }L^{2}(0,T_n),
		\\
		\label{wc:convj:L2}
		\psi_n^{(m)} \xrightharpoonup[m\to\infty]{} \psi_n \quad&\text{weakly in}\ L^{2}(0,T_n),
		\\
		\label{wc:conv:Gv-0'}
		\partial_t \mathbf{v}_{n}^{(m)} \xrightharpoonup[m\to\infty]{} \partial_t\mathbf{v}_{n}\quad&\text{weakly in}\
		L^{2}(0,T_n; {\mathbf{V}^\prime}),
		\\
		\label{wc:convTt:L2}
		\partial_t\vartheta_{n}^{(m)} \xrightharpoonup[m\to\infty]{} \partial_t\vartheta_{n}\quad &\text{weakly in}\ L^{2}(0,T_n;
		W^{-1,2}(\Omega)),
	\end{align}
\end{subequations}
where $\varphi_n(t;[\mathbf{v}_n,\vartheta_n])$ and $\psi_n(t;[\mathbf{v}_n,\vartheta_n])$ are the functionals defined in \eqref{00:eq:ht} and \eqref{00:eq:jt}, respectively.

On the other hand, due to the compact embeddings $\mathbf{V} \hookrightarrow\hookrightarrow\mathbf{H}$ and 
$W_0^{1,2}(\Omega)\hookrightarrow\hookrightarrow L^2(\Omega)$, and the Aubin--Lions--Simon compactness lemma, see Lemma \ref{ALSlemma}, it follows that
\begin{subequations}
	\begin{align}
		\label{str:conv:v}
		&\mathbf{v}_{n}^{(m)}\longrightarrow\mathbf{v}_n\quad\mbox{strongly in }{L}^2(0,T_n;\mathbf{H}) \quad \text{and a.e.} \ \ Q_{T_n} \ \  \mbox{as}\ m\longrightarrow \infty,
		\\
		\label{str:conv:T}
		&\vartheta_{n}^{(m)}\longrightarrow\vartheta_n\quad\mbox{strongly in } L^2(0,T_n;{L}^2(\Omega))\quad \text{and a.e.} \ \ Q_{T_n} \ \  \mbox{as}\ m\longrightarrow \infty.
	\end{align}
\end{subequations}
Using  \eqref{wc:conv:L2p}-\eqref{str:conv:T}, passing to limit in (\ref{eq:pr2}) and (\ref{eq:pr22}) standard way, we see that limit functions $(\mathbf{v}_n,\vartheta_n)$ are weak solutions to the direct problem \eqref{eq:1+}, and by Lemma~\ref{le:1}, they together with the limit functions $\varphi_n(t)=\varphi_n(t;[\mathbf{v}_n,\vartheta_n])$ and $\psi_n(t)=\psi_n(t;[\mathbf{v}_n,\vartheta_n])$ give the weak solutions to the inverse problem \eqref{eq:1}.

Let $\xi_1,$ $\xi_2\in C^\infty_0([0,T_n])$. Multiplying the equation (\ref{eq:pr2}) by $\xi_1(t)$ and equation (\ref{eq:pr22}) by $\xi_2(t)$, and integrating the resulting equations between $0$ and $T_n$, we obtain
\begin{align}\label{eq:lim1:nz}
	\begin{split}
		&
		\int_{Q_{T_n}} \partial_t \mathbf{v}_{n}^{(m)}\cdot\boldsymbol{\phi}_k \xi_1\,dxdt
		+
		\int_{Q_{T_n}}(\mathbf{v}_n^{(m)}\cdot{\nabla})\mathbf{v}_n\cdot\boldsymbol{\phi}_k \xi_1 \,dxdt
		+\nu\int_{Q_{T_n}}{\nabla}\mathbf{v}_{n}^{(m)}:{\nabla}\boldsymbol{\phi}_k \xi_1 \,dxdt
		\\
		&\quad=\int_{Q_{T_n}} \vartheta_{n}^{(m)} \mathbf{g}(x,t)\cdot \boldsymbol{\phi}_k
		\xi_1 \,dxdt+ \int_0^{T_n}
		\varphi_n^{(m)}(t)\int_{\Omega}\boldsymbol{\omega}(x) \cdot \boldsymbol{\phi}_k \xi_1\,dxdt,
	\end{split}
	\\
	\label{eq:lim2:nz}
	\begin{split}
		&
		\int_{Q_{T_n}}\partial_t\vartheta_{n}^{(m)}\zeta_k\xi_2\,dxdt +
		\int_{Q_{T_n}}(\mathbf{v}_n^{(m)}\cdot{\nabla})\vartheta_n^{(m)}\zeta_k\xi_2\,dxdt
		+\mu\int_{Q_{T_n}}{\nabla}\vartheta_n^{(m)}\cdot{\nabla}\zeta_k\xi_2\,dxdt
		\\
		&\quad
		=\int_0^{T_n}\psi_n^{(m)}(t)\int_{\Omega}\eta(x)\zeta_k\xi_2\,dxdt,
	\end{split}
\end{align}
for $k\in\{1,\dots,m\}$.

Fixing $k$, we can pass in equations (\ref{eq:lim1:nz}) and \eqref{eq:lim2:nz} to the limit $m\to\infty$ by using the convergence results (\ref{wc:conv:L2p})-(\ref{str:conv:T}). Then, we obtain
\begin{align}
	\label{eq:lim1:z}
	\begin{split}
		&
		\int_{Q_{T_n}} \partial_t \mathbf{v}_{n}\cdot\boldsymbol{\phi}_k\xi_1\,dxdt +
		\int_{Q_{T_n}}(\mathbf{v}_n\cdot{\nabla})\mathbf{v}_n\cdot\boldsymbol{\phi}_k\xi_1\,dxdt
		+
		\nu\int_{Q_{T_n}}{\nabla}\mathbf{v}_n:{\nabla}\boldsymbol{\phi}_k\xi_1\,dxdt
		\\
		& 
		=\int_{Q_{T_n}} \vartheta_n \mathbf{g}(x,t)\cdot  \boldsymbol{\phi}_k
		\xi_1\,dxdt+ \int_0^{T_n}\varphi_n(t)\int_{\Omega}\boldsymbol{\omega}(x) \cdot \boldsymbol{\phi}_k\xi_1\,dxdt,
	\end{split}
	\\
	\label{eq:lim2:z}
	\begin{split}
		&\int_{Q_{T_n}} \partial_t \vartheta_{n}\zeta_k\xi_2\,dxdt +
		\int_{Q_{T_n}}(\mathbf{v}_n\cdot{\nabla})\vartheta_n\zeta_k\xi_2\,dxdt
		+\mu \int_{Q_{T_n}}{\nabla}\vartheta_n\cdot{\nabla}\zeta_k\xi_2\,dxdt
		\\
		& 
		=\int_0^{T_n}
		\psi_n(t)\int_{\Omega}\eta(x)\zeta_k\xi_2\,dxdt,
	\end{split}
\end{align}
for $k\in\{1,\dots,m\}$.
{ Here, for the convective terms, we passed to the limit using the following convergence:
	\begin{subequations}
		\begin{align}\label{wc:rho:cnv}
			(\mathbf{v}_n^{(m)}\cdot{\nabla})\mathbf{v}_n^{(m)} \xrightharpoonup[m\to\infty]{}  (\mathbf{v}_n\cdot{\nabla})\mathbf{v}_n\quad&\text{weakly in}\ \ 
			{L}^{2}(0,T_n;\mathbf{V}^2), 
			\\
			\label{wc:T:cnv}
			(\mathbf{v}_n^{(m)}\cdot{\nabla})\vartheta_n^{(m)} \xrightharpoonup[m\to\infty]{}  (\mathbf{v}_n\cdot{\nabla})
			\vartheta_n\quad&\text{weakly in }L^{2}(0,T_n;W_0^{2,2}(\Omega)), 
		\end{align}
	\end{subequations}
	which can be proved using the strong convergences \eqref{str:conv:v}, \eqref{str:conv:T}, and Remark \ref{Remark 3.5}. In fact, writing the corresponding integral in \eqref{wc:rho:cnv} as
	\begin{align*}
		& \int_{Q_{T_n}}\big((\mathbf{v}_n^{(m)}\cdot{\nabla})\mathbf{v}_n^{(m)}-(\mathbf{v}_n\cdot{\nabla})\mathbf{v}_n\big)
		\cdot\boldsymbol{\phi}\,dxdt= 
		-\int_{Q_{T_n}}\big((\mathbf{v}_n^{(m)}-\mathbf{v}_n)\cdot{\nabla}\big)\boldsymbol{\phi}\cdot\mathbf{v}_n^{(m)}\,dxdt 
		\\
		&\qquad-\int_{Q_{T_n}}\big(\mathbf{v}_n\cdot{\nabla}\big)\boldsymbol{\phi}\cdot\big(\mathbf{v}_n^{(m)}-\mathbf{v}_n\big)\,dxdt,
	\end{align*}
	We see that the  right-hand side integrals converge to zero by an application of H\"older's inequality, together with Remark \ref{Remark 3.5} and the convergence (\ref{str:conv:v}): 
	\begin{align*}
		&\left\vert\int_{Q_{T_n}}\big((\mathbf{v}_n^{(m)}-\mathbf{v}_n)\cdot{\nabla}\big)\boldsymbol{\phi}\cdot \mathbf{v}_n^{(m)}\,dxdt\right\vert
		\leq
		\|\mathbf{v}_n^{(m)}-\mathbf{v}_n\|_{{L}^2(Q_{T_n})}\|\mathbf{v}_{n}^{(m)}\|_{{L}^4(Q_{T_n})} 
		\|\nabla\boldsymbol{\phi}\|_{{L}^4(Q_{T_n})}
		\longrightarrow0,
		\\
		&
		\left\vert\int_{Q_{T_n}}(\mathbf{v}_n\cdot{\nabla})\boldsymbol{\phi} \cdot(\mathbf{v}_n^{(m)}-\mathbf{v}_n)\,dxdt\right\vert
		\leq
		\|\mathbf{v}_n^{(m)}-\mathbf{v}_n\|_{{L}^2(Q_{T_n})}\|\mathbf{v}_{n}\|_{{L}^4(Q_{T_n})} 
		\|\nabla\boldsymbol{\phi}\|_{{L}^4(Q_{T_n})}
		\longrightarrow0,
	\end{align*}
	as $m \to \infty$ and $\boldsymbol{\phi} \in {L}^{2}(0,T_n;\mathbf{V}^2)$.
	
	Analogically, writing the corresponding integral in  \eqref{wc:T:cnv} as
	\begin{align*}
		& \int_{Q_{T_n}}\big((\mathbf{v}_n^{(m)}\cdot{\nabla})\vartheta_n^{(m)}-(\mathbf{v}_n\cdot{\nabla})\vartheta_n\big)
		\zeta\,dxdt= 
		-\int_{Q_{T_n}}\big((\mathbf{v}_n^{(m)}-\mathbf{v}_n)\cdot{\nabla}\big)\zeta\vartheta_n^{(m)}\,dxdt 
		\\
		&\qquad-\int_{Q_{T_n}}(\mathbf{v}_n\cdot{\nabla})\zeta\,\big(\vartheta_n^{(m)}-\vartheta_n \big)\,dxdt,
	\end{align*}
	and the  right-hand side integrals converge to zero using H\"older's inequality, the convergence (\ref{str:conv:v}), \eqref{str:conv:T}, and Remark \ref{Remark 3.5}: 
	\begin{align*}
		\\
		&
		\left\vert\int_{Q_{T_n}}\big((\mathbf{v}_{n}^{(m)}-\mathbf{v}_n)\cdot{\nabla}\big)\zeta\vartheta_n^{(m)}\,dxdt
		\right\vert
		\leq
		\|\mathbf{v}_n^{(m)}-\mathbf{v}_n\|_{{L}^2(Q_{T_n})}\|\vartheta_{n}^{(m)}\|_{L^4(Q_{T_n})} 
		\|\nabla\zeta\|_{{L}^4(Q_{T_n})}\longrightarrow0,
		\\
		&\left\vert\int_{Q_{T_n}}(\mathbf{v}_n\cdot{\nabla})\zeta\,(\vartheta_n^{(m)}-\vartheta_n)\,dxdt\right\vert
		\leq
		\|\vartheta_n^{(m)}-\vartheta_n\|_{L^2(Q_{T_n})}
		\|\mathbf{v}_n\|_{{L}^4(Q_{T_n})}
		\|\nabla\zeta\|_{L^4(Q_{T_n})}\longrightarrow0,
	\end{align*}
	as $m \to \infty$ and ${\zeta} \in {L}^{2}(0,T_n;W_0^{2,2}(\Omega))$.}

By linearity, the equations (\ref{eq:lim1:nz}) and (\ref{eq:lim2:nz}) hold for any finite linear combination of $\boldsymbol{\phi}_1,\dots,\boldsymbol{\phi}_m$ and ${\zeta}_1,\dots,{\zeta}_m$, respectively, and, by a continuity argument, they are still true for any { $\boldsymbol{\phi}\xi_1 \in {L}^2(0,T_n;\mathbf{V}^2)$ and $\zeta\xi_2\in {L}^2(0,T_n; W_0^{2,2}(\Omega))$ }with $\xi_1$, $\xi_2\in C^\infty_0 (0,T_n)$, respectively. Moreover, all terms in the equations (\ref{eq:lim1:nz}) and (\ref{eq:lim2:nz}) are absolutely continuous as functions of $t$ defined by integrals over $[0, T_0]$. So we obtain the following equalities, which hold for a.e. $t\in [0, T_0]$ and for any {$\boldsymbol{\phi}\in\mathbf{V}$ and $\zeta\in W_0^{1,2}(\Omega)$,} respectively:
\begin{align}
	\label{eq:lim:u}
	\begin{split}
		&
		\int_{\Omega}\partial_t \mathbf{v}_{n}(t) \cdot \boldsymbol{\phi}\,dx + \int_{\Omega}\big(\mathbf{v}_n(t)\cdot{\nabla}\big)\mathbf{v}_n(t)\cdot\boldsymbol{\phi}\,dx+
		\nu \int_{\Omega}{\nabla}\mathbf{v}_n(t):{\nabla}\boldsymbol{\phi}\,dx
		\\
		&
		\quad=\int_{\Omega}  \vartheta_n(t) \mathbf{g}(x,t)
		\cdot \boldsymbol{\phi}dx
		+\varphi_n(t)\int_{\Omega}\boldsymbol{\omega}(x)\cdot\boldsymbol{\phi}\,dx,
	\end{split}
	\\
	\label{eq:weak-form:u:nc}
	\begin{split}
		& \int_{\Omega} \partial_t\vartheta_{n}(t)\zeta\,dx +
		\int_{\Omega}(\mathbf{v}_n(t)\cdot{\nabla})\vartheta_n(t)\zeta\,dx
		+\mu\int_{\Omega}{\nabla}\vartheta_n(t) \cdot{\nabla}\zeta\,dx
		\\
		&\quad=\psi_n(t)\int_{\Omega}\eta(x)\zeta\,dx.
	\end{split}
\end{align}
Thus, the pair of limit functions $(\mathbf{v}_n,\vartheta_n)$ is the weak regular solution to the direct problem \eqref{eq:1+}, and by Lemma~\ref{le:1}, it gives the weak regular solution to the inverse problem \eqref{eq:1}.

Furthermore, due to the weakly lower semicontinuity of norms, we obtain the following estimate from (\ref{fir-est1}) and \eqref{sec-est13}.
\begin{subequations}
	\label{eq:apr1:lim}
	\begin{equation}
		\begin{split}
			&\left\Vert \mathbf{v}_n\right\Vert_{{L}^{\infty}(0,T_n; \mathbf{V})}+\left\Vert \vartheta_n\right\Vert_{{L}^{\infty}(0,T_n; W_0^{1,2}(\Omega))}+\left\Vert \mathbf{v}_n\right\Vert_{{L}^{2}(0,T_n; \mathbf{V})}+\left\Vert \vartheta_n\right\Vert_{L^2(0,T_n;W_0^{1,2}(\Omega))}
			\\
			&\quad+\left\Vert  \partial_t \mathbf{v}_{n}\right\Vert_{{L}^{1}(0,T_n; \mathbf{V}^\prime)}+\left\Vert \partial_t \vartheta_n\right\Vert_{L^1\left(0,T_n;W^{-1,2}(\Omega)\right)}+\left\Vert \Delta \mathbf{v}_n\right\Vert_{L^2(0,T_n;\mathbf{H})}\\&\qquad+\left\Vert \Delta \vartheta_n\right\Vert_{L^2(0,T_n;L^{2}(\Omega))}
			\leq C<\infty,
		\end{split}
	\end{equation}
	and it follows from \eqref{fir-est1b} that
	\begin{align}
		\begin{split}
			&
			\int_0^{T_n}\left|\varphi_n(t)\right|\,dt\leq \int_0^{T_n}|f_n'(t)|\,dt
			+\nu\int_0^{T_n}\left\Vert\nabla \mathbf{v}_{n}\right\Vert_{\mathbf{H}}\left\Vert
			\nabla\boldsymbol\omega\right\Vert_{\mathbf{H}}\,dt
			\\
			&\quad+g_0 T_n\sup\limits_{t\in(0,T_n)}\left\Vert
			\vartheta_{n}(t)\right\Vert_{L^2(\Omega)}\left\Vert \boldsymbol\omega\right\Vert_{\mathbf{H}}
			+\int_0^{T_n}\left\Vert \mathbf{v}_{n}\right\Vert_{\mathbf{H}}
			\left\Vert\nabla\mathbf{v}_{n}\right\Vert_{\mathbf{H}}
			\left\Vert\boldsymbol\omega\right\Vert _{\widetilde{\mathbf{L}}^\infty} \, d t,
		\end{split}
		\\
		%\label{est:jn}
		\begin{split}
			&\int_0^{T_n}\left|\psi_n(t)\right|\,dt\leq \int_0^{T_n}|g_n'(t)|\,dt+\mu \int_0^{T_n}\left\Vert \nabla \vartheta_n\right\Vert_{\mathbf{L}^2(\Omega)}\left\Vert \nabla\eta\right\Vert _{\mathbf{L}^2(\Omega)}\,dt
			\\
			&\quad+\int_0^{T_n}\left\Vert\mathbf{v}_n\right\Vert_{\mathbf{H}}
			\left\Vert\nabla\vartheta_n\right\Vert_{\mathbf{L}^2(\Omega)}\left\Vert\eta\right\Vert _{{L}^\infty(\Omega)}\,dt,
		\end{split}
	\end{align}
\end{subequations}
respectively. The set of estimates \eqref{eq:apr1:lim} gives \eqref{estim-1}.

\subsection{\bf Reconstruction of Pressure}
Define 
\begin{align*}
	\mathbf{G}_n=\mathbf{g}(x,t) \vartheta_n+ \varphi_n(t)\boldsymbol{\omega}(x)  - \partial_t \mathbf{v}_{n}
	-\left(\mathbf{v}_n\cdot \nabla \right) \mathbf{v}_n+\nu \Delta\mathbf{v}_n.
\end{align*}
Then from \eqref{eq:1+a}, we have 
\begin{align*}
	\mathbf{G}_n = \nabla \pi_n.
\end{align*}
Using the Ladyzhenskaya inequality \eqref{SE:case:p*>4} and the estimates \eqref{fir-est1a}, we obtain
\begin{align*}
	\int_0^{T_n} \| \left(\mathbf{v}_n\cdot \nabla \right) \mathbf{v}_n \|_{\mathbf{V}^\prime} \, d t &\leq C \int_0^{T_n}  \| \mathbf{v}_n \|_{\widetilde{\mathbf{L}}^4}^2 \, d t \leq C \int_0^{T_n}  \| \mathbf{v}_n \|_{\mathbf{H}} \| \mathbf{v}_n \|_{\mathbf{V}} \, d t
	\\& \leq C \sqrt{T_n}  \| \mathbf{v}_n \|_{L^\infty(0,T_n;\mathbf{H})} \|  \mathbf{v}_n \|_{L^2(0,T_n;\mathbf{V})} \leq C.
\end{align*}
Moreover, we have
\begin{align*}
	&\int_0^{T_n} \| \Delta \mathbf{v}_n \|_{\mathbf{V}^\prime} \, d t  \leq C  \int_0^{T_n} \|  \mathbf{v}_n \|_{\mathbf{V}} \, d t \leq C \sqrt{T_n} \|  \mathbf{v}_n \|_{L^2(0,T_n;\mathbf{V})} \leq C,\\
	& \int_0^{T_n} \|\varphi_n  \boldsymbol{\omega}\|_{\mathbf{V}^\prime} \, d t \leq    C \| \boldsymbol{\omega}\|_{\mathbf{H}} \int_0^{T_n} |\varphi_n(t) | \, d t \leq C,
	\\&
	\int_0^{T_n} \| \mathbf{g} \vartheta_n \|_{{H}^{-1}(\Omega)} \, d t \leq C g_0  \int_0^{T_n} \|  \vartheta_n \|_{{L}^{2}(\Omega)} \, d t
	\leq C g_0  {T_n} \|  \vartheta_n \|_{L^\infty(0,T_n;{L}^{2}(\Omega))} \leq C,
\end{align*}
and $\partial_t \mathbf{v}_n \in L^1(0,T_n;\mathbf{V}^\prime)$. Hence, we have 
\begin{align*}
	\mathbf{G}_n \in {L^1(0,T_n;\mathbf{V}^\prime)}.
\end{align*}
For every divergence-free test function $\boldsymbol{\phi} \in \mathbf{V}$, the weak formulation of \eqref{eq:1+a} yields
\begin{align*}
	\langle \mathbf{G}_n, \boldsymbol{\phi} \rangle=0.
\end{align*}
Therefore, by the de Rham theorem, there exists a scalar distribution $\pi_n$, unique up to an additive function of time, such that
\begin{align*}
	\mathbf{G}_n= \nabla \pi_n.
\end{align*}
Therefore, 
\begin{align*}
	\nabla \pi_n =\mathbf{g}(x,t) \vartheta_n+ \varphi_n(t)\boldsymbol{\omega}(x)  - \partial_t \mathbf{v}_{n}
	-\left(\mathbf{v}_n\cdot \nabla \right) \mathbf{v}_n+\nu \Delta\mathbf{v}_n.
\end{align*}
Since pressure is determined only up to an additive function of time, we normalize it by $\int_{\Omega} \pi_n(x,t) \, d x=0$. Then $\pi_n$ is uniquely determine.

Using the pressure estimate, we have
\begin{align*}
	\| \nabla \pi_n(\cdot, t) \|_{\mathbf{H}^{-1}(\Omega)} \leq C  \| \mathbf{G}_n (\cdot, t) \|_{\mathbf{H}^{-1}(\Omega)}, 
\end{align*}
and we obtain
\begin{align*}
	\| \pi_n \|_{L^1(0,T_n;L_0^2(\Omega))} \leq C  \| \mathbf{G}_n  \|_{L^1(0,T_n;\mathbf{H}^{-1}(\Omega))} \leq C.
\end{align*}
where $L_0^2(\Omega):= \Big\{ \pi \in L^2(\Omega) : \int_{\Omega} \pi(x) \, d x =0 \Big\}$.  Hence,  $\pi_n \in {L^1(0,T_n;L_0^2(\Omega))} $.

\subsection{\bf Uniqueness}
Let us first prove the uniqueness of weak solutions. Let $(\mathbf{v}_n^1,\vartheta_n^1,\pi_n^1)$ and $(\mathbf{v}_n^2,\vartheta_n^2,\pi_n^2)$ be two weak solutions of the system \eqref{eq:1+} satisfying \eqref{estim-1a+} and \eqref{estim-1b+}. Then $\mathbf{v}_n:=\mathbf{v}_n^1-\mathbf{v}_n^2$, $\vartheta_n:=\vartheta_n^1-\vartheta_n^2$, and $\pi_n:=\pi_n^1-\pi_n^2$ satisfy the following energy integral:
\begin{align}\label{U1}
	&\frac{1}{2}\frac{d}{dt} \Big(\int_ \Omega|\mathbf{v}_n|^2\,d {x}+\int_\Omega|\vartheta_n|^2\,d x\Big)+ \nu\int_{\Omega}|\nabla\mathbf{v}_n|^2\,d x + \mu\int_{\Omega}|\nabla\vartheta_n|^2\,d {x} \nonumber\\&=-\int_{\Omega} \big((\mathbf{v}_n \cdot\nabla)\mathbf{v}_n^1\big) \cdot \mathbf{v}_n \,d{x} -\int_{\Omega} \big((\mathbf{v}_n \cdot\nabla)\vartheta_n^1\big)  \vartheta_n \,d{x} +\int_{\Omega}  \vartheta_n \mathbf{g} \cdot \mathbf{v}_n \,d{x}
	\nonumber\\&\quad+\Big(\varphi_n^1(t;[\mathbf{v}_n^{1}, \vartheta_n^{1}])-\varphi_n^2(t;[\mathbf{v}_n^{2}, \vartheta_n^{2}])\Big)
	\int_\Omega\boldsymbol{\omega}({x}) \cdot\mathbf{v}_n\,d{x}\nonumber\\&\quad+\Big(\psi_n^1(t;[\mathbf{v}_n^{1}, \vartheta_n^{1}])-\psi_n^2(t;[\mathbf{v}_n^{2}, \vartheta_n^{2}])\Big)
	\int_\Omega \eta({x}) \vartheta_n\,d{x}.
\end{align}
Using H\"older's, Ladyzhenskaya's, and Young's inequalities, we obtain 
\begin{align*}
	\Big|\int_{\Omega} \big((\mathbf{v}_n \cdot\nabla)\mathbf{v}_n^1\big) \cdot \mathbf{v}_n \,d{x} \Big| &\leq \|\nabla \mathbf{v}_n^1\|_{\mathbf{H}} \|\mathbf{v}_n\|_{\widetilde{\mathbf{L}}^4}^2 \leq \sqrt{2} \|\nabla \mathbf{v}_n^1\|_{\mathbf{H}} \|\mathbf{v}_n\|_{\mathbf{H}}\|\nabla\mathbf{v}_n\|_{\mathbb{H}}\nonumber
	\\
	&   \leq \frac{\nu}{4} \|\nabla\mathbf{v}_n\|_{\mathbf{H}}^2+ \frac{2}{\nu}\|\nabla \mathbf{v}_n^1\|_{\mathbf{H}}^2 \|\mathbf{v}_n\|_{\mathbf{H}}^2,
\end{align*}
and
\begin{align*}
	\Big| \int_{\Omega} \big((\mathbf{v}_n \cdot\nabla)\vartheta_n^1\big)  \vartheta_n \,d{x} \Big|&=\Big| \int_{\Omega} \big((\mathbf{v}_n \cdot\nabla)\vartheta_n\big)  \vartheta_n^1 \,d{x} \Big| \leq \|\mathbf{v}_n\|_{\widetilde{\mathbf{L}}^4} \|\nabla \vartheta_n\|_{\mathbf{L}^2(\Omega)} \|\vartheta_n^1\|_{L^4(\Omega)} \\&\leq \sqrt{2} \| \mathbf{v}_n\|_{\mathbf{H}}^\frac{1}{2} \|\nabla\mathbf{v}_n\|_{\mathbf{H}}^\frac{1}{2}\|\nabla \vartheta_n\|_{\mathbf{L}^2(\Omega)} \|\vartheta_n^1\|_{L^4(\Omega)} 
	\\&\leq \frac{\mu}{2} \|\nabla \vartheta_n\|_{\mathbf{L}^2(\Omega)}^2+\frac{1}{\mu} \| \mathbf{v}_n\|_{\mathbf{H}}\|\nabla\mathbf{v}_n\|_{\mathbf{H}} \|\vartheta_n^1\|_{L^4(\Omega)}^2 
	\nonumber
	\\&   \leq \frac{\nu}{4} \|\nabla\mathbf{v}_n\|_{\mathbf{H}}^2+ \frac{\mu}{2} \|\nabla \vartheta_n\|_{\mathbf{L}^2(\Omega)}^2+\frac{1}{\mu^2 \nu} \| \mathbf{v}_n\|_{\mathbf{H}}^2\|\vartheta_n^1\|_{L^4(\Omega)}^4.
\end{align*}
We also have
\begin{align*}
	\Big| \int_{\Omega}  \vartheta_n \mathbf{g} \cdot \mathbf{v}_n \,d{x} \Big| \leq \frac{g_0}{2}\Big( \|\vartheta_n\|_{{L}^2(\Omega)}^2+ \| \mathbf{v}_n\|_{\mathbf{H}}^2 \Big).
\end{align*}
Recalling the definitions of the functionals \eqref{00:eq:ht} and \eqref{00:eq:jt} and the over-determination conditions \eqref{eq:4}, the terms, which include $\varphi_n^1-\varphi_n^2$ and $\psi_n^1-\psi_n^2$, vanish in identity \eqref{U1}. Indeed, 
\begin{align*}
	&(\varphi_n^1-\varphi_n^2)\int_\Omega \boldsymbol{\omega}({x}) \cdot\mathbf{v}_n({x},t)\,d {x}=
	(\varphi_n^1-\varphi_n^2)(f_n(t)-f_n(t))=0,
	\\&(\psi_n^1-\psi_n^2)\int_\Omega \eta({x}) \vartheta_n({x},t)\,d {x}=
	(\psi_n^1-\psi_n^2)(g_n(t)-g_n(t))=0.
\end{align*}
Substituting the above estimates in \eqref{U1} and then integrating it from $0$ to $t$, we get
\begin{align}\label{U4}
	& \|\mathbf{v}_n(\cdot,t)\|_{\mathbf{H}}^2+\|\vartheta_n(\cdot,t)\|_{L^2(\Omega)}^2+ \nu\int_0^t \|\nabla \mathbf{v}_n(\cdot,s)\|_{\mathbf{H}}^2 \,d s + \mu\int_0^t \|\nabla \vartheta_n(\cdot,s)\|_{\mathbf{L}^2(\Omega)}^2 \,d s \nonumber
	\\
	&\leq \|\mathbf{v}_n(\cdot,0)\|_{\mathbf{H}}^2 +\|\vartheta_n(\cdot,0)\|_{L^2(\Omega)}^2 + \int_0^t \Big(\frac{4}{\nu} \|\nabla \mathbf{v}_n^1(\cdot,s)\|_{\mathbf{H}}^2 +g_0+ \frac{2}{\mu^2 \nu} \|\vartheta_n^1(\cdot,s)\|_{L^4(\Omega)}^4 \Big)
	\nonumber \\& \qquad \times\Big(\|\mathbf{v}_n(\cdot,s)\|_{\mathbf{H}}^2 +\|\vartheta_n(\cdot,s)\|_{{L}^2(\Omega)}^2\Big)\,ds,
\end{align}
for all $t \in [0,T_n]$. Applying Gronwall's inequality in \eqref{U4}, we obtain
\begin{align}\label{U5}
	\|\mathbf{v}_n(\cdot,t)\|_{\mathbf{H}}^2 &+\|\vartheta_n(\cdot,t)\|_{L^2(\Omega)}^2
	\leq \|\mathbf{v}_n(\cdot,0)\|_{\mathbf{H}}^2 +\|\vartheta_n(\cdot,0)\|_{L^2(\Omega)}^2 \nonumber\\& \times \exp\left(\int_0^{T_n} \Big(\frac{4}{\nu} \|\nabla \mathbf{v}_n^1(\cdot,t)\|_{\mathbf{H}}^2 +g_0+ \frac{2}{\mu^2 \nu} \|\vartheta_n^1(\cdot,t)\|_{L^4(\Omega)}^4 \Big) \,d t \right).
\end{align}
Taking $\mathbf{v}_n^1({x},0)=\mathbf{v}_n^2({x},0)$ and $\vartheta_n^1({x},0)=\vartheta_n^2({x},0)$ in \eqref{U5}, we get $\mathbf{v}_n^1({x},t)=\mathbf{v}_n^2({x},t)$ and $\vartheta_n^1({x},t)=\vartheta_n^2({x},t)$, and hence the uniqueness of the weak solution follows. The uniqueness of weak regular solutions then follows from the uniqueness of weak solutions.

\section{\bf The main result}\label{sec4}\setcounter{equation}{0}
{ This section is devoted to stating and proving the main result of the paper. Auxiliary procedures, such as scaling and shifting, are employed to handle the infinitesimal initial layer and to establish the convergence of the approximate solutions. We then pass to the limit to obtain a solution to the original inverse problem.}

In the present paper, the sequences $\{f_n\}_{n\geq n_0}$ and $\{g_n\}_{n\geq n_0}$ approximate discontinuous functions $f(t)$ and $g(t)$ at $t=0$ on the interval $(0,\frac1n)$: 
\[
\mathrm{supp}(f_n-f)=\mathrm{supp}(g_n-g)=\Big[0,\dfrac1n\Big].
\]
{ Here, we represent the domain $Q_{T_n}$ as the union of two subdomains $\Omega\times[0,\frac1n]$ and $\Omega\times[\frac1n,T_n]$ with intersection at $t=\frac1n$, and formulate the scaled and shifted problems, respectively.} Such an approach was elaborated in \cite{AKS} for the impulsive Navier--Stokes system and in \cite{AKS-2024} for the impulsive Kelvin--Voigt equations and in \cite{Ant-2026} for the impulsive inverse problem for the 2D CBF equations. 

In the domain $\Omega\times[0,\frac1n]$, for $\overline{t}=tn\in[0,1]$, we introduce the scaled functions
\[
\begin{split}
	&\overline{\mathbf{v}}_n({x},{\overline{t}})=\mathbf{v}_n({x},\overline{t}/n), \ \ 
	\overline{\varphi}_n({\overline{t}})=\frac1n\varphi_n(\overline{t}/n), \ \
	\overline{f}_n({\overline{t}})=f_n(\overline{t}/n), \ \
	\overline{\pi}_n({x},{\overline{t}})=\frac1n\pi_n({x},\overline{t}/n),
	\\
	&\overline{\mathbf{g}}_n({x},\overline{t})=\mathbf{g}({x},\overline{t}/n), 
	\
	\overline{\vartheta}_n({x},{\overline{t}})=\vartheta_n({x},\overline{t}/n), 
	\  
	\overline{\psi}_n({\overline{t}})=\frac1n\psi_n(\overline{t}/n), 
	\
	\overline{g}_n({\overline{t}})=g_n(\overline{t}/n),
\end{split}
\]
satisfying the system
\begin{subequations}
	\label{eq3.1}
	\begin{align}
		\begin{split} 
			&\partial_{\overline{t}}\overline{\mathbf{v}}_{n}+
			\frac1n\left(\overline{\mathbf{v}}_n\cdot \nabla \right) \overline{\mathbf{v}}_{n}+\nabla\overline{\pi}_n
			=\frac{\nu}{n}\Delta\overline{\mathbf{v}}_n
			+\frac1n\overline{\mathbf{g}}_n\overline{\vartheta}_{n}+\overline{\varphi}_n({\overline{t}})\,\boldsymbol{\omega}({x}), \ \ \text{ in }Q_1,
		\end{split}
		\\
		&
		\partial_{\overline{t}}\overline{\vartheta}_{n}+
		\frac1n\left(\overline{\mathbf{v}}_n\cdot \nabla \right) \overline{\vartheta}_n
		= \frac\mu n\Delta \overline{\vartheta}_n+\overline{\psi}_n(\overline{t}) \eta(x),\ \ \ \text{ in } Q_1,
		\\
		&\operatorname{div}\overline{\mathbf{v}}_n ({x},{\overline{t}}) =0, \ \ \ \text{ in
		}Q_1,
		\\
		&\int_\Omega
		\overline{\mathbf{v}}_n({x},{\overline{t}})\cdot\boldsymbol{\omega}({x})\,d{x}=\overline{f}_n({\overline{t}}), \ \ \ \text{ in }(0,1),
		\\
		&\int_\Omega
		\overline{\vartheta}_n({x},{\overline{t}})\eta({x})\,d{x}=\overline{g}_n({\overline{t}}), \ \ \ \text{ in }(0,1),
		\\
		&\overline{\mathbf{v}}_n({x},0)=\mathbf{v}_0({x}),\quad \overline{\vartheta}_n({x},0)=\vartheta_0({x}), \ \ \ \text{ in }\Omega,
		\\
		&\overline{\mathbf{v}}_n({x},{\overline{t}})=\mathbf{0},\quad \overline{\vartheta}_n({x},{\overline{t}})=0, \ \ \ \text{ on }\Gamma_1,
	\end{align}
\end{subequations}
where $Q_1:=\Omega\times(0,1)$ and $\Gamma_1:=\partial\Omega\times(0,1)$. {This inverse problem is equivalent to the functional direct problem 
	\begin{subequations}
		\label{eq3.1d}
		\begin{align}
			\begin{split} &\partial_{\overline{t}}
				\overline{\mathbf{v}}_{n}+\frac1n\left(\overline{\mathbf{v}}_n\cdot \nabla \right) \overline{\mathbf{v}}_{n}+\nabla\overline{\pi}_n=
				\frac{\nu}{n}\Delta\overline{\mathbf{v}}_n  +\frac1n\overline{\mathbf{g}}_n \overline{\vartheta}_{n}+
				\overline{\varphi}_n({\overline{t}};[\overline{\mathbf{v}}_n,\overline{\vartheta}_n])\,
				\boldsymbol{\omega}({x}), \ \ \ \text{ in }Q_1,
			\end{split}
			\\
			&
			\partial_{\overline{t}}\overline{\vartheta}_{n}+\frac1n\left(\overline{\mathbf{v}}_n\cdot \nabla \right) \overline{\vartheta}_n
			= \frac\mu n\Delta \overline{\vartheta}_n+\overline{\psi}_n(\overline{t};[\overline{\mathbf{v}}_n,\overline{\vartheta}_n])\eta(x),\ \ \ \text{ in } Q_1,
			\\
			&\operatorname{div}\overline{\mathbf{v}}_n ({x},{\overline{t}})=0, \ \ \ \text{ in
			}Q_1,
			\\
			&\overline{\mathbf{v}}_n({x},0)=\mathbf{v}_0({x}),\quad \overline{\vartheta}_n({x},0)=\vartheta_0({x}), \ \ \ \text{ in }\Omega,
			\\
			&\overline{\mathbf{v}}_n({x},{\overline{t}})=\mathbf{0},\quad \overline{\vartheta}_n({x},{\overline{t}})=0, \ \ \ \text{ on }\Gamma_1,
		\end{align}
		where
		\begin{align}
			\begin{split}
				\overline{\varphi}_n({\overline{t}};[\overline{\mathbf{v}}_n,\overline{\vartheta}_n])&=
				\overline{f}_n'({\overline{t}})
				-\frac1n\int_\Omega
				(\overline{\mathbf{v}}_n({x},{\overline{t}})\cdot\nabla)\boldsymbol{\omega}({x})\cdot\overline{\mathbf{v}}_n({x},{\overline{t}})\,d{x}
				\\
				&\quad+\frac{\nu}{n}\int_\Omega\nabla\overline{\mathbf{v}}_n({x},{\overline{t}}):\nabla\boldsymbol{\omega}({x})\,d{x} -\frac1n\int_\Omega
				\overline{\mathbf{g}}_n(x,{\overline{t}})
				\cdot\boldsymbol{\omega}({x})\overline{\vartheta}_{n}(x,\overline{t})\,dx,
			\end{split}
			\label{eq3.1df}
			\\
			\overline{\psi}_n({\overline{t}};[\overline{\mathbf{v}}_n,\overline{\vartheta}_n])&=
			\overline{g}_n'({\overline{t}})
			+\frac{\mu}{n}\int_\Omega\nabla\overline{\vartheta}_n({x},{\overline{t}})\cdot\nabla\eta({x})\,d{x} 
			-\frac1n\int_\Omega
			(\overline{\mathbf{v}}_n({x},{\overline{t}})\cdot\nabla)\eta({x})\overline{\vartheta}_n({x},{\overline{t}})\,d{x}.
			\label{eq3.1dg}
		\end{align}
	\end{subequations}
	
	In the rest of the domain {$\Omega\times[\frac1n,T_n]$}, for $\widetilde{t}=t-\frac1n\in[0,T]$, we introduce the functions
	\[
	\begin{split}
		&\widetilde{\mathbf{v}}_n({x},\widetilde{t})=\mathbf{v}_n({x},\widetilde{t}+1/n), 
		\ \ 
		\widetilde{\varphi}_n({\widetilde{t}})=\varphi_n(\widetilde{t}+1/n), 
		\ \
		\widetilde{f}_n(\widetilde{t})=f_n(\widetilde{t}+1/n), 
		\\
		&
		\widetilde{\pi}_n({x},\widetilde{t})=\pi_n({x},\widetilde{t}+1/n),
		\ \
		\widetilde{\mathbf{g}}_n({x},\widetilde{t})=\mathbf{g}({x},\widetilde{t}+1/n),
		\ \ 
		\widetilde{\vartheta}_n({x},\widetilde{t})=\vartheta_n({x},\widetilde{t}+1/n), 
		\\
		&\widetilde{\psi}_n(\widetilde{t})=\psi_n(\widetilde{t}+1/n), 
		\ \
		\widetilde{g}_n(\widetilde{t})=g_n(\widetilde{t}+1/n),
	\end{split}
	\]
	satisfying the system
	\begin{subequations}
		\label{eq3.2}
		\begin{align}
			\begin{split} &\partial_{\widetilde{t}}\widetilde{\mathbf{v}}_{n}+
				\left(\widetilde{\mathbf{v}}_n\cdot \nabla \right) \widetilde{\mathbf{v}}_n
				+\nabla\widetilde{\pi}_n=
				\nu\Delta\widetilde{\mathbf{v}}_n+\widetilde{\mathbf{g}}_n
				\widetilde{\vartheta}_n+\widetilde{\varphi}_n(\widetilde{t})\,
				\boldsymbol{\omega}({x}), \ \ \ \text{ in }Q_{T},
			\end{split}
			\\
			&
			\partial_{\widetilde{t}}\widetilde{\vartheta}_{n}+\left(\widetilde{\mathbf{v}}_n\cdot \nabla \right) \widetilde{\vartheta}_n
			=\mu\Delta\widetilde{\vartheta}_n+\widetilde{\psi}_n(\widetilde{t})\eta(x),\ \ \ \text{ in }Q_{T},
			\\
			&\operatorname{div}\widetilde{\mathbf{v}}_n ({x},\widetilde{t}) =0, \ \ \ \text{ in
			}Q_{T},
			\\
			&\int_\Omega
			\widetilde{\mathbf{v}}_n({x},\widetilde{t})
			\cdot\boldsymbol{\omega}({x})\,d{x}
			=\widetilde{f}_n(\widetilde{t}), \ \ \ \text{ in }(0,T),
			\\
			&
			\int_\Omega
			\widetilde{\vartheta}_n({x},\widetilde{t})
			\eta({x})\,d{x}=\widetilde{g}_n(\widetilde{t}),
			\ \ \ \text{ in }(0,T),
			\\
			&\widetilde{\mathbf{v}}_n({x},0)=
			\overline{\mathbf{v}}_n({x},1), \ \ \ \widetilde{\vartheta}_n({x},0)=
			\overline{\vartheta}_n({x},1), \ \ \ \text{ in }\Omega,
			\\
			&\widetilde{\mathbf{v}}_n({x},\widetilde{t})=\mathbf{0}, \ \ \ \widetilde{\vartheta}_n({x},\widetilde{t})=0, \ \ \ \text{ on }\Gamma_{T},
		\end{align}
	\end{subequations}
	which is equivalent to the functional direct problem 
	\begin{subequations}
		\label{eq3.2d}
		\begin{align}
			\begin{split} &\partial_{\widetilde{t}}\widetilde{\mathbf{v}}_{n}+
				\left(\widetilde{\mathbf{v}}_n\cdot \nabla \right) \widetilde{\mathbf{v}}_n+\nabla\widetilde{\pi}_n=
				\nu\Delta\widetilde{\mathbf{v}}_n
				+\widetilde{\mathbf{g}}_n\widetilde{\vartheta}_n
				+
				\widetilde{\varphi}_n(\widetilde{t};[\widetilde{\mathbf{v}}_n,\widetilde{\vartheta}_n])\,
				\boldsymbol{\omega}({x}), \ \ \ \text{ in }Q_{T},
			\end{split}
			\\
			&
			\partial_{\widetilde{t}}\widetilde{\vartheta}_{n}+\left(\widetilde{\mathbf{v}}_n\cdot \nabla \right) \widetilde{\vartheta}_n
			=\mu\Delta\widetilde{\vartheta}_n+\widetilde{\psi}_n(\widetilde{t};[\widetilde{\mathbf{v}}_n,\widetilde{\vartheta}_n])\eta(x), \ \ \ \text{ in }Q_T,
			\\
			&\operatorname{div}\widetilde{\mathbf{v}}_n ({x},\widetilde{t})=0, \ \ \ \text{ in }Q_{T},
			\\
			&\widetilde{\mathbf{v}}_n({x},0)=
			\overline{\mathbf{v}}_n({x},1), \ \ \ \widetilde{\vartheta}_n({x},0)=
			\overline{\vartheta}_n({x},1), \ \ \ \text{ in }\Omega,
			\\
			&\widetilde{\mathbf{v}}_n({x},\widetilde{t})=\mathbf{0}, \ \ \widetilde{\vartheta}_n({x},\widetilde{t})=0, \ \ \ \text{ on }\Gamma_{T},
		\end{align}
		where
		\begin{align}
			\label{sphi_n}
			\begin{split}
				\widetilde{\varphi}_n(\widetilde{t};[\widetilde{\mathbf{v}}_n,\widetilde{\vartheta}_n])&=\widetilde{f}_n'(\widetilde{t})
			-\int_\Omega(\widetilde{\mathbf{v}}_n({x},\widetilde{t})
				\cdot\nabla)\boldsymbol{\omega}({x})\cdot\widetilde{\mathbf{v}}_n({x},\widetilde{t})\,d{x}
				\\
				&\quad	+\nu\int_\Omega\nabla\widetilde{\mathbf{v}}_n({x},\widetilde{t}):\nabla\boldsymbol{\omega}({x})\,d{x}-\int_\Omega\widetilde{\mathbf{g}}_n(x,\widetilde{t})
				\cdot\boldsymbol{\omega}({x})\widetilde{\vartheta}_{n}(x,\widetilde{t})\,dx,
			\end{split}
			\\
			\begin{split}
				\widetilde{\psi}_n(\widetilde{t};[\widetilde{\mathbf{v}}_n,\widetilde{\vartheta}_n])&=
				\widetilde{g}_n'(\widetilde{t})
				+\mu\int_\Omega\nabla\widetilde{\vartheta}_n({x},\widetilde{t})
				\cdot\nabla\eta({x})\,d{x}-\int_\Omega
				(\widetilde{\mathbf{v}}_n({x},\widetilde{t})
				\cdot\nabla)\eta({x})\widetilde{\vartheta}_n({x},\widetilde{t})\,d{x}.
			\end{split}
		\end{align}
	\end{subequations}
	
	{ We are now in a position to state the main result of the paper.}
	\begin{theorem}\label{theorem3.1} Let Conditions \ref{C} be fulfilled.
		Then the following assertions hold true.
		\begin{itemize}
			\item[{\bf 1.}] The family of scaled solutions  $\displaystyle\{(\overline{\mathbf{v}}_n,\overline{\varphi}_n,\overline{\vartheta}_n,\overline{\psi}_n)\}_{n \geq n_0}$ of problem \eqref{eq3.1} is relatively compact in the Cartesian product of the spaces  $\mathrm{C}([0,1];\mathbf{H})\times{L^{2}(0,1)}\times \mathrm{C}([0,1];L^{2}(\Omega))\times L^{2}(0,1)$. Namely, there exist subsequences of $\displaystyle\{\overline{\mathbf{v}}_n\}_{n\geq n_0}$, $\displaystyle\{\overline{\vartheta}_n\}_{n\geq n_0}$, $\displaystyle\{\overline{\varphi}_n\}_{n\geq n_0}$, and $\displaystyle\{\overline{\psi}_n\}_{n\geq n_0}$, still labeled by $n$, and a limit vector function $\displaystyle\overline{\mathbf{v}}\in \mathrm{C}([0,1];\mathbf{H})$ and a limit scalar-function $\displaystyle\overline{\vartheta}\in \mathrm{C}([0,1];L^2(\Omega))$ such that
			\begin{subequations}
				\begin{align}
					\label{eq3.3a}
					\overline{\varphi}_n \xrightarrow[n\to\infty]{}\overline{f}^\prime
					& \ \ \mbox{ strongly in }L^{2}(0,1),
					\\
					\label{eq3.3b} 
					\overline{\psi}_n \xrightarrow[n\to\infty]{} \overline{g}^\prime
					& \ \ \mbox{ strongly in }L^{2}(0,1),
					\\
					\label{eq3.3c} 
					\overline{\mathbf{v}}_n
					\xrightarrow[n\to\infty]{} \overline{\mathbf{v}}
					& \ \ \mbox{ strongly in } \mathrm{C}([0,1];\mathbf{H}),
					\\
					\label{eq3.3d} 
					\overline{\vartheta}_n
					\xrightarrow[n\to\infty]{} \overline{\vartheta}& \  \ \mbox{ strongly in } \mathrm{C}([0,1];L^2(\Omega)),
				\end{align}
			\end{subequations}
			where the pair of functions $(\overline{\mathbf{v}},\overline{\vartheta}) $ resolve 
			the direct problem
			\begin{subequations}
				\label{eq3.4}
				\begin{align}
					\begin{split} & \partial_{\overline{t}} \overline{\mathbf{v}}
						+
						\nabla\overline{\pi}=\overline{f}^\prime({\overline{t}})\,\boldsymbol{\omega}({x}), \ \
						\partial_{\overline{t}} \overline{\vartheta}=\overline{g}^\prime({\overline{t}})\,\eta({x}),
						\ \ \text{ in }Q_1,
					\end{split}
					\\
					&\operatorname{div} \overline{\mathbf{v}}({x},{\overline{t}})=0, \ \ \text{ in }Q_1,
					\\
					&\overline{\mathbf{v}}({x},0)=\mathbf{v}_0({x}), \ \ \overline{\vartheta}({x},0)=\vartheta_0({x}), \ \ \text{ in }\Omega,
					\\
					&\overline{\mathbf{v}}({x},{\overline{t}})=\mathbf{0}, \ \ 
					\overline{\vartheta}({x},{\overline{t}})=0
					\ \ \text{ on
					}\Gamma_1.
				\end{align}
			\end{subequations}
			In the case the pressure does not depend on $x$, i.e., $\nabla\overline{\pi}=0$, problem \eqref{eq3.4} is resolved in an explicit way
			\[
			\overline{\mathbf{v}}({x},{\overline{t}})=\mathbf{v}_0({x})+(\overline{f}({\overline{t}})-f(0))\boldsymbol{\omega}({x}),
			\ \
			\overline{\vartheta}({x},{\overline{t}})=\vartheta_0({x})+(\overline{g}({\overline{t}})-g(0))\eta({x}).
			\]
			We get the final trace
			\begin{equation}
				\label{eq3.5}
				\overline{\mathbf{v}}({x},1)=\mathbf{v}_0({x})+(\overline{f}(1)-f(0))\boldsymbol{\omega}({x}), \ \ \overline{\vartheta}({x},1)=\vartheta_0({x})+(\overline{g}(1)-g(0))\eta({x}).
			\end{equation}
			
			\item[{\bf 2.}] The family of shifted solutions  {$\displaystyle\{(\widetilde{\mathbf{v}}_n,\widetilde{\vartheta}_n,\widetilde{\varphi}_n,\widetilde{\psi}_n)\}_{n\geq n_0}$} of problem \eqref{eq3.2} is strongly and weakly relatively compact in the Cartesian product of four spaces 
			\[
			\mathrm{C}([0,T];\mathbf{H})\times \mathrm{C}([0,T];L^2(\Omega))\times L^{2}(0,T)\times L^{2}(0,T).
			\] 
			Therefore, there exist subsequences of
			$\displaystyle\{\widetilde{\mathbf{v}}_n\}_{n\geq n_0}$,
			$\displaystyle\{\widetilde{\vartheta}_n\}_{n\geq n_0}$,
			$\{\widetilde{\varphi}_n\}_{n\geq n_0}$, and $\{\widetilde{\psi}_n\}_{n\geq n_0}$, still labeled by $n$, and limit
			functions 
			
			\[
			\mathbf{v}\in {\mathrm{C}([0,T];\mathbf{H})}\cap L^2(0,T;\mathbf{V}^2),\quad
			\vartheta\in {\mathrm{C}([0,T];L^2(\Omega))}\cap L^2(0,T;W_0^{2,2}(\Omega)),
			\]
			and  
			$\varphi,\psi\in L^{2}(0,T)$, such that 
			\begin{subequations}
				\label{eq25}
				\begin{align}
					\widetilde{\mathbf{v}}_n \xrightarrow[n\to\infty]{} \mathbf{v}\quad
					&\text{strongly in }{\mathrm{C}([0,T];\mathbf{H})},
					\\
					\widetilde{\mathbf{v}}_n \xrightharpoonup[n\to\infty]{*} \mathbf{v}\quad
					& \text{weakly$^*$ in } L^\infty(0,T;\mathbf{V}),
					\\
					\widetilde{\mathbf{v}}_n \xrightharpoonup[n\to\infty]{} \mathbf{v}\quad
					&\text{weakly in } L^2(0,T;\mathbf{V}^2),
					\\
					\widetilde{\vartheta}_n \xrightarrow[n\to\infty]{} \vartheta\quad
					&\text{strongly in }{\mathrm{C}([0,T];L^2(\Omega))},
					\\
					\widetilde{\vartheta}_n \xrightharpoonup[n\to\infty]{*} \vartheta\quad
					& \text{weakly$^*$ in } L^\infty(0,T;W_0^{1,2}(\Omega)),
					\\
					\widetilde{\vartheta}_n\xrightharpoonup[n\to\infty]{} \vartheta\quad
					&\text{weakly in } L^2(0,T;W_0^{2,2}(\Omega)),
					\\
					\widetilde{\varphi}_n \xrightharpoonup[n\to\infty]{} \varphi\quad
					&\text{weakly in } L^{2}(0,T),
					\\
					\widetilde{\psi}_n \xrightharpoonup[n\to\infty]{} \psi\quad
					&\text{weakly in } L^{2}(0,T),
				\end{align}
			\end{subequations}
			where the quadruple $(\mathbf{v},\vartheta,\varphi,\psi)$ is a \textbf{weak regular solution} to the following inverse problem:
			\begin{subequations}
				\label{eq3.6}
				\begin{align}
					\begin{split} &\partial_{t} \mathbf{v}+
						(\mathbf{v}\cdot\nabla)\mathbf{v}+\nabla\pi=
						\nu\Delta\mathbf{v}+\mathbf{g}\vartheta+\varphi(t)\,\boldsymbol{\omega}({x}), \ \ \ \text{
							in }Q_T,
					\end{split}
					\\
					\begin{split} &\partial_{t}\vartheta+
						(\mathbf{v}\cdot\nabla)\vartheta=
						\mu\Delta\vartheta+\psi(t)\,\eta({x}), \ \ \ \text{ in }Q_T,
					\end{split}
					\\
					&\operatorname{div}\mathbf{v}({x},t)=0, \ \ \ \text{ in }Q_T,
					\\
					&\int_\Omega \mathbf{v}({x},t)\cdot\boldsymbol{\omega}({x})\,d{x}=f(t), \ 
					\ \ \text{ in }(0,T),
					\\
					&
					\int_\Omega \vartheta({x},t)\eta({x})\,d{x}=g(t), \ 
					\ \ \text{ in }(0,T),
					\\
					&\mathbf{v}({x},0)=\overline{\mathbf{v}}({x},1), \ \ \ \
					\vartheta({x},0)=\overline{\vartheta}({x},1), \ \ \ \text{ in }\Omega,
					\\
					&\mathbf{v}({x},t)=\mathbf{0}, \ \ \ \ \vartheta({x},t)=0, \ \ \ \text{ on }\Gamma_T,
				\end{align}
			\end{subequations}
			with a new initial velocity $\mathbf{v}({x},0)$ and a new initial temperature $\vartheta({x},0)$, which correspond to the final data of problem \eqref{eq3.4}. In general case, see \eqref{eq3.5}; it follows that $\mathbf{v}({x},0)\neq\mathbf{v}_0(x)$ and $\vartheta({x},0)\neq\vartheta_0(x)$.
			\item[{\bf 3.}] The families of functions $\displaystyle\{\mathbf{v}_n\}_{n\geq n_0}$  and $\displaystyle\{\vartheta_n\}_{n\geq n_0}$ are relatively compact in $L^2(0,T;\mathbf{H})$ and $L^2(0,T;L^2(\Omega))$, respectively, and $\{\mathbf{1}_{(1/n,T)}\varphi_n\}_{n\geq n_0}$
			and $\{\mathbf{1}_{(1/n,T)}\psi_n\}_{n\geq n_0}$ are weakly relatively compact in $L^{2}(0,T)$.
			Therefore, there exist subsequences of
			$\displaystyle\{\mathbf{v}_n\}_{n\geq n_0}$, $\displaystyle\{\vartheta_n\}_{n\geq n_0}$, $\{\mathbf{1}_{(1/n,T)}\varphi_n\}_{n\geq n_0}$, and
			$\{\mathbf{1}_{(1/n,T)}\psi_n\}_{n\geq n_0}$, still labeled by $n$, such that
			\begin{subequations}
				\label{eq27}
				\begin{align}
					\mathbf{v}_n \xrightarrow[n\to\infty]{} \mathbf{v}\quad
					&\text{strongly in } L^2(0,T;\mathbf{H}),
					\\
					\mathbf{v}_n \xrightharpoonup[n\to\infty]{*} \mathbf{v}\quad
					& \text{weakly$^*$ in } L^\infty(0,T;\mathbf{V}),
					\\
					\mathbf{v}_n \xrightharpoonup[n\to\infty]{}
					\mathbf{v}
					\quad
					&\text{weakly in } L^2(0,T;\mathbf{V}^2),
					\\
					\vartheta_n \xrightarrow[n\to\infty]{} \vartheta\quad
					&\text{strongly in } L^2(0,T;L^2(\Omega)),
					\\
					\vartheta_n \xrightharpoonup[n\to\infty]{*} \vartheta\quad
					& \text{weakly$^*$ in } L^\infty(0,T;W_0^{1,2}(\Omega)),
					\\
					\vartheta_n \xrightharpoonup[n\to\infty]{} \vartheta
					\quad
					&\text{weakly in } L^2(0,T;W_0^{2,2}(\Omega)),
					\\
					{\mathbf{1}_{(\frac1n,T)}\varphi_n \xrightharpoonup[n\to\infty]{}  \varphi}\quad
					&\text{weakly in } L^{2}(0,T),
					\\
					\mathbf{1}_{(\frac1n,T)}\psi_n \xrightharpoonup[n\to\infty]{} \psi\quad
					&\text{weakly in } L^{2}(0,T),
				\end{align}
			\end{subequations}
			
			where the quadruple $(\mathbf{v},\vartheta,\varphi,\psi)$ is a \textbf{weak regular solution} to problem
			
			\eqref{eq3.6}.
		\end{itemize}
	\end{theorem}
	
	\begin{remark}
		The problem \eqref{eq3.1} is inverse, but the limit problem \eqref{eq3.4} is direct. This is explained by the fact that the terms \eqref{eq3.1df} and \eqref{eq3.1dg} in the limit as $n\to\infty$ do not depend on an unknown solution.
	\end{remark}
	
	\begin{remark}
		
		Formally, the original solution is a \textbf{weak regular solution}, and a scaled solution is a \textbf{weak regular solution}, but the shifted solution is a \textbf{weak solution}, since the shifted initial data lacks regularity as the final trace of a scaled solution.
	\end{remark}
	
	The proof of this theorem is split into subsections \ref{sec3.1}--\ref{sec3.3}.
	
	\subsubsection*{\bf Scaling and shifting}
	
	Following the carefully developed procedure outlined previously (see, for example, \cite{Ant-2026,AKS,AKS-2024,AKSS-2024,AKSS-2024-2,KS-NARWA-2022,KS-JEPE-2022,KS-NA-TMA-2023}), let us consider some preliminary aspects before we address the limiting passage as $n \to \infty$ in \eqref{eq2.1}.
	
	Assuming that the test functions $\boldsymbol{\phi}$ and $\zeta$ in the equations \eqref{eq2.1a} and \eqref{eq2.1b} vanish in a neighborhood of the subsection $\{t=T\}$, let us separate the integrals in \eqref{eq2.1} over segments $(0,1/n)$ and $(1/n,T_n)$ from each other:
	\begin{subequations}
		\label{eq5.1}
		\begin{align}
			& \int_0^{1/n} \int_\Omega \Big(-\mathbf{v}_n \cdot \partial_t \boldsymbol{\phi} +(\mathbf{v}_n \cdot \nabla)\mathbf{v}_n\cdot\boldsymbol{\phi} + \nu \,
			\nabla\mathbf{v}_n:
			\nabla \boldsymbol{\phi}
			-\vartheta_n\mathbf{g}\cdot\boldsymbol{\phi}
			-\varphi_n(t)\boldsymbol{\omega}({x})\cdot\boldsymbol{\phi}\Big)\, d{x} dt\nonumber
			\\
			&
			-\int_\Omega \mathbf{v}_0({x}) \cdot \boldsymbol{\phi}({x},0)\, d{x} 
			+ \int_{1/n}^{T_n} \int_\Omega \Big(-\mathbf{v}_n \cdot \partial_t \boldsymbol{\phi} 
			+(\mathbf{v}_n \cdot \nabla)\mathbf{v}_n\cdot\boldsymbol{\phi}
			+\nu\nabla\mathbf{v}_n: \nabla \boldsymbol{\phi} 
			\nonumber \\
			& 
			-\vartheta_n\mathbf{g}\cdot\boldsymbol{\phi}
			-\varphi_n(t)\boldsymbol{\omega}({x})\cdot\boldsymbol{\phi}\Big)\, d{x} dt =0, \label{eq5.1a}
			\\
			& \int_0^{1/n} \int_\Omega \Big(-\vartheta_n\partial_t\zeta 
			+(\mathbf{v}_n \cdot \nabla)\vartheta_n\zeta
			+ \mu \nabla\vartheta_n\cdot\nabla \zeta -\psi_n(t)\eta({x})\zeta \Big)\, d{x} dt \nonumber
			\\& 
			- \int_\Omega \vartheta_0({x})\zeta({x},0)\, d{x}  + \int_{1/n}^{T_n} \int_\Omega 
			\Big(-\vartheta_n\partial_t\zeta 
			+(\mathbf{v}_n \cdot \nabla)\vartheta_n\zeta
			+ \mu \nabla\vartheta_n\cdot\nabla \zeta 
			\nonumber 
			\\
			&
			-\psi_n(t)\eta({x})\zeta \Big)\, d{x} dt =0. \label{eq5.1b}
		\end{align}
	\end{subequations}
	Thus, \eqref{eq5.1} takes the form
	\begin{subequations}
		\label{eq5.4}
		\begin{align}
			& \int_{Q_1}
			\Big(-\overline{\mathbf{v}}_n({x},{\overline{t}}) \cdot \partial_{\overline{t}} \boldsymbol{\phi}({x},n^{-1}{\overline{t}})+
			\frac{1}{n} (\overline{\mathbf{v}}_n({x},{\overline{t}}) \cdot\nabla )\overline{\mathbf{v}}_n({x},{\overline{t}})
			\cdot
			\boldsymbol{\phi}({x},n^{-1}{\overline{t}}) 
			\nonumber 
		\\[0.6ex]
		&\quad	+\frac{\nu}{n}\nabla\overline{\mathbf{v}}_n({x},{\overline{t}}):\nabla\boldsymbol{\phi}({x},n^{-1}{\overline{t}})
			-\frac1n\overline{\vartheta}_n({x},{\overline{t}})
			\overline{\mathbf{g}}_n({x},{\overline{t}})
			\cdot\boldsymbol{\phi}({x},n^{-1}{\overline{t}})
			\nonumber 
		\\[0.6ex]
		&\quad	-\overline{\varphi}_n({\overline{t}})\boldsymbol{\omega}({x})
		\cdot\boldsymbol{\phi}({x},n^{-1}{\overline{t}})	\Big)\, d{x} d{\overline{t}} -\int_\Omega \mathbf{v}_0({x}) \cdot \boldsymbol{\phi}({x},0)\, d{x}  
			\nonumber 
			\\[0.6ex]
			&
			+\int_{Q_T}
			\Big(-\widetilde{\mathbf{v}}_n ({x},\widetilde{t}) \cdot \partial_{\widetilde{t}}\, \boldsymbol{\phi}({x},\widetilde{t}+1/n)+ (\widetilde{\mathbf{v}}_n({x},\widetilde{t}) \cdot\nabla) \widetilde{\mathbf{v}}_n({x},\widetilde{t})\cdot
			\boldsymbol{\phi}({x},\widetilde{t}+1/n) \nonumber 
			\\[0.6ex]
			&\quad +\nu \nabla\widetilde{\mathbf{v}}_n({x},\widetilde{t}) : \nabla \boldsymbol{\phi}({x},\widetilde{t} +1/n)-\widetilde{\varphi}_n(\widetilde{t})\boldsymbol{\omega}({x})\cdot \boldsymbol{\phi}({x},\widetilde{t} +1/n)\nonumber
			\\[0.6ex]
			& \quad-\widetilde{\vartheta}_n({x},\widetilde{t})
			\widetilde{\mathbf{g}}_n({x},\widetilde{t})
			\cdot\boldsymbol{\phi}({x},\widetilde{t}+1/n)
			\Big)\, d{x} d\widetilde{t}=0, 
			\label{eq5.4a}
			\\[0.6ex]
			&\int_{Q_1} \Big(-\overline{\vartheta}_n({x},{\overline{t}})\partial_{\overline{t}}\zeta(x,n^{-1}{\overline{t}})
			+\frac1n(
			\overline{\mathbf{v}}_n({x},{\overline{t}}) \cdot \nabla)\overline{\vartheta}_n({x},{\overline{t}})\zeta(x,n^{-1}{\overline{t}})
			\nonumber
		\\[0.6ex]
		&\quad	+ 
			\frac{\mu}n\nabla\overline{\vartheta}_n({x},{\overline{t}})\cdot\nabla \zeta(x,n^{-1}{\overline{t}}) 
		-\overline{\psi}_n(\overline{t})\eta({x})\zeta(x,n^{-1}{\overline{t}}) \Big)\, d{x} d\overline{t}  - \int_\Omega \vartheta_0({x})\zeta({x},0)\, d{x} 
			\nonumber
			\\[0.6ex]
			& \quad+ \int_{Q_T} 
			\Big(-\widetilde{\vartheta}_n({x},\widetilde{t})\partial_{\widetilde{t}}\zeta(x,\widetilde{t}+1/n) 
			+(\widetilde{\mathbf{v}}_n({x},\widetilde{t}) \cdot \nabla)\widetilde{\vartheta}_n({x},\widetilde{t})\zeta(x,\widetilde{t}+1/n)
			\nonumber
			\\[0.6ex]
			&\quad+ \mu \nabla\widetilde{\vartheta}_n({x},\widetilde{t})\cdot\nabla \zeta(x,\widetilde{t}+1/n) -\widetilde{\psi}_n(\widetilde{t})\eta({x})\zeta(x,\widetilde{t}+1/n)
			\Big)\, d{x} d\widetilde{t}=0.
			\label{eq5.4b}
		\end{align}
	\end{subequations}
	
	In accordance  with the further limiting passage as $n \to \infty$ in \eqref{eq5.4}, we take the test vector function $\boldsymbol{\phi}=\boldsymbol{\phi}_n({x},t)$ and the test scalar function $\zeta=\zeta_n({x},t)$, depending on $n\in \mathbb{N}$ in the following form:
	\begin{subequations}
		\label{eq5.5}
		\begin{align}
			\label{eq5.5a}
			&\boldsymbol{\phi}_n({x},t)=\left\{
			\begin{array}{ll}
				\displaystyle \overline{\boldsymbol{\phi}}({x},{\overline{t}}) \equiv \overline{\boldsymbol{\phi}}({x},nt) & \text{ for } t\in [0,1/n], \text{ i.e., for } {\overline{t}} \in [0,1], \\[2ex]
				\widetilde{\boldsymbol{\phi}}({x},\widetilde{t}) \equiv \widetilde{\boldsymbol{\phi}}({x},t-1/n) & \text{ for } t\in [1/n,T_n],   \text{ i.e., for } \widetilde{t} \in [0,T],
			\end{array}
			\right.
			\\
			&
			\zeta_n({x},t)=\left\{
			\begin{array}{ll}
				\displaystyle \overline{\zeta}({x},{\overline{t}}) \equiv \overline{\zeta}({x},nt) & \text{ for } t\in [0,1/n], \text{ i.e., for } {\overline{t}} \in [0,1], \\[2ex]
				\widetilde{\zeta}({x},\widetilde{t}) \equiv \widetilde{\zeta}({x},t-1/n) & \text{ for } t\in [1/n,T_n],   \text{ i.e., for } \widetilde{t} \in [0,T],
			\end{array}
			\right.
		\end{align}
	\end{subequations}
	where $\overline{\boldsymbol{\phi}} = \overline{\boldsymbol{\phi}}({x},{\overline{t}})$ and $\widetilde{\boldsymbol{\phi}}=\widetilde{\boldsymbol{\phi}}({x},\widetilde{t})$ are arbitrary smooth vector-valued test  functions defined on $\overline{\Omega} \times [0,1]$ and $\overline{\Omega} \times [0,T]$, respectively, and $\overline{\zeta} = \overline{\zeta}({x},{\overline{t}})$ and $\widetilde{\zeta}=\widetilde{\zeta}({x},\widetilde{t})$ are arbitrary smooth scalar-valued  test  functions defined on $\overline{\Omega} \times [0,1]$ and $\overline{\Omega} \times [0,T]$, respectively, such that $\overline{\boldsymbol{\phi}} = \widetilde{\boldsymbol{\phi}} \equiv 0$ and $\overline{\zeta} = \widetilde{\zeta} \equiv 0$ hold in a neighbourhood of $\partial \Omega$, $\widetilde{\boldsymbol{\phi}} \equiv 0$ and $\widetilde{\zeta} \equiv 0$ hold in a neighbourhood of the plane $\{\widetilde{t}=T\}$, and the matching conditions
	\begin{equation}
		\label{eq5.6}
		\overline{\boldsymbol{\phi}}({x},1^{-}) = \widetilde{\boldsymbol{\phi}}({x},0^+)
		\quad \text{and} \quad
		\overline{\zeta}({x},1^{-}) = \widetilde{\zeta}({x},0^+)
	\end{equation}
	hold. {We notice that condition \eqref{eq5.6} yields the weak derivatives 
		\begin{align*}
			\partial_t \boldsymbol{\phi}_n \in L^2(0,T_n;\mathbf{H}) \ \ \ \text{and} \ \ \  \partial_t \zeta_n \in L^2(0,T_n;L^2(\Omega)),
		\end{align*}
		which imply that 
		\begin{align*}
			\boldsymbol{\phi}_n \! \in \mathrm{C}([0,T_n];\mathbf{V}) \ \ \ \text{and} \ \ \  \zeta_n \ \in \mathrm{C}([0,T_n];W_0^{1,2}(\Omega)),
		\end{align*}
		respectively, see \cite[Chapter 5, Section 5.9]{Evans}. Since 
		\begin{align*}
			\partial_{\overline{t}} \overline{\boldsymbol{\phi}}_n \in L^2(0,1;\mathbf{H}),  \  \partial_{\overline{t}} \overline{\zeta}_n \in L^2(0,1;L^2(\Omega)), \ 
			\partial_{\widetilde{t}} \widetilde{\boldsymbol{\phi}}_n \in L^2(0,T;\mathbf{H}),   \ \partial_{\widetilde{t}} \widetilde{\zeta}_n \in L^2(0,T;L^2(\Omega)),
		\end{align*}
		it follows, by combining the two time intervals, that 
		\begin{align*}
			\partial_t \boldsymbol{\phi}_n \in L^2(0,T_n;\mathbf{H}) \ \ \ \text{and} \ \ \  \partial_t \zeta_n \in L^2(0,T_n;L^2(\Omega)).
		\end{align*}
		Similarly, we  have
		\begin{align*}
			\overline{\boldsymbol{\phi}}_n \in L^2(0,1;\mathbf{V}^2),  \   \overline{\zeta}_n \in L^2(0,1;W_0^{2,2}(\Omega)), \ 
			\widetilde{\boldsymbol{\phi}}_n \in L^2(0,T;\mathbf{V}^2),  \   \widetilde{\zeta}_n \in L^2(0,T;W_0^{2,2}(\Omega)). 
		\end{align*}
		Hence, by combining the two time intervals,  we obtain
		\begin{align*}
			\boldsymbol{\phi}_n \in L^2(0,T_n;\mathbf{V}^2) \ \ \ \text{and} \ \ \  \zeta_n \in L^2(0,T_n;W_0^{2,2}(\Omega)).
		\end{align*}
		Moreover, in view of the matching condition \eqref{eq5.6}, the one-sided traces are understood in the $\mathbf{H}$-sense and $L^2(\Omega)$-sense, respectively, and there is no discontinuity in $\boldsymbol{\phi}_n$ and $\zeta_n$ at $t=\frac1n$, which corresponds to $\overline{t}=1$ and $\widetilde{t}=0$. Therefore, $\boldsymbol{\phi}_n$ and $\zeta_n$ are admissible test functions for \eqref{eq5.4a} and \eqref{eq5.4b}, respectively. }
	Inserting \eqref{eq5.5} into \eqref{eq5.4}, we obtain
	\begin{subequations}
		\label{eq5.7}
		\begin{align}
			& \int_{Q_1} \Big(-\overline{\mathbf{v}}_n({x},{\overline{t}}) \cdot \partial_{\overline{t}} \overline{\boldsymbol{\phi}}({x},{\overline{t}}) - \frac{1}{n} (\overline{\mathbf{v}}_n({x},{\overline{t}}) \cdot \nabla)\overline{\boldsymbol{\phi}}({x},{\overline{t}})  
			\cdot \overline{\mathbf{v}}_n({x},{\overline{t}}) 
			+ \frac{\nu}{n} \nabla\overline{\mathbf{v}}_n({x},{\overline{t}}): \nabla \overline{\boldsymbol{\phi}}({x},{\overline{t}})
			\nonumber
			\\[0.8ex]
			&-\frac1n\overline{\vartheta}_n({x},{\overline{t}})
			\overline{\mathbf{g}}_n({x},{\overline{t}})
			\cdot\overline{\boldsymbol{\phi}}({x},\overline{t})
			-
			\overline{\varphi}_n({\overline{t}})\boldsymbol{\omega}({x})
			\cdot\overline{\boldsymbol{\phi}}({x},{\overline{t}}) \Big)\,d{x}d{\overline{t}}-\int_\Omega \mathbf{v}_0({x}) \cdot \overline{\boldsymbol{\phi}}({x},0)\, d{x}  
			\nonumber 
			\\[0.8ex]
			&+\int_{Q_T} \Big(-\widetilde{\mathbf{v}}_n ({x},\widetilde{t}) \cdot \partial_{\widetilde{t}}\, \widetilde{\boldsymbol{\phi}}({x},\widetilde{t}) - (\widetilde{\mathbf{v}}_n({x},\widetilde{t}) \cdot\nabla)\widetilde{\boldsymbol{\phi}}({x},\widetilde{t})\cdot
			\widetilde{\mathbf{v}}_n({x},\widetilde{t})
			+\nu\nabla\widetilde{\mathbf{v}}_n({x},{\widetilde{t}}) : \nabla \widetilde{\boldsymbol{\phi}}({x},\widetilde{t})
			\nonumber 
			\\[0.8ex]
			& -\widetilde{\vartheta}_n({x},\widetilde{t})
			\widetilde{\mathbf{g}}_n({x},\widetilde{t})
			\cdot\widetilde{\boldsymbol{\phi}}({x},\widetilde{t})
			-\widetilde{\varphi}_n(\widetilde{t})\boldsymbol{\omega} ({x})
			\cdot\widetilde{\boldsymbol{\phi}}({x},\widetilde{t})\Big)\, d{x} d\widetilde{t} =0,   \label{eq5.7a}
			\\
			& \int_{Q_1} \Big(-\overline{\vartheta}_n({x},{\overline{t}}) 
			\partial_{\overline{t}} \overline{\zeta}({x},{\overline{t}}) 
			+ 
			\frac{1}{n} (\overline{\mathbf{v}}_n({x},{\overline{t}}) \cdot
			\nabla)\overline{\vartheta}_n({x},{\overline{t}})
			\overline{\zeta}({x},{\overline{t}})  + \frac{\mu}{n} \nabla\overline{\vartheta}_n({x},{\overline{t}})\cdot 
			\nabla \overline{\zeta}({x},{\overline{t}})
			\nonumber
			\\[0.8ex]
			&
			-\overline{\psi}_n({\overline{t}})\eta({x})
			\overline{\zeta}({x},{\overline{t}}) \Big)\,d{x}d{\overline{t}}
			-\int_\Omega \vartheta_0({x}) \overline{\zeta}({x},0)\, d{x}
			+\int_{Q_T} \Big(-\widetilde{\vartheta}_n ({x},\widetilde{t})
			\partial_{\widetilde{t}}\, \widetilde{\zeta}({x},\widetilde{t}) 
			\nonumber 
			\\[0.8ex]
			& 
			-
			(\widetilde{\mathbf{v}}_n({x},\widetilde{t})\cdot\nabla)
			\widetilde{\zeta}({x},\widetilde{t})
			\widetilde{\vartheta}_n({x},\widetilde{t})+\mu\nabla\widetilde{\vartheta}_n({x},{\widetilde{t}})\cdot 
			\nabla \widetilde{\zeta}({x},\widetilde{t})-\widetilde{\psi}_n(\widetilde{t})\eta({x})\widetilde{\zeta}({x},{\widetilde{t}})
			\Big)\, d{x} d\widetilde{t} =0.   \label{eq5.7b}
		\end{align}
	\end{subequations}
	Furthermore, we notice that
	\begin{equation} \label{eq5.8}
		\overline{\mathbf{v}}_n({x},1^-) = \widetilde{\mathbf{v}}_n({x},0^+),
		\quad 
		\overline{\vartheta}_n({x},1^-) = \widetilde{\vartheta}_n({x},0^+),
		\quad 
		\mbox{ in } \Omega,
	\end{equation}
	in the sense of distributions, due to the regularity properties of $\overline{\mathbf{v}}_n\in \mathrm{C}([0,1],\mathbf{H})$, $\widetilde{\mathbf{v}}_n\in \mathrm{C}([0,T];\mathbf{V}^\prime)$, $\overline{\vartheta}_n\in \mathrm{C}([0,1];L^2(\Omega))$, and $\widetilde{\vartheta}_n\in \mathrm{C}([0,T];W^{-1,2}(\Omega))$. 
	
	The remainder of the proof of Theorem \ref{theorem3.1} is based also on a thorough study of \eqref{eq5.7}, taking into account \eqref{eq5.8}.
	
	\subsection{\bf Proof of item 1 in Theorem \ref{theorem3.1}}
	\label{sec3.1} 
	%We assume that $\{\overline{\mathbf{v}}_n\}_{n\geq n_0}$ and
	%$\{\overline{\varphi}_n({\overline{t}})\}_{n\geq n_0}$ are uniformly bounded in. 
	
	{ We need to introduce the notion of a weak regular solution on an infinitesimal initial layer. We recall that a weak regular solution is a weak solution satisfying additional regularity, see Definitions \ref{def:ws} and \ref{def:wrs}.}
	\begin{definition}
		\label{def3.1} A quadruple of functions $(\overline{\mathbf{v}}_n,\overline{\vartheta}_n,\overline{\varphi}_n,\overline{\psi}_n)$ is said to be a \textbf{weak regular solution} of the nonlinear inverse problem \eqref{eq3.1} with initial condition $\mathbf{v}_0 \in \mathbf{V}$ and $\vartheta_0 \in W_{0}^{1,2}(\Omega)$ if
		\begin{align*}
			&\overline{\mathbf{v}}_n \in L^\infty(0,1;\mathbf{V})\cap L^2(0,1;\mathbf{V}^2),
			\quad
			\overline{\vartheta}_n \in L^\infty(0,1;W_0^{1,2}(\Omega))\cap L^2(0,1;
			W_0^{2,2}(\Omega)),
			\\
			& \partial_{\overline{t}}\overline{\mathbf{v}}_n\in L^{2}(0,1;\mathbf{V}^\prime), \ \ \ \partial_{\overline{t}}\overline{\vartheta}_n\in L^{2}(0,1;
			W^{-1,2}(\Omega)), \ \ \
			\overline{\varphi}_n\in L^{2}(0,1), \ \ \ 
			\overline{\psi}_n\in L^{2}(0,1),
		\end{align*}
		and $(\overline{\mathbf{v}}_n,\overline{\vartheta}_n,\overline{\varphi}_n,\overline{\psi}_n)$ satisfies 
		\begin{align*}
			& \frac{d}{d{\overline{t}}}\int_\Omega
			\overline{\mathbf{v}}_n({x},{\overline{t}})\cdot\boldsymbol{\phi}({x})\,d{x}
			+\frac{\nu}{n}\int_\Omega\nabla\overline{\mathbf{v}}_n({x},{\overline{t}}):\nabla\boldsymbol{\phi}({x})\,d{x}-\frac{1}{n}\int_\Omega\overline{\vartheta}_n(x,\overline{t})
			\overline{\mathbf{g}}_n({x},{\overline{t}})\cdot\boldsymbol{\phi}({x})\,d{x}
			\\
			& -\frac1n\int_\Omega
			(\overline{\mathbf{v}}_n({x},{\overline{t}})\cdot\nabla)\boldsymbol{\phi}({x})\cdot\overline{\mathbf{v}}_n({x},{\overline{t}})\,d{x} =\overline{\varphi}_n({\overline{t}})\int_\Omega
			\boldsymbol{\omega}({x})\cdot \boldsymbol{\phi}({x})\,d{x} \quad \text{for a.e. }{\overline{t}}\in(0,1),
			\\
			& \frac{d}{d{\overline{t}}}\int_\Omega
			\overline{\vartheta}_n({x},{\overline{t}})\zeta({x})\,d{x}
			+\frac{\mu}{n}\int_\Omega\nabla\overline{\vartheta}_n({x},{\overline{t}})\cdot\nabla\zeta({x})\,d{x}
			-\frac1n\int_\Omega
			(\overline{\mathbf{v}}_n({x},{\overline{t}})\cdot\nabla)\zeta({x})
			\vartheta_n({x},{\overline{t}})\,d{x} 
			\\
			&
			=\overline{\psi}_n({\overline{t}})\int_\Omega
			\eta({x})\zeta({x})\,d{x} \quad \text{for a.e. }{\overline{t}}\in(0,1),
		\end{align*}
		for all test functions $\boldsymbol{\phi}\in\mathbf{V}$, $\zeta\in W_0^{1,2}(\Omega)$, and the corresponding compatibility conditions
		\begin{align*}
			&\int_\Omega\mathbf{v}_0({x})\cdot\boldsymbol{\omega}({x})\,d{x}=\overline{f}_n(0),
			\quad
			\int_\Omega\vartheta_0({x})\eta({x})\,d{x}=\overline{g}_n(0).
		\end{align*}
	\end{definition}
	
	\begin{proposition}[Estimate for a scaled solution]
		\label{proposition.2.1} Under Conditions \ref{C}, the problem \eqref{eq3.1} has a \textbf{weak regular solution} in the sense of Definition \ref{def3.1}, which satisfies the estimates
		\begin{align}
			&\underset{{\overline{t}}\in(0,1)} {\rm esssup} \Big (
			\|\overline{\mathbf{v}}_n(\cdot,{\overline{t}})\|_{\mathbf{H}}^2 +\|\nabla \overline{\mathbf{v}}_n(\cdot,{\overline{t}})\|_{\mathbf{H}}^2 \Big)
			+\frac1n\|\overline{\mathbf{v}}_n\|_{L^2(0,1;\mathbf{V})}^2+\frac1n\|\overline{\mathbf{v}}_n\|_{L^2(0,1;\mathbf{V}^2)}^2
			\nonumber
			\\
			&\quad +\|\partial_{\overline{t}}\overline{\mathbf{v}}_n\|_{L^2(0,1;\mathbf{V}^\prime)} 
			+\|\overline{\varphi}_n\|_{L^{2}(0,1)}\leq C,
			\label{eq:3est-global-1}
			\\
			&\underset{{\overline{t}}\in(0,1)} {\rm esssup} \Big (
			\|\overline{\vartheta}_n(\cdot,{\overline{t}})\|_{L^{2}(\Omega)}^2 +\|\nabla \overline{\vartheta}_n(\cdot,{\overline{t}})\|_{\mathbf{L}^{2}(\Omega)}^2 \Big)
			+\frac1n\|\overline{\vartheta}_n\|_{L^2(0,1;W_0^{1,2}(\Omega))}^2\nonumber
			\\
			&\quad+\frac1n\|\overline{\vartheta}_n\|_{L^2(0,1;W_0^{2,2}(\Omega))}^2
			+
			\|\partial_{\overline{t}}\overline{\vartheta}_n\|_{L^2(0,1;W^{-1,2}(\Omega))} 
			+\|\overline{\psi}_n\|_{L^{2}(0,1)}\leq C,
			\label{eq:3est-global-1+}
		\end{align}
		where a constant $C$ does not depend on $n$. 
	\end{proposition}
	The proof of Proposition \ref{proposition.2.1} relies on Theorem \ref{thm1}.\\
	\noindent{\bf Proof}.  Applying the rescaling and shifting  transformations  in the estimates obtained in Theorem \ref{thm1} and discarding the nonnegative expressions containing $\widetilde{\mathbf{v}}_n$ and $\widetilde{\vartheta}_n$, we derive the following estimates for the family $\displaystyle \bigl\{ \overline{\mathbf{v}}_n\bigr\}_{n\geq n_0}$ and $\displaystyle \bigl\{ \overline{\vartheta}_n\bigr\}_{n\geq n_0}$:
	\begin{align}
		\label{eq3.15}
		&\underset{0 \leq {\overline{t}} \leq 1}{\rm esssup}  \Big( \|\overline{\mathbf{v}}_n(\cdot,{\overline{t}})\|_{\mathbf{H}}^2+\| \nabla \overline{\mathbf{v}}_n(\cdot,{\overline{t}})\|_{\mathbf{H}}^2 \Big)+\frac{1}{n}
		\|\overline{\mathbf{v}}_n\|_{L^2(0,1;\mathbf{V})}^2
		+\frac{1}{n}
		\|\overline{\mathbf{v}}_n\|_{L^2(0,1;\mathbf{V}^2)}^2
		\nonumber
		\\
		&
		\leq \underset{0 \leq {t} \leq T_n }{\rm esssup}\Big(\|\mathbf{v}_n(\cdot,t)\|_{\mathbf{H}}^2+\|\nabla\mathbf{v}_n(\cdot,t)\|_{\mathbf{H}}^2\Big)+\|\mathbf{v}_n\|_{L^2(0,T_n;\mathbf{V})}^2
		+\|\mathbf{v}_n\|_{L^2(0,T_n;\mathbf{V}^2)}^2
		\leq C,
		\\
		&\underset{0 \leq {\overline{t}} \leq 1}{\rm esssup}  \Big( \|\overline{\vartheta}_n(\cdot,{\overline{t}})\|_{L^2(\Omega)}^2+\| \nabla \overline{\vartheta}_n(\cdot,{\overline{t}})\|_{\mathbf{L}^2(\Omega)}^2 \Big)+\frac{1}{n}
		\|\overline{\vartheta}_n\|_{L^2(0,1;W_0^{1,2}(\Omega))}^2
		+\frac{1}{n}
		\|\overline{\vartheta}_n\|_{L^2(0,1;W_0^{2,2}(\Omega))}^2
		\nonumber
		\\
		&
		\leq \underset{0 \leq {t} \leq T_n }{\rm esssup}\Big(\|\vartheta_n(\cdot,t)\|_{L^2(\Omega)}^2
		+\|\nabla\vartheta_n(\cdot,t)\|_{\mathbf{L}^2(\Omega)}^2\Big)
		+\|\vartheta_n\|_{L^2(0,T_n;W_0^{1,2}(\Omega))}^2
		\nonumber
		\\
		&\qquad+\|\vartheta_n\|_{L^2(0,T_n;W_0^{2,2}(\Omega))}^2
		\leq C.
	\end{align}
	From these estimates, we get the uniform estimates for $\overline{\varphi}_n,\overline{\psi}_n\in L^{2}(0,1)$ (see Lemma \ref{1Globest:1} also):
	\begin{align*}
		&\|\overline{\varphi}_n\|_{L^{2}(0,1)}=
		\Big\|\overline{f}_n^{\prime}
		(\cdot)+\frac{\nu}{n}\int_\Omega\nabla\overline{\mathbf{v}}_n({x},\cdot)
		:\nabla\boldsymbol{\omega}({x})\,d{x} 
		+\frac{1}{n}\int_\Omega\overline{\vartheta}_n({x},\cdot)
		\overline{\mathbf{g}}_n({x},\cdot)\cdot
		\boldsymbol{\omega}({x})\,d{x}
		\nonumber
		\\
		&\quad
		-\frac1n\int_\Omega
		(\overline{\mathbf{v}}_n({x},\cdot)\cdot\nabla)\boldsymbol{\omega}({x})\cdot\overline{\mathbf{v}}_n({x},\cdot)\,d{x}
		\Big\|_{L^{2}(0,1)}
		\\
		&\leq\|\overline{f}_n'\|_{L^{2}(0,1)}
		+\frac{\nu}{n}
		\Big\|\int_\Omega \nabla\overline{\mathbf{v}}_n({x},\cdot)
		:\nabla\boldsymbol{\omega}({x})\,d{x}\Big\|_{L^{2}(0,1)}
		+\frac{1}{n}\Big\|\int_\Omega\overline{\vartheta}_n({x},\cdot)
		\overline{\mathbf{g}}_n({x},\cdot)\cdot
		\boldsymbol{\omega}({x})\,d{x}\Big\|_{L^{2}(0,1)}
		\\
		&\quad
		+\frac1n\Big\|\int_\Omega
		(\overline{\mathbf{v}}_n({x},\cdot)\cdot\nabla)\boldsymbol{\omega}({x})\cdot\overline{\mathbf{v}}_n({x},\cdot)\,d{x}\Big\|_{L^{2}(0,1)}\leq C,
	\end{align*}
	and
	\begin{align*}
		\|\overline{\psi}_n\|_{L^{2}(0,1)}&=
		\Big\|\overline{g}_n^{\prime}
		(\cdot)+\frac{\mu}{n}\int_\Omega \nabla\overline{\vartheta}_n({x},\cdot)
		\cdot\nabla\eta({x})\,d{x}
		-\frac1n\int_\Omega
		(\overline{\mathbf{v}}_n({x},\cdot)\cdot\nabla)\eta({x})
		\overline{\vartheta}_n({x},\cdot)\,d{x}
		\Big\|_{L^{2}(0,1)}
		\\
		&\leq
		\|\overline{g}_n'\|_{L^{2}(0,1)}
		+\frac{\mu}{n}
		\Big\|\int_\Omega \nabla\overline{\vartheta}_n({x},\cdot)\cdot\nabla\eta({x})\,d{x}\Big\|_{L^{2}(0,1)}
		\\
		&\quad
		+\frac1n
		\Big\|\int_\Omega
		(\overline{\mathbf{v}}_n({x},\cdot)\cdot\nabla)\eta({x})
		\vartheta_n({x},\cdot)\,d{x}\Big\|_{L^{2}(0,1)}\leq C,
	\end{align*}
	where $C$ does not depend on $n$. Similarly to the Lemmas \ref{Lemma 3.5} and \ref{Lemma 3.6}, we can easily estimate $\partial_{\overline{t}}\overline{\mathbf{v}}_n$ and $\partial_{\overline{t}}\overline{\vartheta}_n$. 
	\null\nobreak\hfill\ensuremath{\square}
	
	In order to prove the relative compactness of $\{\overline{\mathbf{v}}_n\}_{n\geq n_0}$ and $\{\overline{\vartheta}_n\}_{n\geq n_0}$ in $\mathrm{C}([0,1];\mathbf{H})$ and 
	$\mathrm{C}([0,1];L^2(\Omega))$, respectively, and apply the Aubin--Lions--Simon lemma; see Lemma \ref{ALSlemma}, we prove equi-continuity:
	\begin{align*}
		&\int_\Omega|\overline{\mathbf{v}}_n(x,\overline{t}+h)-\overline{\mathbf{v}}_n(x,\overline{t})
		|^2\,dx=
		\int_{\overline{t}}^{\overline{t}+h}
		\langle \partial_s\overline{\mathbf{v}}_n(\cdot,s),\overline{\mathbf{v}}_n(\cdot,\overline{t}+h)-\overline{\mathbf{v}}_n(\cdot,\overline{t})\rangle\,d s
		\\
		&\leq 
		2  \, \underset{{\overline{t}}\in(0,1)} {\rm esssup}
		\left\Vert\overline{\mathbf{v}}_n(\cdot,\overline{t})\right\Vert_{\mathbf{V}}  \int_{\overline{t}}^{\overline{t}+h}
		\left\|\partial_s \overline{\mathbf{v}}_n(\cdot,s)
		\right\|_{\mathbf{V}^\prime} \,d s
		\\
		&\leq 
		2 \, \underset{{\overline{t}}\in(0,1)} {\rm esssup}
		\left\Vert\overline{\mathbf{v}}_n(\cdot,\overline{t})\right\Vert_{\mathbf{V}} \, h^{\frac1{2}} \, \left(\int_0^1
		\left\|\partial_s \overline{\mathbf{v}}_n(\cdot,s)
		\right\|_{\mathbf{V}^\prime}^2\,d s \right)^{\frac{1}{2}}
		\leq C h^{\frac1{2}},
	\end{align*}
	and 
	\begin{align*}
		&\int_\Omega|\overline{\vartheta}_n(x,\overline{t}+h)-\overline{\vartheta}_n(x,\overline{t})
		|^2\,dx=
		\int_{\overline{t}}^{\overline{t}+h}
		\langle \partial_s \overline{\vartheta}_n(\cdot,s),
		\overline{\vartheta}_n(\cdot,\overline{t}+h)-\overline{\vartheta}_n(\cdot,\overline{t})\rangle\,d s
		\\
		&\leq 
		2 \, \underset{{\overline{t}}\in(0,1)} {\rm esssup}
		\left\Vert\overline{\vartheta}_n(\cdot,\overline{t})\right\Vert_{W_0^{1,2}(\Omega)}\, \int_{\overline{t}}^{\overline{t}+h}
		\left\|\partial_s \overline{\vartheta}_n(\cdot,s)
		\right\|_{W^{-1,2}(\Omega)}
		\,d s
		\\
		&\leq 2 \, \underset{{\overline{t}}\in(0,1)} {\rm esssup}
		\left\Vert\overline{\vartheta}_n(\cdot,\overline{t})\right\Vert_{W_0^{1,2}(\Omega)}\, h^{\frac1{2}} \, \left(\int_0^1
		\left\|\partial_s \overline{\vartheta}_n(\cdot,s)
		\right\|_{W^{-1,2}(\Omega)}^2\,d s\right)^{\frac{1}{2}}
		\leq C h^{\frac1{2}}.
	\end{align*}
	To finish the proof of item 1 of Theorem \ref{theorem3.1}, we need to prove \eqref{eq3.3a} and \eqref{eq3.3b}, which follow from the representation
	\begin{align*}
		{\overline{\varphi}_n({\overline{t}})-
			\overline{f}'({\overline{t}})}&=
		\frac{\nu}{n}\int_\Omega  \nabla\overline{\mathbf{v}}_n({x},{\overline{t}}):\nabla\boldsymbol{\omega}({x})\,d{x}
		-\frac1n\int_\Omega
		(\overline{\mathbf{v}}_n({x},{\overline{t}})\cdot\nabla)\boldsymbol{\omega}({x})\cdot\overline{\mathbf{v}}_n({x},{\overline{t}})\,d{x} 
		\\
		&\quad
		-\frac1n\int_\Omega
		\overline{\vartheta}_n({x},{\overline{t}})
		\overline{\mathbf{g}}_n({x},{\overline{t}})\cdot\boldsymbol{\omega}({x})
		\,d{x}+\overline{f}_n'({\overline{t}})-\overline{f}'({\overline{t}}),
		\\
		{\overline{\psi}_n({\overline{t}})-
			\overline{g}'({\overline{t}})}&=
		\frac{\mu}{n}\int_ \Omega\nabla\overline{\vartheta}_n({x},{\overline{t}})\cdot\nabla\eta({x})\,d{x}
		-\frac1n\int_\Omega
		(\overline{\mathbf{v}}_n({x},{\overline{t}})\cdot\nabla)\eta({x})\overline{\vartheta}_n({x},{\overline{t}})\,d{x} 
		\\
		&\quad
		+\overline{g}_n'({\overline{t}})-\overline{g}'({\overline{t}}).
	\end{align*}
	Applying the H\"older inequality, we need to prove the convergence of $\{\overline{\varphi}_n\}_{n\geq n_0}$ and $\{\overline{\psi}_n\}_{n\geq n_0}$ only in $L^2(0,1)$. We deduce that
	\begin{align*}
		&\left(\int_0^1|\overline{\varphi}_n({\overline{t}})-\overline{f}'({\overline{t}})|^2\,d{\overline{t}}\right)^{\frac12}
		\leq
		\frac{\nu}{\sqrt{n}}\left(\frac1n\int_{Q_1}
		|\nabla\overline{\mathbf{v}}_n({x},{\overline{t}})|^2\,d{x}d{\overline{t}}\right)^{\frac12}
		\left(\int_\Omega|\nabla\boldsymbol\omega (x)|^2\,d{x}\right)^{\frac12}
		\nonumber
		\\
		&
		\quad+\frac1n \|\nabla\boldsymbol\omega \|_{\widetilde{\mathbf{L}}^\infty} \, \underset{{\overline{t}}\in(0,1)} {\rm esssup}
		\int_\Omega|\overline{\mathbf{v}}_n({x},{\overline{t}})|^2\,d{x} 
		+
		\frac1n g_0\|\boldsymbol\omega\|_{\mathbf{H}}
		\, \underset{{\overline{t}}\in(0,1)}{\rm esssup} \, 
		\|\overline{\vartheta}_n(\cdot,{\overline{t}})\|_{L^2(\Omega)} 
		\\
		&\quad
		+\left(\int_0^1|\overline{f}_n'({\overline{t}})
		-\overline{f}'({\overline{t}})|^2\,d{\overline{t}}\right)^{\frac12},
		\\
		&\left(\int_0^1|\overline{\psi}_n({\overline{t}})-\overline{g}'({\overline{t}})|^2\,d{\overline{t}}\right)^{\frac12}
		\leq
		\frac{\mu}{\sqrt{n}}\left(\frac1n\int_{Q_1}
		|\nabla\overline{\vartheta}_n({x},{\overline{t}})|^2\,d{x}d{\overline{t}}\right)^{\frac{1}2}
		\left(\int_\Omega|\nabla\eta(x)|^2\,d{x}\right)^{\frac12}
		\nonumber
		\\
		&
		\quad+\frac1n \|\nabla\eta \|_{\mathbf{L}^\infty(\Omega)} \, \underset{{\overline{t}}\in(0,1)} {\rm esssup}
		\|\overline{\mathbf{v}}_n(\cdot,{\overline{t}})\|_{\mathbf{H}}  \, 
		\underset{{\overline{t}}\in(0,1)} {\rm esssup}
		\|\overline{\vartheta}_n(\cdot,{\overline{t}})\|_{L^2(\Omega)} 
		+\left(\int_0^1|\overline{g}_n'({\overline{t}})
		-\overline{g}'({\overline{t}})|^2\,d{\overline{t}}\right)^{\frac12}.
	\end{align*}
	From the above bounds and the condition ${\bf C}_{6}$ in Conditions \ref{C}, we infer that
	\begin{align}
		\left(\int_0^1|\overline{\varphi}_n({\overline{t}})-\overline{f}'({\overline{t}})|^2\,d{\overline{t}}\right)^{\frac12}
		\longrightarrow 0, \ \  
		\left(\int_0^1|\overline{\psi}_n({\overline{t}})-\overline{g}'({\overline{t}})|^2\,d{\overline{t}}\right)^{\frac12}
		\longrightarrow 0, \ \ \text{as} \ n \rightarrow \infty.
		\label{Equicon}
	\end{align}
	
	Therefore, there exist subsequences of $\{\overline{\mathbf{v}}_n\}_{n\geq n_0}$ and $\{\overline{\vartheta}_n\}_{n\geq n_0}$, still labeled by $n$, which converge weakly to some functions. We need to show that these functions are $\overline{\mathbf{v}}({x},{\overline{t}})$ and $\overline{\vartheta}({x},{\overline{t}})$ which satisfy problem \eqref{eq3.4}.
	
	We put in \eqref{eq5.7} the test functions $\widetilde{\boldsymbol{\phi}}\equiv0$ and 
	$\widetilde{\zeta}\equiv0$ in the domain $Q_T$ and take
	$\overline{\boldsymbol{\phi}}\theta_k(1-\overline{t})$ and 
	$\overline{\zeta}\theta_k(1-\overline{t})$, where $\theta_k(s)=sk\mathbf{1}_{(0,1/k]}(s)+
	\mathbf{1}_{(1/k,1]}(s)$ approximates the Heaviside function. Here, we assume that 
	$\overline{\boldsymbol{\phi}}$ and 
	$\overline{\zeta}$ have non-zero traces at $\overline{t}=1$. These rearrangements of the test functions lead to 
	\begin{align*}
		& \int_{Q_1} \theta_k(1-\overline{t})\Big(-\overline{\mathbf{v}}_n({x},{\overline{t}}) \cdot \partial_{\overline{t}} \overline{\boldsymbol{\phi}}({x},{\overline{t}})
		-\frac{1}{n} (\overline{\mathbf{v}}_n({x},{\overline{t}}) \cdot
		\nabla)\overline{\boldsymbol{\phi}}({x},{\overline{t}})\cdot
		\overline{\mathbf{v}}_n({x},{\overline{t}})
		\\
	& \quad	+ \frac{\nu}{ n}
		\nabla\overline{\mathbf{v}}_n({x},{\overline{t}}):
		\nabla \overline{\boldsymbol{\phi}}({x},{\overline{t}})
		-\overline{\varphi}_n({\overline{t}})\boldsymbol{\omega}({x})\cdot\overline{\boldsymbol{\phi}}({x},{\overline{t}}) 
		-\frac1n\overline{\vartheta}_n({x},{\overline{t}})
		\overline{\mathbf{g}}_n({x},{\overline{t}})\cdot\overline{\boldsymbol{\phi}}({x},{\overline{t}}) 
		\Big)\,d{x}d{\overline{t}}
		\nonumber\\&\quad-\int_\Omega \mathbf{v}_0({x}) \cdot \overline{\boldsymbol{\phi}}({x},0)\,d{x}+k\int_{1-\frac1k}^1\int_\Omega \overline{\mathbf{v}}_n({x},\overline{t}) \cdot \overline{\boldsymbol{\phi}}({x},\overline{t})\,d{x}d\overline{t}=0,
		\\
		&
		\int_{Q_1} \theta_k(1-\overline{t})\Big(-\overline{\vartheta}_n({x},{\overline{t}}) 
		\partial_{\overline{t}} \overline{\zeta}({x},{\overline{t}})
		-\frac{1}{n} (\overline{\mathbf{v}}_n({x},{\overline{t}})\cdot \nabla) \overline{\zeta}({x},{\overline{t}})\overline{\vartheta}_n({x},{\overline{t}})
		+\frac{\mu}{n}
		\nabla\overline{\vartheta}_n({x},{\overline{t}})\cdot
		\nabla \overline{\zeta}({x},{\overline{t}})
		\\
		& \quad
		-\overline{\varphi}_n({\overline{t}})\eta({x})\overline{\zeta}({x},{\overline{t}})
		\Big)\,d{x}d{\overline{t}}
		-\int_\Omega \vartheta_0({x})\overline{\zeta}({x},0)\,d{x} +k\int_{1-\frac1k}^1\int_\Omega \overline{\vartheta}_n({x},\overline{t})\overline{\zeta}({x},\overline{t})\,d{x}d\overline{t}=0,
	\end{align*}
	with subsequent passing to the limit as $k\to\infty$: 
	\begin{align*}
		& \int_{Q_1} \Big(-\overline{\mathbf{v}}_n({x},{\overline{t}}) \cdot \partial_{\overline{t}} \overline{\boldsymbol{\phi}}({x},{\overline{t}})
		-\frac{1}{n} (\overline{\mathbf{v}}_n({x},{\overline{t}}) \cdot
		\nabla)\overline{\boldsymbol{\phi}}({x},{\overline{t}})\cdot
		\overline{\mathbf{v}}_n({x},{\overline{t}})
		+ \frac{\nu}{ n}
		\nabla\overline{\mathbf{v}}_n({x},{\overline{t}}):
		\nabla \overline{\boldsymbol{\phi}}({x},{\overline{t}})
		\\
		& \quad
		-\overline{\varphi}_n({\overline{t}})\boldsymbol{\omega}({x})\cdot\overline{\boldsymbol{\phi}}({x},{\overline{t}}) 
		-\frac1n\overline{\vartheta}_n({x},{\overline{t}})
		\overline{\mathbf{g}}_n({x},{\overline{t}})\cdot\overline{\boldsymbol{\phi}}({x},{\overline{t}}) 
		\Big)\,d{x}d{\overline{t}}
		\nonumber\\&\quad-\int_\Omega \mathbf{v}_0({x}) \cdot \overline{\boldsymbol{\phi}}({x},0)\,d{x} {+\int_\Omega \overline{\mathbf{v}}_n({x},1) \cdot \overline{\boldsymbol{\phi}}({x},1)\,d{x}}=0,
		\\
		&
		\int_{Q_1} \Big(-\overline{\vartheta}_n({x},{\overline{t}}) 
		\partial_{\overline{t}} \overline{\zeta}({x},{\overline{t}})
		-\frac{1}{n} (\overline{\mathbf{v}}_n({x},{\overline{t}})\cdot \nabla) \overline{\zeta}({x},{\overline{t}})\overline{\vartheta}_n({x},{\overline{t}})
		+\frac{\mu}{n}
		\nabla\overline{\vartheta}_n({x},{\overline{t}})\cdot
		\nabla \overline{\zeta}({x},{\overline{t}})
		\\
		& \quad
		-\overline{\varphi}_n({\overline{t}})\eta({x})\overline{\zeta}({x},{\overline{t}})
		\Big)\,d{x}d{\overline{t}}
		-\int_\Omega \vartheta_0({x})\overline{\zeta}({x},0)\,d{x} {+\int_\Omega \overline{\vartheta}_n({x},1)\overline{\zeta}({x},1)\,d{x}}=0.
	\end{align*}
	Let us accomplish the passage to the limit in the problem \eqref{eq3.1}. Proposition \ref{proposition.2.1} involves 
	\begin{align*}
		&\left|
		\frac{\nu}{n}\int_{Q_1}
		\nabla\overline{\mathbf{v}}_n:\nabla\overline{\boldsymbol{\phi}}\,
		d{x}d{\overline{t}}\right|
		\leq \frac{\nu}{\sqrt{n}}
		\left( \frac{1}{n}\int_{Q_1}|\nabla\overline{\mathbf{v}}_n|^{2}\,d{x}d{\overline{t}}\right)^{1/2}
		{\|\overline{\boldsymbol{\phi}}\|_{L^2(0,1;\mathbf{V})}}\longrightarrow 0,
		\\
		&\left|\frac{1}{n}\int_{Q_1}
		(\overline{\mathbf{v}}_n\cdot\nabla)\overline{\boldsymbol{\phi}}\cdot \overline{\mathbf{v}}_n\,d{x}d{\overline{t}} \right|\leq
		\frac{C}{\sqrt{n}} \,
		\underset{{\overline{t}}\in(0,1)}{\mathrm{esssup}} \, \|\overline{\mathbf{v}}_n(\cdot,{\overline{t}})\|_{\mathbf{H}} \,
		\left( \frac{1}{n}\int_{Q_1}|\nabla\overline{\mathbf{v}}_n|^{2}\,d{x}d{\overline{t}}\right)^{1/2}
		\\
		& \qquad \qquad \qquad \qquad \qquad \qquad \times{\|\overline{\boldsymbol{\phi}}\|_{L^2(0,1;\mathbf{V})}}\longrightarrow0,
		\\
		&\left|\frac{\mu}{n}\int_{Q_1}
		\nabla\overline{\vartheta}_n \cdot \nabla\overline{\zeta}\,
		d{x}d{\overline{t}}\right|
		\leq \frac{\mu}{\sqrt{n}}
		\left( \frac{1}{n}\int_{Q_1}|\nabla\overline{\vartheta}_n|^{2}\,d{x}d{\overline{t}}\right)^{1/2}
		{\|\overline{\zeta}\|_{L^2(0,1;W_0^{1,2}(\Omega))}}\longrightarrow 0,
		\\
		&\left|\frac{1}{n}\int_{Q_1}
		(\overline{\mathbf{v}}_n\cdot\nabla)\overline{\zeta}\cdot \overline{\vartheta}_n\,d{x}d{\overline{t}} \right|\leq
		\frac{C}{\sqrt{n}} \, 
		\underset{{\overline{t}}\in(0,1)}{\mathrm{esssup}} \, \|\overline{\mathbf{v}}_n(\cdot,{\overline{t}})\|_{\mathbf{H}}^{\frac{1}{2}}
		\, \left( \frac{1}{n}\int_{Q_1}|\nabla\overline{\mathbf{v}}_n|^{2}\,d{x}d{\overline{t}}\right)^{\frac{1}{4}} 
		\\
		&\quad \qquad \times\underset{{\overline{t}}\in(0,1)}{\mathrm{esssup}} \, \|\overline{\vartheta}_n(\cdot,{\overline{t}})\|_{L^2(\Omega)}^{\frac{1}{2}}
		\, \left( \frac{1}{n}\int_{Q_1}|\nabla\overline{\vartheta}_n|^{2}\,d{x}d{\overline{t}}\right)^{1/4} {\|\overline{\zeta}\|_{L^2(0,1;W_0^{1,2}(\Omega))}}\longrightarrow0,
	\end{align*}
	as $n\to \infty$, {for all $\overline{\boldsymbol{\phi}}\in L^2(0,1;\mathbf{V})$ and
		$\overline{\zeta}\in L^2(0,1;W_0^{1,2}(\Omega))$.} Therefore, in the limit, we get
	\begin{align*}
		& \int_{Q_1} \Big(-\overline{\mathbf{v}}({x},{\overline{t}}) \cdot \partial_{\overline{t}} \overline{\boldsymbol{\phi}}({x},{\overline{t}})
		-\overline{\varphi}(\overline{t})\boldsymbol{\omega}({x})\cdot\overline{\boldsymbol{\phi}}({x},{\overline{t}}) \Big)\,d{x}d{\overline{t}}
		-\int_\Omega \mathbf{v}_0({x}) \cdot \overline{\boldsymbol{\phi}}({x},0)\, d{x} \\& \quad{+\int_\Omega \overline{\mathbf{v}}({x},1) \cdot \overline{\boldsymbol{\phi}}({x},1)\,d{x}}=0,
		\\
		& \int_{Q_1} \Big(-\overline{\vartheta}({x},{\overline{t}})\partial_{\overline{t}} \overline{\zeta}({x},{\overline{t}})
		-\overline{\psi}(\overline{t})\eta({x})\overline{\zeta}({x},{\overline{t}}) \Big)\,d{x}d{\overline{t}}
		-\int_\Omega \vartheta_0({x})\overline{\zeta}({x},0)\, d{x} \\& \quad {+\int_\Omega \overline{\vartheta}({x},1)\overline{\zeta}({x},1)\,d{x}}=0,
	\end{align*}
	which is a weak formulation of the problem \eqref{eq3.4} with $\overline{\varphi}(\overline{t})=\overline{f}^\prime(\overline{t})$ and $\overline{\psi}(\overline{t})=\overline{g}^\prime(\overline{t})$.
	\null\nobreak\hfill\ensuremath{\square}
	
	\subsection{\bf Proof of item 2 in Theorem \ref{theorem3.1}}
	\label{sec3.2}
	{In the present subsection, we formulate the definition of weak solutions to problem \eqref{eq3.2}, since the shifted initial conditions are less regular, i.e., $\widetilde{\mathbf{v}}_n(x,0)=\overline{\mathbf{v}}_n(x,1)\in\mathbf{H}$ and $\widetilde{\vartheta}_n(x,0)=\overline{\vartheta}_n(x,1)\in L^2(\Omega)$.}
	\begin{definition}
		\label{def3.2w} A quadruple of functions  $(\widetilde{\mathbf{v}}_n,\widetilde{\varphi}_n,\widetilde{\vartheta}_n,\widetilde{\psi}_n)$ is said to be a \textbf{weak solution} of the nonlinear inverse problem
		\eqref{eq3.2} with initial conditions $\overline{\mathbf{v}}_n({x},1) \in \mathbf{H}$ and $\overline{\vartheta}_n({x},1) \in L^2(\Omega)$, which are valid in the weak sense, if
		\begin{align*}
			&\widetilde{\mathbf{v}}_n\in L^\infty(0,T;\mathbf{H})\cap L^2(0,T;\mathbf{V}),
			\quad
			\partial_{\widetilde{t}}\widetilde{\mathbf{v}}_n\in 
			L^{2}(0,T;\mathbf{V}^\prime), \quad \widetilde{\varphi}_n\in L^{2}(0,T),
			\\
			&\widetilde{\vartheta}_n\in L^\infty(0,T;L^2(\Omega))\cap 
			L^2(0,T;W_0^{1,2}(\Omega)),
			\ \ \partial_{\widetilde{t}}\widetilde{\vartheta}_n\in 
			L^2(0,T;W^{-1,2}(\Omega)),\ \ \widetilde{\psi}_n\in L^2(0,T),
		\end{align*}
		and $(\widetilde{\mathbf{v}}_n,
		\widetilde{\vartheta}_n,
		\widetilde{\varphi}_n,
		\widetilde{\psi}_n)$ satisfies 
		\begin{align*}
			& \frac{d}{d\widetilde{t}}\int_\Omega
			\widetilde{\mathbf{v}}_n({x},\widetilde{t})
			\cdot\boldsymbol{\phi}({x})\,d{x}
			+\nu\int_\Omega\nabla\widetilde{\mathbf{v}}_n({x},\widetilde{t})
			:\nabla\boldsymbol{\phi}({x})\,d{x} 
			-\int_\Omega
			(\widetilde{\mathbf{v}}_n({x},\widetilde{t})
			\cdot\nabla)\boldsymbol{\phi}({x})\cdot\widetilde{\mathbf{v}}_n({x},\widetilde{t})
			\,d{x}
			\\
			&
			=\frac1n\int_\Omega\widetilde{\vartheta}_n({x},\widetilde{t})
			\widetilde{\mathbf{g}}_n({x},\widetilde{t})
			\cdot\boldsymbol{\phi}({x})\,dx +\widetilde{\varphi}_n(\widetilde{t})\int_\Omega 
			\boldsymbol{\omega}({x})\cdot\boldsymbol{\phi}({x})\,d{x} 
			\,\,\, \text{ for a.e. }\widetilde{t}\in(0,T),
			\\
			& \frac{d}{d\widetilde{t}}\int_\Omega
			\widetilde{\vartheta}_n({x},\widetilde{t})
			\zeta({x})\,d{x}
			+\mu\int_\Omega\nabla\widetilde{\vartheta}_n({x},\widetilde{t})
			\cdot\nabla\zeta({x})\,d{x} 
			-\int_\Omega
			(\widetilde{\mathbf{v}}_n({x},\widetilde{t})
			\cdot\nabla)\zeta({x})\widetilde{\vartheta}_n({x},\widetilde{t})\,d{x}
			\\
			&
			=\widetilde{\psi}_n(\widetilde{t})
			\int_\Omega\eta({x})\zeta({x})\,d{x} \, \, \, \text{ for a.e. }\widetilde{t}\in(0,T),
		\end{align*}
		for every test functions $\boldsymbol{\phi}\in\mathbf{V}$ and $\zeta \in W_0^{1,2}(\Omega)$, and the corresponding compatibility conditions
		\begin{align*}
			&\int_\Omega
			\overline{\mathbf{v}}_n({x},1)\cdot
			\boldsymbol{\omega}({x})\,d{x}=\widetilde{f}_n(0),
			\quad 
			\int_\Omega\overline{\vartheta}_n({x},1)\eta({x})\,d{x}=\widetilde{g}_n(0).
		\end{align*}
	\end{definition}
	
	\begin{proposition}[\textbf{Estimate for a shifted solution}]
		\label{proposition.2.2} Under Conditions \ref{C}, the problem \eqref{eq3.2} has {\bf a weak solution} in the sense of Definition \ref{def3.2w}, which satisfies the estimate
		\begin{align}
			&\underset{\widetilde{t}\in(0,T)}{\rm esssup} \|\widetilde{\mathbf{v}}_n(\cdot,\widetilde{t})
			\|_{\mathbf{H}}+
			\|\widetilde{\mathbf{v}}_n\|_{L^2(0,T;\mathbf{V})}
			+\|\widetilde{\varphi}_n\|_{L^2(0,T)}
			+\|\partial_{\widetilde{t}}\widetilde{\mathbf{v}}_n\|_{L^2(0,T;\mathbf{V}^\prime)}
			\leq C,
			\label{eq:est-global-11a}
			\\
			&\underset{\widetilde{t}\in(0,T)}{\rm esssup} \|\widetilde{\vartheta}_n(\cdot,\widetilde{t})
			\|_{L^{2}(\Omega)}+
			\|\widetilde{\vartheta}_n\|_{L^2(0,T;W_0^{1,2}(\Omega))} 
			+\|\widetilde{\psi}_n\|_{L^2(0,T)}
			\nonumber
			\\
			&\quad
			+
			\|\partial_{\widetilde{t}}\widetilde{\vartheta}_n\|_{L^{2}(0,T;W^{-1,2}(\Omega))}
			\leq C,
			\label{eq:est-global-11b}
		\end{align}
		where the constant $C$ does not depend on $n$. 
	\end{proposition}
	\noindent{\bf Proof.}
	{ Similar to Lemmas \ref{1Globest:1}, \ref{Lemma 3.5}, and \ref{Lemma 3.6} (see also the proof of Proposition \ref{proposition.2.1}), we can obtain the following estimates:
		\begin{align*}
			&\|\widetilde{\mathbf{v}}_n\|
			_{L^\infty(0,T;\mathbf{H})}+
			\|\widetilde{\mathbf{v}}_n\|_{L^2(0,T;\mathbf{V})}
			+\|\partial_{\widetilde{t}}\widetilde{\mathbf{v}}_n\|_{L^2(0,T;\mathbf{V}^\prime)}+\|\widetilde{\varphi}_n\|_{L^2(0,T)}\leq C,
			\\
			&\|\widetilde{\vartheta}_n\|_{L^\infty(0,T;L^{2}(\Omega))}+
			\|\widetilde{\vartheta}_n\|_{L^2(0,T;W_0^{1,2}(\Omega))}
			+\|\partial_{\widetilde{t}}\widetilde{\vartheta}_n\|_{L^2(0,T;W^{-1,2}(\Omega))}+\|\widetilde{\psi}_n\|_{L^2(0,T)}\leq C,
		\end{align*}
		where the constant $C$ does not depend on $n$.}
	\null\nobreak\hfill\ensuremath{\square}
	
	\begin{corollary}
		\label{cor4.1}
		The families $\{\widetilde{\mathbf{v}}_n\}_{n\geq n_0}$ and  $\{\widetilde{\vartheta}_n\}_{n\geq n_0}$ are  integral equi-continuous in time in $L^2(0,T;\mathbf{H})$ and  $L^2(0,T;L^2(\Omega))$, respectively.
	\end{corollary}
	
	\noindent{\bf Proof.}
	The integral equi-continuity in $\widetilde{t}$ for $\{\widetilde{\mathbf{v}}_n\}_{n\geq n_0}$ and for $\{\widetilde{\vartheta}_n\}_{n\geq n_0}$ follows from Corollary \ref{cor3.1}:
	\begin{align*}
		&\int_0^{T-h}\int_\Omega|\widetilde{\mathbf{v}}_n({x},
		\widetilde{t}+h)-
		\widetilde{\mathbf{v}}_n({x},\widetilde{t})|^2\,d{x}d\widetilde{t}
		\\
		&\quad=\int_{\frac1n}^{T_n-h}\int_\Omega
		|\mathbf{v}_n({x},t+h)-
		\mathbf{v}_n({x},t)|^2\,d{x}dt\le C\omega(h),
		\\
		&\int_0^{T-h}\int_\Omega
		|\widetilde{\vartheta}_n({x},
		\widetilde{t}+h)-
		\widetilde{\vartheta}_n({x},\widetilde{t})|^2\,d{x}d\widetilde{t}
		\\
		&\quad=\int_{\frac1n}^{T_n-h}\int_\Omega
		|\vartheta_n({x},t+h)-\vartheta_n({x},t)|^2\,d{x}dt\le C\omega(h),
	\end{align*}
	where $C$ does not depend on $n$ and $h$, and $\omega(h)$ is defined in Corollary \ref{cor3.1}.
	\null\nobreak\hfill\ensuremath{\square}
	
	\begin{corollary}
		The family $\{\widetilde{\mathbf{v}}_n\}_{n\geq n_0}$ coincides in the limit with  $\{\mathbf{v}_n\}_{n\geq n_0}$ in $L^2(0,T;\mathbf{H})$, and the family $\{\widetilde{\vartheta}_n\}_{n\geq n_0}$ coincides in the limit with  $\{\vartheta_n\}_{n\geq n_0}$ in $L^2(0,T;L^{2}(\Omega))$.
	\end{corollary}
	\noindent{\bf Proof.}
	This follows from Corollaries \ref{cor3.1} and \ref{cor4.1}.
	\null\nobreak\hfill\ensuremath{\square}
	
	\begin{corollary}
		\label{cor4.3}
		The family $\{\widetilde{\mathbf{v}}_n\}_{n\geq n_0}$ is relatively compact in 
		$\mathrm{C}([0,T];\mathbf{V}^\prime)$, and the family $\{\widetilde{\vartheta}_n\}_{n\geq n_0}$ is relatively compact in $\mathrm{C}([0,T];W^{-1,2}(\Omega))$.
	\end{corollary}
	\noindent{\bf Proof.} According to the Aubin--Lions--Simon lemma, see Lemma \ref{ALSlemma}, that 
	$\{\widetilde{\mathbf{v}}_n\}_{n\geq n_0}$ is relatively compact in $\mathrm{C}([0,T];\mathbf{V}^\prime)$, and 
	$\{\widetilde{\vartheta}_n\}_{n\geq n_0}$ is relatively compact in $\mathrm{C}([0,T];W^{-1,2}(\Omega))$.
	\null\nobreak\hfill\ensuremath{\square}
	\vspace{3mm}
	
	{We need to pass to the limit as $n\to\infty$. Before, we put in \eqref{eq5.7} in the domain ${Q}_1$ the test functions $(\overline{\boldsymbol{\phi}},\overline{\zeta})\equiv0$  and rearrange test functions $(\widetilde{\boldsymbol{\phi}},\widetilde{\zeta})$ on the plane $\widetilde{t}=0$ as the product $(\widetilde{\boldsymbol{\phi}},\widetilde{\zeta})$ and the approximation of the Heaviside function $\theta_k(\widetilde{t})=k\widetilde{t}\mathbf{1}_{[0,1/k]}(\widetilde{t})+\mathbf{1}_{(1/k,T)}(t)$:
		\begin{align*}
			&\int_{Q_T}\theta_k(\widetilde{t})\Big(-\widetilde{\mathbf{v}}_n ({x},\widetilde{t}) \cdot \partial_{\widetilde{t}}\, \widetilde{\boldsymbol{\phi}}({x},\widetilde{t})
			- (\widetilde{\mathbf{v}}_n({x},\widetilde{t})\cdot\nabla)
			\widetilde{\boldsymbol{\phi}}
			\cdot \widetilde{\mathbf{v}}_n({x},\widetilde{t}) 
			+\nu\nabla\widetilde{\mathbf{v}}_n ({x},\widetilde{t}) : \nabla \widetilde{\boldsymbol{\phi}}({x},\widetilde{t})
			\\[0.8ex]
			&\quad-\widetilde{\vartheta}_n({x},\widetilde{t})\widetilde{\mathbf{g}}_n({x},\widetilde{t})\cdot\widetilde{\boldsymbol{\phi}}({x},\widetilde{t})
			-\widetilde{\varphi}_n(\widetilde{t})\boldsymbol{\omega}({x})\cdot\widetilde{\boldsymbol{\phi}}({x},\widetilde{t})
			\Big)\, d{x} d\widetilde{t} 	
			\\[0.8ex]
			&\quad-k\int_0^{\frac1k}\int_\Omega\widetilde{\mathbf{v}}_n({x},\widetilde{t})
			\cdot\widetilde{\boldsymbol{\phi}}({x},\widetilde{t})\,d{x}d\widetilde{t}=0,
			\\
			&\int_{Q_T} \theta_k(\widetilde{t})\Big(-\widetilde{\vartheta}_n ({x},\widetilde{t}) \partial_{\widetilde{t}}\, \widetilde{\zeta}({x},\widetilde{t})
			- (\widetilde{\mathbf{v}}_n({x},\widetilde{t})\cdot\nabla)\widetilde{\zeta}({x},\widetilde{t})
			\widetilde{\vartheta}_n({x},\widetilde{t}) 
			+\mu\nabla\widetilde{\vartheta}_n ({x},\widetilde{t}) 
			\cdot \nabla \widetilde{\zeta}({x},\widetilde{t})
			\\[0.8ex]
			&\quad
			-\widetilde{\psi}_n(\widetilde{t})
			\eta({x})\widetilde{\zeta}({x},\widetilde{t})\Big)\, d{x} d\widetilde{t}
			-k\int_0^{\frac1k}\int_\Omega\widetilde{\vartheta}_n({x},\widetilde{t})\widetilde{\zeta}({x},\widetilde{t})\,d{x}d\widetilde{t}=0.
		\end{align*}
		After, we pass to the limit as $k\to\infty$, and encompass the traces at $\widetilde{t}=0$ by the rules \eqref{eq5.8}:
		\begin{align*}
			&\int_{Q_T} \Big(-\widetilde{\mathbf{v}}_n ({x},\widetilde{t}) \cdot \partial_{\widetilde{t}}\, \widetilde{\boldsymbol{\phi}}({x},\widetilde{t})
			- (\widetilde{\mathbf{v}}_n({x},\widetilde{t})\cdot\nabla)
			\widetilde{\boldsymbol{\phi}}
			\cdot \widetilde{\mathbf{v}}_n({x},\widetilde{t}) 
			+\nu\nabla\widetilde{\mathbf{v}}_n ({x},\widetilde{t}) : \nabla \widetilde{\boldsymbol{\phi}}({x},\widetilde{t})
			\\[0.8ex]
			&\quad-\widetilde{\vartheta}_n({x},\widetilde{t})\widetilde{\mathbf{g}}_n({x},\widetilde{t})\cdot\widetilde{\boldsymbol{\phi}}({x},\widetilde{t})
			-\widetilde{\varphi}_n(\widetilde{t})\boldsymbol{\omega}({x})\cdot\widetilde{\boldsymbol{\phi}}({x},\widetilde{t})
			\Big)\, d{x} d\widetilde{t}-\int_\Omega\overline{\mathbf{v}}_n({x},1)
			\cdot\widetilde{\boldsymbol{\phi}}({x},0)\,d{x}=0,
			\\
			&\int_{Q_T} \Big(-\widetilde{\vartheta}_n ({x},\widetilde{t}) \partial_{\widetilde{t}}\, \widetilde{\zeta}({x},\widetilde{t})
			- (\widetilde{\mathbf{v}}_n({x},\widetilde{t})\cdot\nabla)\widetilde{\zeta}({x},\widetilde{t})
			\widetilde{\vartheta}_n({x},\widetilde{t}) 
			+\mu\nabla\widetilde{\vartheta}_n ({x},\widetilde{t}) 
			\cdot \nabla \widetilde{\zeta}({x},\widetilde{t})
			\\[0.8ex]
			&\quad
			-\widetilde{\psi}_n(\widetilde{t})
			\eta({x})\widetilde{\zeta}({x},\widetilde{t})\Big)\, d{x} d\widetilde{t}
			-\int_\Omega\overline{\vartheta}_n({x},1)\widetilde{\zeta}({x},0)\,d{x}=0.
	\end{align*}}
	Using the energy estimates \eqref{eq:est-global-11a} and \eqref{eq:est-global-11b}, along with the Banach--Alaoglu theorem, we can extract subsequences of $\{\widetilde{\mathbf{v}}_n\}_{n \geq n_0}$ and $\{\widetilde{\vartheta}_n\}_{n \geq n_0}$, still labeled by $n$, such that the following convergences hold:
	\begin{subequations}
		\label{Weak conv}
		\begin{align}
			\widetilde{\mathbf{v}}_n
			\xrightharpoonup[n\to\infty]{}
			&\mathbf{v}  \ \ \text{ weakly in } L^2(0,T;\mathbf{V}),
			\\
			\widetilde{\mathbf{v}}_n \xrightharpoonup[n\to\infty]{*}
			&\mathbf{v}  \ \ \text{ weakly$^*$ in } L^\infty(0,T;\mathbf{H}),
			\\
			\widetilde{\vartheta}_n
			\xrightharpoonup[n\to\infty]{}
			&\vartheta \ \ \text{ weakly in } L^2(0,T;W_0^{1,2}(\Omega)),
			\\
			\widetilde{\vartheta}_n \xrightharpoonup[n\to\infty]{*}
			&\vartheta \ \ \text{ weakly$^*$ in } L^\infty(0,T;L^{2}(\Omega)).
		\end{align}
	\end{subequations}
	We obtain an integral identity in the limit:
	\begin{align*}
		&\int_{Q_T} \Big(-\mathbf{v}({x},\widetilde{t}) \cdot \partial_{\widetilde{t}}\, \widetilde{\boldsymbol{\phi}}({x},\widetilde{t})
		- (\mathbf{v}({x},\widetilde{t}) \cdot\nabla)\widetilde{\boldsymbol{\phi}}({x},\widetilde{t})\cdot\mathbf{v}({x},\widetilde{t}) 
		+\nu\nabla\mathbf{v}({x},\widetilde{t}):\nabla\widetilde{\boldsymbol{\phi}}({x},\widetilde{t})
		\\[0.8ex]
		&\quad
		-\vartheta({x},\widetilde{t})
		\mathbf{g}({x},\widetilde{t})
		\cdot\widetilde{\boldsymbol{\phi}}({x},\widetilde{t})
		-\varphi(\widetilde{t};[\mathbf{v},\vartheta])\boldsymbol{\omega}({x})
		\cdot\widetilde{\boldsymbol{\phi}}({x},\widetilde{t})
		\Big)\, d{x} d\widetilde{t} 
		-\int_\Omega\overline{\mathbf{v}}({x},1)
		\cdot\widetilde{\boldsymbol{\phi}}({x},0)\,d{x}
		=0,
		\\[0.8ex]
		&\int_{Q_T} \Big(-\vartheta({x},\widetilde{t})
		\partial_{\widetilde{t}}\, \widetilde{\zeta}({x},\widetilde{t})
		-
		(\mathbf{v}({x},\widetilde{t})\cdot\nabla)\widetilde{\zeta}({x},\widetilde{t})
		\vartheta({x},\widetilde{t})
		+\mu\nabla\vartheta({x},\widetilde{t})\cdot\nabla \widetilde{\zeta}({x},\widetilde{t})
		\\[0.8ex]
		&\quad
		-\psi(\widetilde{t};[\mathbf{v},\vartheta])\eta({x})\widetilde{\zeta}
		\Big)\,dxd\widetilde{t} 
		-\int_\Omega\overline{\vartheta}({x},1)
		\widetilde{\zeta}({x},0)\,d{x}=0,
	\end{align*}
	with
	\begin{align*}
		\varphi(\widetilde{t};[\mathbf{v},\vartheta])&=f'(\widetilde{t})
		+\nu\int_\Omega 
		\nabla\mathbf{v}(x,\widetilde{t}):\nabla\boldsymbol{\omega}({x})\,d{x}
		-\int_\Omega
		(\mathbf{v}({x},\widetilde{t})\cdot \nabla)
		\boldsymbol{\omega}({x})\cdot\mathbf{v}({x},\widetilde{t})\,d{x}
		\\
		&\quad  -\int_\Omega
		\vartheta({x},\widetilde{t})
		\mathbf{g}({x},\widetilde{t})\cdot\boldsymbol{\omega}({x})\,d{x},
		\\
		\psi(\widetilde{t};[\mathbf{v},\vartheta])&=g'(\widetilde{t})
		+\mu\int_\Omega\nabla\vartheta({x},\widetilde{t})\cdot\nabla\eta({x})\,d{x}
		-\int_\Omega
		(\mathbf{v}({x},\widetilde{t})\cdot\nabla)\eta(x)\widetilde{\vartheta}({x},\widetilde{t})\,d{x}.
	\end{align*}
	
	It is clear that the limit solutions $(\mathbf{v},\vartheta,\varphi,\psi)$, where $\mathbf{v}=\lim\limits_{n\to\infty}\widetilde{\mathbf{v}}_n$, 
	$\vartheta=\lim\limits_{n\to\infty}\widetilde{\vartheta}_n$, $\varphi=\lim\limits_{n\to\infty}\widetilde{\varphi}_n$, and  $\psi=\lim\limits_{n\to\infty}\widetilde{\psi}_n$, satisfy the problem \eqref{eq3.6}.
	
	\subsection{\bf Proof of item 3 in Theorem \ref{theorem3.1}}
	\label{sec3.3}
	In item 2, we have shown that the limit solution $\mathbf{v}(x,t)=\lim\limits_{n\to\infty}\widetilde{\mathbf{v}}_n({x},t)$ and $\vartheta(x,t)=\lim\limits_{n\to\infty}\widetilde{\vartheta}_n({x},t)$ satisfies the new initial data:
	\[
	\mathbf{v}({x},0)=\lim\limits_{n\to\infty}
	\widetilde{\mathbf{v}}_n({x},0)=
	\overline{\mathbf{v}}({x},1),\quad 
	\vartheta({x},0)=\lim\limits_{n\to\infty}
	\widetilde{\vartheta}_n({x},0)=
	\overline{\vartheta}({x},1),
	\]
	in a \textbf{strong sense}, see Corollary \ref{cor4.3}. 
	
	{ We have proven that the family $\{\partial_t\mathbf{v}_{n}\}_{n\geq n_0}$ is uniformly bounded in $L^1(0,T;\mathbf{V}^\prime)$, and the family $\{\partial_t\vartheta_{n}\}_{n\geq n_0}$ is uniformly bounded in $L^1(0,T;W^{-1,2}(\Omega))$. These facts, along with the uniform boundedness of  $\{\mathbf{v}_{n}\}_{n\geq n_0}$ in $L^2(0,T;\mathbf{V})$ and $\{\vartheta_{n}\}_{n\geq n_0}$ in $L^2(0,T;W_0^{1,2}(\Omega))$, is sufficient to apply the Aubin--Lions--Simon lemma, see Lemma \ref{ALSlemma}, and to prove the relative compactness of $\{\mathbf{v}_n\}_{n\geq n_0}$ in $L^2(0,T;\mathbf{H})$, and $\{\vartheta_n\}_{n\geq n_0}$ in $L^2(0,T;L^2(\Omega))$. There is a subsequence that strongly converges in $L^2(0,T;\mathbf{H})$ to $\mathbf{v}$, which is the limit of $\{\widetilde{\mathbf{v}}_n\}_{n\geq n_0}$, and there is a subsequence that strongly converges in $L^2(0,T;L^2(\Omega))$ to $\vartheta$, which is the limit of $\{\widetilde{\vartheta}_n\}_{n\geq n_0}$:
		\[
		\int_0^T\|\mathbf{v}_n(\cdot,\widetilde{t})-
		\widetilde{\mathbf{v}}_n(\cdot,\widetilde{t})\|_{\mathbf{H}}^2\,d\widetilde{t}\longrightarrow 0,
		\quad \int_0^T\|\vartheta_n(\cdot,\widetilde{t})-
		\widetilde{\vartheta}_n(\cdot,\widetilde{t})\|_{L^2(\Omega)}^2\,d\widetilde{t}\longrightarrow0
		\text{ as }n\longrightarrow\infty,
		\]
		due to Corollary \ref{cor3.1}. Moreover, $\mathbf{1}_{(1/n,T)}(t)\varphi_{n}$ and  $\mathbf{1}_{(1/n,T)}(t)\psi_{n}$ converge weakly in $L^2(0,T)$ to $\varphi$ and $\psi$, respectively.}
	
	\section*{Conclusion}
	The existence result, see Theorem \ref{thm1} for a fixed $n\in\mathbb{N}$, is based on several well-known facts about inverse problems for the Boussinesq system with the integral over-determination conditions. In fact, we deal with the family of inverse problems, where the right-hand sides $f_n$ and $g_n$ approximate discontinuous functions as $n\to\infty$. In order to pass to the limit, we have obtained several uniform estimates.
	The main novelty - Theorem \ref{theorem3.1} - is linked with passing to the limit as $n\to\infty$. We apply the scaling  $\overline{\mathbf{v}}_n(x,\overline{t})=
	\mathbf{v}_n({x},\overline{t}/n)$ and $\overline{\vartheta}_n({x},\overline{t})=\vartheta_n({x},
	\overline{t}/n)$ to obtain the auxiliary inverse problem on the infinitesimal initial layer. 
	In the limit, we obtain new initial data $\mathbf{v}({x},0)=\overline{\mathbf{v}}(x,1)$ and $\vartheta({x},0)=\overline{\vartheta}({x},1)$ defined by \eqref{eq3.5} for the limit problem \eqref{eq3.6}.

\subsection*{\bf Acknowledgments}
S.N. Antontsev was supported by the Russian Science Foundation project no. 25-11-00056. P. Kumar gratefully acknowledges Prof. Ankik Kumar Giri, Department of Mathematics, Indian Institute of Technology Roorkee, India for kindly hosting his visit and providing a conducive research environment during the visit. I. Kuznetsov and S. Aitzhanov were supported by the Ministry of Science and Higher Education of Kazakhstan under project no. AP26199323. I. Kuznetsov was supported by the Ministry of Science and Higher Education of the Russian Federation under project no. FWGG-2026-0001.

\section*{Data Availability Statement}
The authors state that no data were generated or analysed in this
paper, so data availability is not applicable.

\section*{Conflict of interest} The authors state that there is no
conflict of interest.

\bibliographystyle{ieeetr}
\bibliography{bibl}

\end{document}